\documentclass[11pt]{article} 
 
\usepackage{latexsym} 
\usepackage{amsfonts} 
\usepackage{amsmath} 
\usepackage{amsthm,stackrel,esint} 
\usepackage{amssymb}
\usepackage{mathrsfs} 
\usepackage{dsfont,yfonts} 
\usepackage{bbold} 
\usepackage[misc]{ifsym}
\usepackage{caption} 
\usepackage{epsfig} 
\usepackage{float} 
\usepackage{subfigure,mathtools} 
\usepackage{psfrag} 
\usepackage{graphicx} 
\usepackage{epsfig} 
\usepackage{amsbsy} 
\usepackage{hyperref,color}
\usepackage{xcolor}
\usepackage{enumitem}
\usepackage{quoting}
\usepackage{scalerel,stackengine}
\usepackage[numbers]{natbib}
\usepackage{setspace}

\stackMath
\newcommand\reallywidehat[1]{%
	\savestack{\tmpbox}{\stretchto{%
			\scaleto{%
				\scalerel*[\widthof{\ensuremath{#1}}]{\kern-.6pt\bigwedge\kern-.6pt}%
				{\rule[-\textheight/2]{1ex}{\textheight}}%WIDTH-LIMITED BIG WEDGE
			}{\textheight}% 
		}{0.5ex}}%
	\stackon[1pt]{#1}{\tmpbox}%
}
\DeclareMathAlphabet{\mathpgoth}{OT1}{pgoth}{m}{n}
\DeclareMathAlphabet{\mathpzc}{OT1}{pzc}{m}{it}

\DeclareMathOperator{\meas}{meas}

\newtheorem{theorem}{Theorem} 
\newtheorem*{prop*}{Theorem} 
 
\newtheorem{coro}[theorem]{Corollary} 
 
\newtheorem{lemma}[theorem]{Lemma} 
\newtheorem{prop}[theorem]{Proposition} 
\newtheorem{rmk}[theorem]{Remark}

\numberwithin{theorem}{section}

\newcommand{\zerarcounters}{\setcounter{equation}{0}\setcounter{theorem}{0}}

\newcommand{\ZZZ}{\mathds{Z}} 
\newcommand{\CCC}{\mathds{C}} 
\newcommand{\NNN}{\mathds{N}} 
 
\newcommand{\RRR}{\mathds{R}} 
\newcommand{\TTT}{\mathds{T}} 
\newcommand{\uno}{\mathds{1}} 
\newcommand{\gsa}{{\mathtt{g}(s,\al)}}

\newcommand{\uu}{\mathcal{U}} 
\newcommand{\UU}{\mathcal{U}} 
\newcommand{\uuu}{{U}} 
\newcommand{\UUU}{{U}} 
\newcommand{\aaa}{\ol{a}}
 
\newcommand{\CCCCC}{\mathfrak{C}} 
\newcommand{\EEE}{\mathcal{V}_{i_0}}
\newcommand{\norm}{{|\hspace{-0.025cm}|\hspace{-0.025cm}|}} 
\newcommand{\IIII}{{\calI}}
\newcommand{\IIIIO}{{\matC_0}}
\newcommand{\IIIIC}{{\matC}}
 
\newcommand{\calA}{{\mathcal A}} 
\newcommand{\BB}{{\mathcal B}} 
 
\newcommand{\DD}{{\mathcal D}} 
 
\newcommand{\calF}{{\mathcal F}}

\newcommand{\calI}{{\mathcal I}}

\newcommand{\MM}{{\mathcal M}}

\newcommand{\TT}{{\mathcal T}} 
\newcommand{\calU}{{\mathcal U}} 
\newcommand{\WW}{{\mathpzc{W}}} 
\newcommand{\VV}{{\mathpzc{V}}}

\newcommand{\tx}{\mathtt x}

\newcommand{\calmF}{{\mathscr F}}

\DeclareFontFamily{U}{BOONDOX-calo}{\skewchar\font=45 }
\DeclareFontShape{U}{BOONDOX-calo}{m}{n}{
  <-> s*[1.05] BOONDOX-r-calo}{}
\DeclareMathAlphabet{\mathcalboondox}{U}{BOONDOX-calo}{m}{n}
\DeclareMathAlphabet{\mathbcalboondox}{U}{BOONDOX-calo}{b}{n}
\DeclareMathAlphabet{\mathpzc}{OT1}{pzc}{m}{it} 
 
\newcommand{\gota}{{\mathfrak a}}

\newcommand{\gote}{{\mathfrak e}}

\newcommand{\gotr}{{\mathfrak r}}

\newcommand{\gotu}{{\mathfrak u}}

\newcommand{\gotB}{{\mathbcalboondox B}}

\newcommand{\gD}{{\mathbcalboondox D}}

\newcommand{\gotG}{{\mathfrak G}}

\newcommand{\gQ}{{\mathbcalboondox Q}} 

\newcommand{\gotU}{{\mathfrak U}}

\newcommand{\gotWW}{{\mathbcalboondox W}}

\newcommand{\matA}{{\mathscr A}} 
 
\newcommand{\matC}{{\mathscr C}} 
\newcommand{\matD}{{\mathscr D}} 
\newcommand{\matE}{{\mathscr E}} 
\newcommand{\matF}{{\mathscr F}} 
\newcommand{\calG}{{\mathscr G}} %\newcommand{\matG}{{\mathscr G}} 

\newcommand{\matK}{{\mathscr K}}

\newcommand{\matS}{{\mathscr S}} 
\newcommand{\matT}{{\mathscr T}} 
\newcommand{\matU}{{\mathscr U}} 
\newcommand{\matV}{{\mathscr V}}

\newcommand{\matW}{{\mathscr W}} 
\newcommand{\matY}{{\mathscr Y}} 
 
\newcommand{\und}{\underline}

\newcommand{\ol}{\overline} 
 
\newcommand{\prova}{\noindent\textit{Proof. }} 
\newcommand{\io}{\infty} 
\newcommand{\e}{\varepsilon} 
\newcommand{\al}{\alpha} 
\newcommand{\de}{\delta} 
\newcommand{\be}{\beta}

\newcommand{\x}{\xi} 
 
\newcommand{\ka}{\kappa} 
\newcommand{\g}{\gamma} 
\newcommand{\om}{\omega} 
\newcommand{\h}{\eta} 
\newcommand{\ze}{\zeta} 
\newcommand{\la}{\lambda} 
\newcommand{\f}{\varphi} 
\newcommand{\s}{\sigma} 
 
\newcommand{\del}{\partial}

\newcommand{\der}{{\rm d}} 
\newcommand{\ii}{{\rm i}}

\newcommand{\jap}[1]{\langle #1 \rangle}

\newcommand{\pow}{q}

\makeatletter
\newcommand{\oset}[3][0ex]{%
  \mathrel{\mathop{#3}\limits^{
    \vbox to#1{\kern-1\ex@
    \hbox{$\scriptstyle#2$}\vss}}}}
\makeatother

\def\tilde#1{\widetilde{#1}}
 
\def\ins#1#2#3{\vbox to0pt{\kern-#2 \hbox{\kern#1 #3}\vss}\nointerlineskip} 

\newcommand{\DgN}{\matK_N(\g)}

\begin{document}
	
%%%%%%%%%%%%%%%%%%%%%%%%%%%%%%%%%%%%%%%%%%%%%%%%%%%%%%%%%%%%%%%%%%%
\title{\bf Asymptotically full measure sets\\of almost-periodic solutions for the NLS equation}
%%%%%%%%%%%%%%%%%%%%%%%%%%%%%%%%%%%%%%%%%%%%%%%%%%%%%%%%%%%%%%%%%%%
	
\author{\textbf{
Luca Biasco, Livia Corsi, Guido Gentile, Michela Procesi}\\
\\
\small Dipartimento di Matematica e Fisica, Universit\`a Roma Tre, Roma, 00146, Italy\\
\footnotesize 
\Letter~luca.biasco@uniroma3.it, livia.corsi@uniroma3.it, guido.gentile@uniroma3.it, michela.procesi@uniroma3.it
}

\date{} 

\maketitle 

%%%%%%%%%%%%%%%%%%%%%%%%%%%%%%%%%%%%%%%%%%%%%%%%%%%%%%%%%%%%%%%%%%%%%%%%%% 
%%%%%%%%%%%%%%%%%%%%%%%%%%%%%%%%%%%%%%%%%%%%%%%%%%%%%%%%%%%%%%%%%%%%%%%%%% 
\begin{abstract} 
We study the  dynamics of solutions for  a family of nonlinear Schr\"odinger equations on the circle, 
with a smooth convolution potential and Gevrey regular initial data. Our main result is the construction of
an asymptotically full measure set of small-amplitude time almost-periodic solutions,
which are dense on invariant tori.
In regions corresponding to positive actions, we prove that such maximal invariant tori are
Banach manifolds, which provide a Cantor foliation of the phase space.  As a consequence,
we establish that, for many small initial data, the Gevrey norm of the solution remains approximately constant
for all time and hence the elliptic fixed point at the origin is Lyapunov statistically stable.
This is first result in KAM Theory for PDEs that regards
the persistence of a large measure set of invariant tori and hence may be viewed
as a strict extension to the infinite-dimensional setting of the classical KAM theorem.
	
\medskip\medskip


\smallskip\smallskip

\end{abstract}
%%%%%%%%%%%%%%%%%%%%%%%%%%%%%%%%%%%%%%%%%%%%%%%%%%%%%%%%%%%%%%%%%%%%%%%%%% 
%%%%%%%%%%%%%%%%%%%%%%%%%%%%%%%%%%%%%%%%%%%%%%%%%%%%%%%%%%%%%%%%%%%%%%%%%% 

\medskip

%%%%%%%%%%%%%%%%%%%%%%%%%%%%%%%%%%%%%%%%%%%%%%%%%%%%%%%%%%%%%%%%%%%%%%%%%% 
\begin{spacing}{0.92}
\tableofcontents
\end{spacing}
%\renewcommand{\baselinestretch}{1.00}\normalsize
%%%%%%%%%%%%%%%%%%%%%%%%%%%%%%%%%%%%%%%%%%%%%%%%%%%%%%%%%%%%%%%%%%%%%%%%%% 

\vspace{.6cm}

%%%%%%%%%%%%%%%%%%%%%%%%%%%%%%%%%%%%%%%%%%%%%%%%%%%%%%%%%%%%%%%%%%%%%%%%%% 
%%%%%%%%%%%%%%%%%%%%%%%%%%%%%%%%%%%%%%%%%%%%%%%%%%%%%%%%%%%%%%%%%%%%%%%%%% 
\zerarcounters 
\section{Introduction} \label{intro} 
%%%%%%%%%%%%%%%%%%%%%%%%%%%%%%%%%%%%%%%%%%%%%%%%%%%%%%%%%%%%%%%%%%%%%%%%%% 
%%%%%%%%%%%%%%%%%%%%%%%%%%%%%%%%%%%%%%%%%%%%%%%%%%%%%%%%%%%%%%%%%%%%%%%%%% 
	
Understanding the qualitative behaviour of the solutions of nonlinear Hamiltonian PDEs is a pivotal and fascinating question
in Analysis and Mathematical Physics. Once the local or global well-posedness of a given equation is established, the behaviour
of the solutions can still be extremely varied. In the case of compact manifolds, for example,
one expects coexistence of both regular and chaotic dynamics.
This diversity is reflected in the existing literature, where results come into two different flavours:
either broad statements, such as upper bounds on the growth in time of Sobolev norms, holding for all initial data, 
possibly within some small ball around the origin, or
precise statements regarding the time evolution of solutions restricted to very special initial data.
This means that one sways between providing statements that apply universally but might be too general to be insightful
and making very specific claims which might be  limited to extremely special settings.
In other words, resorting to a hyperbole, one is fighting against the well-known dichotomy of whether
to say nothing about everything or everything about nothing. 

The aim of the present paper is to provide a statistical description of the global dynamics of a
particular class of PDEs, parametrised by family of potentials, by investigating the behaviour of
a significant portion of sufficiently regular initial data.
Following the path laid out by Kolmogorov in his famous address at the 1954 ICM conference,
our goal is finding ``which of the properties of dynamical systems are `typical' for `arbitrary'
[Hamiltonian functions]'' \cite[p.~358]{kol-ICM}, with ``typical'' and ``arbitrary'' meaning, in our context,
for the majority of initial data and Hamiltonians, respectively, in the sense of measure.
Still following Kolmogorov, ``the approach from standpoint of measure theory appears
to be quite reasonable and natural as viewed from physics, [although]
its application is hampered by the absence of a natural measure in function spaces'' \cite[pp.~358-359]{kol-ICM}.
Kolmogorov's insight gave birth to KAM theory,  which, in its classical formulation,
concerns Hamiltonian systems which are close to an integrable  one and shows
that statistically  the dynamics is qualitatively the same as that of the integrable system,
namely a linear dynamics which is dense on a Lagranian invariant torus.

KAM theory for Hamiltonian PDEs has ben studied by many authors (see Section \ref{S1.2} for a brief overview),
but most of the results are about periodic or quasi-periodic solutions corresponding to finite-dimensional invariant tori. 
%It is important to note that, in the infinite-dimensional context, both periodic and quasi-periodic 
Such solutions are not expected to be typical w.r.t.~any reasonable measure. Indeed, 
if one considers integrable Hamiltonian PDEs, such as the cubic NLS or the KdV on the circle,
typical solutions are almost-periodic functions\footnote{Here and henceforth,
following Bohr \cite{bohr,bohr1} and Bochner \cite{boch},
the set of almost-periodic functions on a Banach space $(X,\|\cdot\|_X)$ 
is meant as the closure w.r.t.~the uniform topology
of the set of trigonometric polynomials with values in $X$.}
{the images of whose hulls}\footnote{\label{hull}{The hull of an almost-periodic function $f(t)$ on a
Banach space $(X,\|\cdot\|_X)$ is the set of functions defined as
the closure w.r.t.~the uniform topology of the set $\{ f_\tau(t) = f(t + \tau) : \tau\in\RRR\}$
(see for instance \cite{fink} or \cite{JM}).}}
are invariant tori. It is reasonable to surmise that, for close-to-integrable PDEs,
a generic perturbation preserves such a property, if any.

A KAM result on persistence of maximal tori for close-to-integrable PDEs,
as general as in the finite-dimensional case, is still out of reach,
even when confined to semi-linear PDEs on the circle.
A possible approach consists in introducing two significant simplifications:
first, one restricts the analysis to small solutions and
hence one considers perturbations of a linear integrable system; secondly,
one considers families of PDEs parametrized by a smooth convolution potential.
Even in this simplified setting, results regarding existence of  almost-periodic solutions
are relatively few, and only hold for a zero measure set of initial data
(again we refer to Sections \ref{S1.2} and \ref{S1.3} for a brief overview).
In this paper, we follow such an approach, but,
differently from the previous literature, we randomize both the potential and the initial data.
In the context of Hamiltonian PDEs, the idea of randomizing the initial data 
was ushered in by Bourgain \cite{Bomanca1,BouB,Bijm}, %and/or the potential,
and, since then, has been used to prove local well-posedness for equations with very low regularity \cite{BZ1},
and has opened up a vast area of research; see for instance \cite{CO,D1,NS}.

%%%%%%%%%%%%%%%%%%%%%%%%%%%%%%%%%%%%%%%%%%%%%%%%%%%%%%%%%%%%%%%%%%%%%%%%%% 
\subsection{Main results}\vspace{-.2cm}
%%%%%%%%%%%%%%%%%%%%%%%%%%%%%%%%%%%%%%%%%%%%%%%%%%%%%%%%%%%%%%%%%%%%%%%%%% 

We consider a family of nonlinear Schr\"odinger (NLS) equations on the circle, parametrized
by a smooth convolution potential.
Following a long-established tradition, dating back to Bourgain \cite{B2},
we consider explicitly the case of a quintic nonlinearity,
even though the results we describe below extend to any analytic nonlinearity, up to notational intricacies.
Our main result can be stated informally as follows.\footnote{Throughout this section,
we use words such as ``many'' and ``most" (and the like) without specifying their precise meaning,
by referring to Section \ref{setup} for a rigorous definition, once the proper functional spaces and measures have been introduced.}
\vspace{-0.2cm}
\begin{quoting}
	\textit{For most choices of the convolution potential, many small Gevrey initial data give rise to
almost-periodic solutions 
whose hulls describe invariant tori.}
\end{quoting}
\vspace{-0.2cm}

More precisely, we prove that, for many convolution potentials and, correspondingly, for many Gevrey initial data,
it is possible to construct an immersion of the torus  $\TTT^\ZZZ$ in the phase space such
that the dynamics stays on the torus for all times.
As a consequence, for many initial data close to the origin the Gevrey norm is approximately constant in time
and hence the elliptic fixed point at the origin is Lyapunov statistically stable w.r.t.~such a norm.

The invariant tori we construct are in general infinite-dimensional, and might well not be submanifolds
in the phase space. To deal with such an issue, we introduce an appropriate functional setting allowing us
to construct regions of the phase space where there is a Cantor foliation into invariant maximal tori,
in the spirit of Kuksin and P\"oschel approach \cite{kupo}.
In this respect, our result is a true generalization of a finite-dimensional KAM result
inasmuch it ensures that most of the phase space is filled by invariant tori.
A more thorough, still informal, statement of our result is the following one.
\vspace{-0.2cm}
\begin{quoting}
	\textit{A large measure set of Gevrey initial data evolve globally on invariant tori, and in a region of the phase
		space, corresponding to positive actions, such invariant tori are submanifolds and are the leaves of a Cantor foliation.
	}
\end{quoting}
\vspace{-0.1cm}
%\emph{a large measure set of Gevrey initial data evolve globally on invariant tori, and in a region of the phase
%space, corresponding to positive actions, such invariant tori are submanifolds and are the leaves of a Cantor foliation.}
%
By their very nature, the existence of these invariant tori ensures stability.
In finite dimension this inference
is trivial, but in the infinite-dimensional context its occurrence
heavily depends on the chosen functional space
and on the fact that the tori are immersed in such a space.
As a counterpart, while in finite dimension diffusive
solutions are extremely hard to find, in the infinite dimensional context, our result suggests to look for diffusive
solutions in spaces with regularity lower than Gevrey; this is consistent with the results in \cite{CKSTT,HGP}.

%%%%%%%%%%%%%%%%%%%%%%%%%%%%%%%%%%%%%%%%%%%%%%%%%%%%%%%%%%%%%%%%%%%%%%%%%% 
\subsection{Context and background}\label{S1.2}\vspace{-.2cm}
%%%%%%%%%%%%%%%%%%%%%%%%%%%%%%%%%%%%%%%%%%%%%%%%%%%%%%%%%%%%%%%%%%%%%%%%%% 

Almost-global stability\footnote{
{That is, stability for finite but long time, as opposed to global stability (or stability \emph{tout court}), which refers to all times.}}
of solutions of semilinear PDEs close
to an elliptic fixed point has been studied extensively, 
starting from the early result by Bourgain~\cite{Bou96b}. Later, Bambusi and Gr\'ebert introduced
a method based on Birkhoff normal form to prove long but finite time stability of Sobolev norms for small initial data~\cite{BG}.
The latter approach has been generalized in many contexts, for instance
for more regular initial data \cite{yuan,BMP3,cong}, even in higher dimensions \cite{FG},
and in the case of quasi-linear or fully nonlinear PDEs~\cite{Delort,Berti-Delort,FGI20}.
All these results actually deal with families of PDEs depending on parameters,
and the stability result is obtained for many values of the parameters;
more recently, Bernier, Faou and Gr\'ebert considered the NLS without external parameters and
proved long time stability for many initial data~\cite{BFG}.
We stress that, because of the very way the problem is set up,
all the stability results mentioned above hold only for finite time.
Indeed, all of them but the last one
hold for all initial data in some open neigborhood of the origin, so infinite time stability can not be expected -- even in finite dimension.
On the other hand in \cite{BFG}, still relying on a Birkhoff normal form approach, 
it is the initial data that are modulated, instead of the parameters.

As a general consideration, if one aims to obtain information on the solutions
over an arbitrarily long -- possibly infinite -- time scale, then one either manages to provide an upper bound on the growth of
the Sobolev norm which is non-uniform in time, as, for instance, in \cite{Bou96b,Staff,Pl,BPVT} and references therein,
or one  has to restrict the analysis to special solutions.
In fact, many authors proved the existence of solutions exhibiting unstable behaviour in various contexts; see
for instance \cite{CKSTT,HPTV,HGP,GHHMP} and references therein.
Also the existence of recurrent solutions -- that is solutions which are either periodic
or quasi-periodic or almost-periodic in time -- has
been widely studied starting from the pioneering results from the early
nineties~\cite{K,K2,K3,W,CW,B1,P,kupo}, which laid down the foundations and
have been thereafter fruitfully extended to a number of cases,
either in higher dimension settings~\cite{Bann,GY,EK,GXY}
or for unbounded nonlinearities~\cite{IPT,BBM,BBHM}, only to mention a few.
All the literature mentioned above concerns only
the first two kinds of recurrent solutions, the case of almost-periodic solutions being much harder to deal with.

%It is important to note that, in the infinite-dimensional context,
%both periodic and quasi-periodic solutions are not expected to be ``typical'' w.r.t.~any reasonable 
%measure. Indeed, if one considers integrable Hamiltonian PDEs, such as the cubic NLS on the circle,
%most solutions are almost-periodic, 
%and it is reasonable to surmise that, for close-to-integrable PDEs, a generic perturbation preserves such a property, if any.
The search of almost-periodic solutions is a difficult task due to the presence
of extremely ``bad'' small divisors, and this is the reason why the literature on this topic
is relatively scarce, if compared with the periodic and quasi-periodic counterparts,
and, in fact, it reduces to a handful of papers: an almost exhaustive list is given by
\cite{Bou96a,Bou96b,Po,yuan,BMP3,BMP1,congyuan,MP,CMP,BMP2,CGP,CGP2,cong}.
Moreover, in all such papers, infinitely many external parameters are needed, and one is usually only able to prove
the existence of few almost-periodic solutions, with either very high regularity or very special form.
A remarkable exception is \cite{BGR}, where the authors prove the existence of 
almost-periodic solutions for the NLS equation on the circle with no external parameters. 
A key feature of the strategy of \cite{BGR} is to construct solutions which are approximately supported on 
a set of Fourier modes which, although infinite, is sparse in $\ZZZ$; a similar strategy is used in \cite{BMP2}
to obtain low regularity almost-periodic solutions for the NLS equation with external parameters.
These solutions are dense on infinite-dimensional but not maximal tori and hence they are not typical,
and in fact cover a zero measure set.
%however, the corresponding invariant tori
%cover a zero measure set and are not maximal -- and hence not ``typical'' -- since they involve
%a set of Fourier modes which, although infinite, is sparse in $\ZZZ$,
%as in the case of the almost-periodic solutions with low regularity studied in \cite{BMP2}.}
We may also cite \cite{mofra}, where no parameter modulation is needed due to the special structure
of the equation considered therein. Regarding almost-periodic solutions
still in the infinite-dimensional context but not directly related to PDEs,
we mention \cite{frospewa,Po90,CP,CGP2}, where ODEs on a lattice are studied.

%%%%%%%%%%%%%%%%%%%%%%%%%%%%%%%%%%%%%%%%%%%%%%%%%%%%%%%%%%%%%%%%%%%%%%%%%% 
\subsection{Outline of the paper}\label{S1.3}\vspace{-.2cm} %\medskip
%%%%%%%%%%%%%%%%%%%%%%%%%%%%%%%%%%%%%%%%%%%%%%%%%%%%%%%%%%%%%%%%%%%%%%%%%% 

Let us now enter more in detail on the content of our paper,
and give a more accurate description of the results we obtain -- though referring to Section \ref{setup} for a formal statement.
We study the NLS equation on the circle given by
\begin{equation} \label{1.1}
\ii u_t - u_{xx} + \VV * u + \e |u|^4u =0 , \qquad x\in\TTT, 
\end{equation}
where $\VV$ is a convolution potential, which we assume to have finite but arbitrarily large regularity,
and we consider initial data $u(x,0) = \WW(x)$ with Gevrey regularity.
Thus we are in the same context as in \cite{CGP}, where the existence of some almost-periodic solutions
is proved.  The purpose of this paper is to show that such solutions have positive measure and
provide a Cantor foliation in a large region of the phase space.

We briefly recap the scheme followed in \cite{CGP}, keeping in mind our purpose of counting the solutions.
In \cite{CGP} the solutions are parametrized by the solutions
to the NLS equation linearized about the origin, in the following sense.
If we expand both the convolution potential and the looked-for solution into their Fourier series,
we can rewrite \eqref{1.1} as
\begin{equation} \label{1.2}
\ii (u_{j})_t + \left( j^2 + V_j \right) u_j + \e \left( |u|^4u \right)_j = 0 , \qquad j \in\ZZZ,
\end{equation}
where
\begin{equation} \label{1.2bis}
V_j := \frac{1}{2\pi}\int_0^{2\pi} e^{-\ii j x}\VV(x) \, \der x , \qquad
u_j := \frac{1}{2\pi}\int_0^{2\pi} e^{-\ii j x}u(x,t) \, \der x , %\qquad j\in\ZZZ ,
\end{equation}
and hence
\begin{equation} \label{1.2ter}
\left( |u|^4 u \right)_j := \frac{1}{2\pi}\int_0^{2\pi} e^{-\ii j x} \left( |u(x,t)|^4 u(x,t) \right) \der x .%\qquad j\in\ZZZ .
\end{equation}
Then, following a classical approach by Moser, we introduce, for reasons to become clear below, the auxiliary equation
\begin{equation} \label{1.4}
\ii (u_{j})_t + \left( \om_j + \eta_j \right) u_j + \e \left( |u|^4u \right)_j = 0 , \qquad j \in \ZZZ ,
\end{equation}
with $\om\in \RRR^\ZZZ$ such that $\{\om_j - j^2\}_{j\in \ZZZ}\in \ell^\io(\RRR)$
and $\eta=\eta(\e) \in \ell^{\io}(\RRR)$ such that $\eta(0)=0$.
We call \eqref{1.4} the modified equation, $\om$ the frequency and $\eta$ the counterterm.
Any solution of the linear equation 
\begin{equation} \label{1.3}
\ii (u_{j})_t +  \om_j u_j = 0  ,  \qquad j \in \ZZZ ,
\end{equation}
obtained from \eqref{1.4} by setting $\e=0$, is of the form
\begin{equation} \label{1.5}
u_{\rm lin}(x,t) = 
\sum_{j\in\ZZZ} c_j e^{\ii j x} e^{\ii \om_j t } ,
\end{equation}
with $c$ a suitable sequence in $\CCC^\ZZZ$ decaying fast enough for the series \eqref{1.5} to be summable.
Theorem 2.13 in \cite{CGP} shows that the following happens:\footnote{
This is a Moser-like counterterm theorem in the infinite-dimensional setting.
{A formal statement, after the appropriate notation has been introduced,  is provided by Theorem \ref{moser} in Section \ref{secmoser},
with emphasis on the properties we aim to rely upon in the forthcoming discussion.}}
consider any frequency $\om$ satisfying an appropriate non-resonance condition,
consider moreover any sufficiently small Gevrey solution \eqref{1.5} of the linear equation \eqref{1.3};
then there exists a counterterm $\eta$ such that the modified equation \eqref{1.4}
has an almost-periodic solution with frequency $\omega$ bifurcating from the solution \eqref{1.5}.
Therefore, for many values of  $\om$, the corresponding solution \eqref{1.5}
may be continued into a nonlinear solution to \eqref{1.4}, provided the sequence $\eta$ is suitably chosen.
On the other hand, if we fix $c$ small enough and decaying fast enough,
then, for $V$ in a large measure set, it is possible to fix a non-resonant vector $\omega$
in such a way that
\begin{equation} \label{1.6}
\om_j + \eta_j=j^2 +V_j , \qquad j\in\ZZZ , 
\end{equation}
so that the solution to the nonlinear equation \eqref{1.4} is a solution to the original NLS equation \eqref{1.2} as well.
Moreover, Proposition 2.21 in \cite{CGP} ensures that each component $\eta_j$ of the counterterm
admits a suitable asymptotic expansion in terms of $j$ which allows to choose the convolution potential $\VV$ in \eqref{1.1}
with arbitrary finite regularity.\footnote{See also Proposition \ref{prop2.21} in Section \ref{secmoser},
where the result is recalled --  and stated in the form most suited to our purposes.}
In conclusion, we find that  the NLS equation \eqref{1.1}, with $\VV$ a regular function,
admits, for most choices of $\VV$, almost-periodic solutions
which are parametrized by the sequences $c$ which identify the solutions to the linearized equation.

However, from the construction described above, it remains unclear whether such a parametrization is injective,
namely whether different linear solutions give rise to distinct  almost-periodic solutions for the nonlinear system.
In this setting, one can say at best that the solutions are  uncountably many
 {(see \cite{Po} for a similar predicament).}
In finite dimension, in the presence of some twist condition, 
the problem is typically bypassed by proving that different linear solutions 
give rise to nonlinear solutions with different frequencies,
{since the nonlinear solutions ``inherit" the frequencies of the linear solutions.}
In the infinite-dimensional setting, this is not an easy 
task and, in any case, it would not guarantee the almost-periodic solution to have positive measure.
The  key of our strategy is to show that in fact there is a bi-Lipschitz (non-symplectic) invertible map from a ball 
in the space of initial data to the space of linear solutions. 
More precisely, we proceed as follows.
We revisit the two results in \cite{CGP} mentioned above, that is Theorem 2.13 and Proposition 2.21, 
and prove that one can keep accurate control  on the regularity of the solution and of the counterterm w.r.t.~the linear solution. This control,
which comes from the very explicit construction via ``diagrammatic expansion'' of \cite{CGP},  allows us to 
 deviate from the analysis performed in \cite{CGP} and consider a different implicit function problem (see \eqref{system} in Section \ref{eppoi}),
which, besides ensuring the compatibility condition \eqref{1.6} to be satisfied,
takes into account also the initial datum by requiring that $u(x,0)=\WW(x)$.
This allows us to define a bi-Lipschitz map from the space of linear solutions
to the space of initial data, which, accordingly, may be used as further parameters
(see \eqref{systemboh} in Section \ref{eppoi}).
Thus, we are able to express the frequencies in terms not only of the potentials but also of the initial data
and to prove that, for many choices of the convolution potential and of the initial datum,
the corresponding frequency satisfies the non-resonance condition.
This ensures the existence of the invariant tori on which the almost-solutions lie,
provided both the convolution potential and the initial datum are taken in large measure sets.
Finally, we prove the Lyapunov statistical stability by showing that,
for convolution potentials in a large measure set,
the set of initial data such that the corresponding frequency is non-resonant has large positive measure as well.

To complete the comparison with the classical KAM theorem, we study the geometrical structure of the solutions.
In fact, by construction, their hulls have invariant tori as images.
However, as already mentioned,
such tori may not be submanifold in the infinite-dimensional setting. We resort once more to the explicit expression
provided by Theorem 2.13 in \cite{CGP} and, in addition,
we make use of the last ingredient of our analysis: the construction of
a bi-Lipschitz map expressing, at fixed potential $V$, the initial data
in terms of the amplitudes of the linear solution (see \eqref{WWW} in Section \ref{outro}).
By  exploiting such results, we  construct a region of the phase space Cantor-foliated in invariant tori.
Moreover, we prove that those tori  which do not belong to the foliation are still immersed tori.

%%%%%%%%%%%%%%%%%%%%%%%%%%%%%%%%%%%%%%%%%%%%%%%%%%%%%%%%%%%%%%%%%%%%%%%%%% 
\subsection{Methods and tools}\vspace{-.2cm} %\medskip
%%%%%%%%%%%%%%%%%%%%%%%%%%%%%%%%%%%%%%%%%%%%%%%%%%%%%%%%%%%%%%%%%%%%%%%%%% 

We conclude this introductory excursus with a few words about the main techniques we rely upon in our approach.
The proof of Theorem 2.13 in \cite{CGP},
{which is one of the key ingredients of our analysis,
is based on the so-called tree formalism, first introduced by Feldman and Trubowitz \cite{FT}
and, more systematically, by Gallavotti~\cite{Galla}}
and inspired to the Feynman graphs used in Quantum Field Theory \cite{GMP,GM0}.
In the context of recurrent solutions for Hamiltonian PDEs,
such a method has already been exploited in \cite{GM1,GMP,GP,GP2,MP}.
More recently trees expansions have also been used  by Deng and Hani to obtain a full derivation
of the wave kinetic equation \cite{DH1,DH2}.

The tree formalism is based on a graphical representation of the formal expansion of the solution;
this allows, first, to identify readily the terms where small divisors accumulate, and, then, to show that,
when summed together, these terms turn out to compensate each other in such a way that 
the overall contribution they provide is not too large. The advantage of using the tree formalism is that,
once the convergence of the series is proved, the graphical representation provides a fairly explicit expression of the solution. 
In \cite{CGP}, a main point of the analysis where such an expression is used is in the proof of Proposition 2.21,
when one needs to show that, if one suitably chooses the frequency, then the counterterm shares the same asymptotic expansion as the frequency.
In the present paper we heavily exploit the explicit expression of the solution both to obtain
the bi-Lipschitz map from linear  solutions to initial data  needed to solve the implicit function problem \eqref{system},
and to prove the topological properties needed to construct of the Cantor foliation.

%%%%%%%%%%%%%%%%%%%%%%%%%%%%%%%%%%%%%%%%%%%%%%%%%%%%%%%%%%%%%%%%%%%%%%%%%% 
\subsection{Conclusions}\vspace{-.2cm} %\medskip
%%%%%%%%%%%%%%%%%%%%%%%%%%%%%%%%%%%%%%%%%%%%%%%%%%%%%%%%%%%%%%%%%%%%%%%%%% 

The results of this paper provide the first extension of the KAM theorem to a PDE,
specifically the NLS equation with parameters, 
regarding both the geometric structure of the maximal invariant tori and their measure in the phase space.

The main new idea is to parametrize the tori in terms of neither the actions nor the linear solutions but, instead,
the initial data -- a strategy dictated by the need to overcome the lack of control on the twist.
The idea by itself is very simple,
however the implementation in the infinite-dimensional context is rather delicate
and requires an extremely accurate quantitative control
on the solutions and their Lipschitz dependence on all the parameters. 
To obtain such a control we take advantage in a substantial way of the explicit construction of the
solutions and of the refined estimates obtained in \cite{CGP}.

\medskip\medskip

\noindent
\emph{Acknowledgements}.
We  thank E.~Haus and M.~Berti for many useful comments on the manuscript.
L.B., L.C.~and M.P.~have been  supported by the  research project  PRIN 2020XBFL ``Hamiltonian and dispersive PDEs", 
L.C.~has been supported by the research PRIN 2022HSSYPN ``Turbulent Effects vs Stability in Equations from Oceanography'' (TESEO),
G.G.~has been partially supported by the research project PRIN 20223J85K3
``Mathematical Interacting Quantum Fields'' of the Italian Ministry of Education and Research (MIUR).

%%%%%%%%%%%%%%%%%%%%%%%%%%%%%%%%%%%%%%%%%%%%%%%%%%%%%%%%%%%%%%%%%%%
%%%%%%%%%%%%%%%%%%%%%%%%%%%%%%%%%%%%%%%%%%%%%%%%%%%%%%%%%%%%%%%%%%%
\zerarcounters 
\section{Set up and main results}
\label{setup} 
%%%%%%%%%%%%%%%%%%%%%%%%%%%%%%%%%%%%%%%%%%%%%%%%%%%%%%%%%%%%%%%%%%%
%%%%%%%%%%%%%%%%%%%%%%%%%%%%%%%%%%%%%%%%%%%%%%%%%%%%%%%%%%%%%%%%%%%
	
Consider the Cauchy problem associated to the NLS equation
\begin{equation}\label{nls}
\begin{cases}
\ii u_t - u_{xx} + \VV * u + \e |u|^4u =0 , \qquad x\in\TTT, \\
u(x,0) = \WW(x) ,
\end{cases}
\end{equation}
where $\e$ is a small real parameter, $\VV$ is a convolution potential and $\WW$ is the initial datum.
We make the following assumptions on the regularity of $\VV$ and $\WW$.

Let $\calF$ denote the Fourier transform operator and let
${V}=\{V_j\}_{j\in\ZZZ}:=\calF(\VV)$ be the sequence of the Fourier coefficients of $\VV$,
defined according to \eqref{1.2bis} -- if they exist.
Then, if we set, for $k\in\NNN \cup \{0\}$,\footnote{For $k=0$ the norm in \eqref{LNinfty} reduces to the sup-norm $\|\cdot\|_{\io}$.}
\begin{equation}\label{LNinfty}
\ell^{k,\io}(\RRR) := \Big\{ \tx=\{\tx_j\}_{j\in\ZZZ} \,:\, \tx_j\in\RRR, \, \|\tx \|_{k,\io}:=\sup_{j\in\ZZZ}|\tx _j|\jap{j}^k < \io \Big\},
\end{equation}
with
\begin{equation} \nonumber %\label{jap}
\jap{j}:=\max\{1,|j|\} ,
\end{equation}
we assume that
$\VV$ admits a Fourier series expansion and that $V\in\ell^{N,\io}(\RRR)$ for some $N\in \NNN$.
Note that such a condition is satisfied if the convolution potential $\VV$ is of class $C^N$.

Regarding the initial datum $\WW$, we fix $\al\in(0,1)$ and $s>0$, and assume that
the sequence $W = \{W_j\}_{j\in\ZZZ} := \calF(\WW)$ belongs to $\mathtt g(s,\al)$, where
\begin{equation}
\label{gsa0}
\mathtt g({s,\alpha})  :=\Bigl\{\tx =\{\tx _j\}_{j\in\ZZZ}\in \ell^\io(\CCC) : 
\|\tx \|_{s,\alpha}:= \sup_{j\in\ZZZ} |\tx_j| e^{s\jap{j}^\alpha} <\io\Bigr\} .
\end{equation}

With a slight abuse of language, in the rest of the paper, for any fixed $N$, $\al$ and $s$,
we call \emph{potential} the sequence $V \in \ell^{N,\io}(\RRR)$ and \emph{initial datum} the sequence $W\in\mathtt g(s,\al)$,
with $\mathtt g(s,\al)$ being referred to as the \emph{phase space} of the system described by the equation \eqref{nls}.

%%%%%%%%%%%%%%%%%%%%%%%%%%%%%%%%%%%%%%%%%%%%%%%%%%%%%%%%%%%%%%%%%%%%%%%%%% 
\subsection{Almost-periodic solutions for a large measure set of potentials}
\vspace{-.2cm}
%%%%%%%%%%%%%%%%%%%%%%%%%%%%%%%%%%%%%%%%%%%%%%%%%%%%%%%%%%%%%%%%%%%%%%%%%% 

Our purpose is to show that, for ``most" choices of the potential $V$ in a ball of $\ell^{N,\io}(\RRR)$,
a ``large'' subset of initial data $W$ in a ball of the phase space $\mathtt g (s,\al)$ gives rise to almost-periodic solutions.

In order to make precise what we mean above by ``most'' and ``large'', we have to introduce
for both spaces $\ell^{k,\io}(\RRR)$ and $\mathtt g({s,\alpha})$ a proper measure.\footnote{{It would be more appropriate
to speak about \emph{probability measure} rather than measure, since we assign value 1 to the measure
of the entire space, however, for simplicity, throughout the paper we call measures \emph{tout court} the probability measures introduced here.}}
To do that, we need some notation. Given a normed space $(X,\|\cdot\|_X)$,
for any $\rho>0$ we set $\matU_\rho(X):=\{\mathtt{x}\in X : \|\mathtt{x}\|_X<\rho\}$ 
and let $\ol{\matU}_\rho(X)$ denote the closure of $\matU_\rho(X)$ w.r.t.~the norm $\|\cdot\|_X$;
{more generally,
for $X\subset \CCC^\ZZZ$ and any $\mathtt{x}_0 \in \CCC^\ZZZ$,
we define $\matU_{\rho,\mathtt{x}_0}(X):=\mathtt x_0+ \matU_{\rho}(X)$= $ \{\mathtt{x}\in \CCC^\ZZZ : \|\mathtt{x} - \mathtt{x}_0\|_X<\rho\}$ 
and call $\ol{\matU}_{\rho,\mathtt{x}_0}(X)$ the closure of $\matU_{\rho,\mathtt{x}_0}(X)$.}\footnote{{By construction,
one has $\matU_{\rho}(X)=\matU_{\rho,\mathtt{0}}(X)$.}}
If  $X$ is a  real weighted sequence space based on $\ell^\io(\RRR)$, i.e.
\[
X=\Bigl\{ \tx\in \ell^\io(\RRR): \quad \|\tx\|_X:= \sup_{j\in\ZZZ} \mathtt w_j |\tx_j|<\io \Bigr\} ,
\]
for some sequence of positive \emph{weights} $\{\mathtt w_j\}_{j\in\ZZZ}$,
then,
{for any $\mathtt{x}_0\in \RRR^\ZZZ$, the ball $\ol{\matU}_{\rho,\mathtt{x}_0}(X)$} is closed also w.r.t.~the (weaker) product topology. 
This allows us  to endow $\ol{\matU}_{\rho,\mathtt{x}_0}(X)$ with the measure induced by the product measure
on the ball $\ol{\matU}_1(\ell^\io(\RRR))$, by proceeding as follows. For any set
\[
A = \prod_{j\in\ZZZ} A_j\subseteq  \prod_{j\in\ZZZ}
{[\mathtt{x}_{0j} - \rho^{-1} \mathtt w_j , \mathtt{x}_{0j} + \rho^{-1}\mathtt w_j ]} ,
\]
with $A_j$ a Lebesgue-measurable set in $\RRR$ for all $j\in\ZZZ$,
we define the measure of $A$ as
\begin{equation} \nonumber %\label{eq:measur1e}
\meas_{X,\rho}(A) := \lim_{h\to\io} \prod_{j=-h}^{h} \meas_{X,\rho} (A_j) ,\qquad \meas_{X,\rho}(A_j) :=\frac{m(A_j)}{2\rho}\mathtt w_j,
\end{equation}
where $m(\cdot)$ denotes the Lebesgue measure. 
Similarly, if  $X$ is a  complex weighted sequence space based on $\ell^\io(\CCC)$, i.e.
\[
X=\Bigl\{ \tx\in \ell^\io(\CCC): \quad \|\tx\|_X:= \sup_{j\in\ZZZ}  |\tx_j| \mathtt w_j <\io \Bigr\} ,
\]
we endow
{$\ol{\matU}_{\rho,\mathtt{x}_0}(X)$}
with the measure induced by the product measure on the complex ball $\ol{\matU}_1(\ell^\io(\CCC))$;
this means that, for any set
\[
A = \prod_{j\in\ZZZ} A_j\subseteq  \prod_{j\in\ZZZ} \bigl\{ \tx_j\in\CCC:
{ |\tx_j -\tx_{0j} | \rho^{-1}\mathtt w_j \le 1 } 
\bigr\}  ,
\]
with $A_j$ a Lebesgue-measurable set in $\CCC$ for all $j\in\ZZZ$,
we define the measure of $A$ as
\begin{equation} \nonumber %\label{eq:measure2}
\meas_{X,\rho}(A) 
:= \lim_{h\to\io} \prod_{j=-h}^{h} \meas_{X,\rho} (A_j) ,\qquad \meas_{X,\rho}(A_j)=\frac{m(A_j)}{\pi\rho^2}\mathtt w^2_j,
\end{equation}
with $m(\cdot)$ still denoting the Lebesgue measure (in $\CCC$).

%%%%%%%%%%%%%%%%%%%%%%%%%%%%%%%%%%%%%%%%%%%%%%%%%%%%%%%%%%%%%%%%%%%%%%%%%% 
\begin{rmk}\label{pesiemisure}
\emph{
{Actually, we are interested in a few specific cases. Mainly, we deal with the spaces
$X\!=\!\ell^{N,\io}(\RRR)$ and $X=\mathtt{g}(s,\al) \subset \ell^{\io}(\CCC)$,
with weights $\mathtt w_j\!=\!\jap{j}^N$ and $\mathtt w_j=e^{s\jap{j}^\al}$, respectively,
and with $\mathtt{x}_0=\mathtt 0$ in both cases. Then, for any $\rho>0$, we write
$\mu_{1,\rho}:=\meas_{\ell^{N,\io}(\RRR),\rho}$ and $\mu_{2,\rho}:=\meas_{\mathtt{g}(s,\al),\rho}$,
and, for any $\rho_1,\rho_2>0$, we endow
$\ol{\matU}_{\!\rho_1}(\ell^{N,\io}) \times \ol{\matU}_{\!\rho_2}(\mathtt{g}(s,\al))$
with the measure
$\ol\mu_{\rho_1,\rho_2}:=\mu_{1,\rho_1} \times \mu_{2,\rho_2}$.
We also consider the space $\ell^{\io}(\RRR)$, with weights $\mathtt{w}_j=1$ and $\mathtt{x}_0=\varsigma:=\{j^2\}_{j\in\ZZZ}$,
and, in such case, we endow $\ol{\matU}_{1/2,\varsigma}(\ell^{\io}(\RRR))$ with the measure $\mu_0 :=\meas_{\ell^\io(\RRR),1/2}$.}
}
\end{rmk}
%%%%%%%%%%%%%%%%%%%%%%%%%%%%%%%%%%%%%%%%%%%%%%%%%%%%%%%%%%%%%%%%%%%%%%%%%% 

In order to state our first main result we need two last definitions.

Let $(X,\|\cdot\|_X)$ be a Banach space. The set of the \emph{almost-periodic functions} on $(X,\|\cdot\|_X)$
is the closure w.r.t.~the uniform topology of the set of trigonometric polynomials with values in $X$ \cite{bohr,bohr1}.
In fact, we look for an almost-periodic solution to \eqref{nls} which is
%A function $t\in\RRR \mapsto f(t)\in X$ is called an \emph{almost-periodic function} on $X$ if it is 
the restriction of a function of infinitely many angles $\f\in\TTT^\ZZZ$ computed at $\f=\om t$,
for some $\om\in\RRR^\ZZZ$. The vector $\om$ is called  the \emph{frequency} -- or frequency vector -- of the solution.

Let $\Omega\subset\RRR$ be open and let $f\!:\Omega \to \CCC$ be a $C^\io$ function. We say that $f$ is \emph{Gevrey of index $\s$} if
for any compact set $K \subseteq \Omega$ there exists a constant $C$ such that
\[
\sup_{x\in K}  \left| \frac{\partial^k}{\partial x^k} f(x) \right| \le C^{k} k!^{\s} \qquad \forall k \in \NNN  .
\]
%Let $f\!:\TTT\times\RRR \to\CCC$ be a $C^\io$ function. We say that $f(x,t)$ is \emph{Gevrey of index $\al_1$ in $x$ and index $\al_2$ in $t$} if
%(1) for all $x\in\TTT$ the map $t\mapsto f(x,t)$ is Gevrey with index $\al_1$, and
%(2) for all $t\in\RRR$ the map $x\mapsto f(x,t)$ is Gevrey with index $\al_2$ \cite{Gevrey}.

%%%%%%%%%%%%%%%%%%%%%%%%%%%%%%%%%%%%%%%%%%%%%%%%%%%%%%%%%%%%%%%%%%%%%%%%%% 
\begin{rmk}\label{Wgevrey}
\emph{
According to the last definition above, the assumption on the initial datum
$W$ is equivalent to requiring the function $\WW$ to be Gevrey of index $1/\al$.
}
\end{rmk}
%%%%%%%%%%%%%%%%%%%%%%%%%%%%%%%%%%%%%%%%%%%%%%%%%%%%%%%%%%%%%%%%%%%%%%%%%% 

Then, in Sections \ref{eppoi} and \ref{misura} we prove the following result.
	
%%%%%%%%%%%%%%%%%%%%%%%%%%%%%%%%%%%%%%%%%%%%%%%%%%%%%%%%%%%%%%%%%%%%%%%%%% 
\begin{theorem}\label{main}
Fix $s>0$,  $\al\in(0,1)$ and $N\geq 3$, and consider the Cauchy problem \eqref{nls}
with potential $V\in \ell^{N,\io}(\RRR)$ and initial datum $W\in \mathtt g(s,\al)$. 
There exists $\e_*=\e_*(s,\al,N) >0$ such that, 
for  any $\e\in (-\e_*, \e_*)$,
there exists a set {of potentials} $\mathcal G=\mathcal G(\e;\al,s,N)\subset \ol{\matU}_{\!{1/4}}(\ell^{N,\io}(\RRR))$ 
with positive measure such that the following holds. For any potential $V\in \mathcal G$, there exists a set {of initial data} 
$\mathcal T_V =\TT_V(\e;s,\al,N) \subset \ol{\matU}_{\!{1/2}}(\mathtt g(s,\al))$
with positive measure
such that, for any initial datum $W\in \mathcal T_V$, there exists
a unique solution $u(x,t)$
to the Cauchy problem \eqref{nls}, satisfying the following properties:
\vspace{-.2cm}
\begin{enumerate}
\itemsep0em
\item $u(x,t)$ is globally defined,
\item for all $t\in\RRR$ the map $x \mapsto u(x,t)$ is in $\calF^{-1}(\mathtt{g}(s,\al))$ and Gevrey of index $1/\al$,
\item for all $x\in \TTT$ the map $t \mapsto u(x,t)$ is almost-periodic and Gevrey of index $2/\al$.
\vspace{-.2cm}
\end{enumerate}
Moreover the measures of the sets $\mathcal G$ and $\mathcal T_V$ are both asymptotically full,
in the sense that both $\mu_{1,1/4}(\mathcal G)$ and $\mu_{2,1/2}(\mathcal T_V)$ tend to 1 as $\e$ tends to 0.
\end{theorem}
%%%%%%%%%%%%%%%%%%%%%%%%%%%%%%%%%%%%%%%%%%%%%%%%%%%%%%%%%%%%%%%%%%%%%%%%%% 

%%%%%%%%%%%%%%%%%%%%%%%%%%%%%%%%%%%%%%%%%%%%%%%%%%%%%%%%%%%%%%%%%%%%%%%%%% 
\begin{rmk} \label{N3}
\emph{
{All the results we discuss in this paper hold for any finite value of $N$. 
Here we consider only the case $N\ge 3$, since  we are interested in the case of convolution potentials with high regularity.
The case $N\le 2$ requires slightly different expressions \cite{CGP},
but it can be treated essentially in the same way.}
}
\end{rmk}
%%%%%%%%%%%%%%%%%%%%%%%%%%%%%%%%%%%%%%%%%%%%%%%%%%%%%%%%%%%%%%%%%%%%%%%%%% 

%%%%%%%%%%%%%%%%%%%%%%%%%%%%%%%%%%%%%%%%%%%%%%%%%%%%%%%%%%%%%%%%%%%%%%%%%% 
\begin{rmk}\label{remarkonregularity}
\emph{
{In \cite{BMP3} it is proved  that the Cauchy problem \eqref{nls} admits, in the class of regularity considered here,}
a unique solution (local well-posedness), where it is also showed that the solution
is defined at least for {sub-}exponentially long times (almost-global well-posedness).
In the light of Remark \ref{Wgevrey}, property 2 of Theorem \ref{main}, together with the measure estimates, ensures that,
for a large measure set of initial data, the time evolution
preserves their regularity for all times.
}
\end{rmk}
%%%%%%%%%%%%%%%%%%%%%%%%%%%%%%%%%%%%%%%%%%%%%%%%%%%%%%%%%%%%%%%%%%%%%%%%%% 

%%%%%%%%%%%%%%%%%%%%%%%%%%%%%%%%%%%%%%%%%%%%%%%%%%%%%%%%%%%%%%%%%%%%%%%%%% 
\begin{rmk}\label{riscalo}
\emph{
The NLS equation \eqref{nls} is covariant under the scaling (see Remark \ref{verifica} below)
\[
(u,\WW,\e)\mapsto (\lambda u, \lambda \WW, \lambda^{-4}\e) , \qquad \lambda>0 .
\]
As a consequence, we can reformulate Theorem \ref{main} as follows.
For any  $\rho\in(0, 1/2]$ set $\mathcal G_\rho = \mathcal G_\rho(\e;s,\al,N):= \mathcal G((2\rho)^4 \e;s,\al,N)$ and 
$\TT_{V,\rho}=\TT_{V,\rho}(\e;s,\al,N):= 2\rho\TT_V((2\rho)^4\e;s,\al,N)$;
then, %\subset\ol{\matU}_\rho(\mathtt{g}(s,\al))$
for any $\e\in (-\e_*, \e_*)$,
any $V\in \mathcal G_\rho \subset\ol{\matU}_{1/4}(\ell^{N,\io})$ %(\e;\al,s,N):= \mathcal G(\rho^4 \e;\al,s,N)$
and any
$W\in \TT_{V,\rho} \subset 
{\ol{\matU}_\rho(\mathtt{g}(s,\al))}$
there is a solution $u(x,t)$ to the  Cauchy problem \eqref{nls}
satisfying the properties 1, 2 and 3 in Theorem \ref{main}.
Moreover, the measures of the set $\mathcal G_\rho$ in the ball  $\ol{\matU}_{\!{1/4}}(\ell^{N,\io}(\RRR))$ %(\e;\al,s,N)$
and of the set $\TT_{V,\rho}$ in the ball
{$\ol\matU_\rho(\mathtt{g}(s,\al))$}
are both asymptotically full, in the sense that they tend to 1 as $\rho$ tends to $0^+$.
}
\end{rmk}
%%%%%%%%%%%%%%%%%%%%%%%%%%%%%%%%%%%%%%%%%%%%%%%%%%%%%%%%%%%%%%%%%%%%%%%%%% 

%%%%%%%%%%%%%%%%%%%%%%%%%%%%%%%%%%%%%%%%%%%%%%%%%%%%%%%%%%%%%%%%%%%%%%%%%% 
\begin{rmk}\label{random}
\emph{
We may regard Theorem \ref{main} in a probabilistic setting by randomising both the convolution potentials and the initial data w.r.t.~the uniform distribution.
More precisely, we may consider a random convolution potential and a random initial datum \cite{Bijm,DH1}
\[
\VV (x) = \sum_{j\in\ZZZ} \frac{1}{\jap{j}^N} v_j e^{\ii jx} ,
\qquad 
\WW (x) = \sum_{j\in\ZZZ} e^{-s\jap{j}^\alpha} w_j  e^{\ii jx} ,
\]
where, for all $j\in\ZZZ$, the coefficients $v_j$ are real i.i.d.~random variables uniformly distributed in $[-1/4,1/4]$,
while the coefficients $w_j$ are complex i.i.d.~random variables uniformly distributed in $\{z\in\CCC:|z|\le 1/2\}$. 
Then, Theorem \ref{main} ensures that almost-periodic solutions occur with positive probability.
In principle, one might also use other distributions, for instance one might assume the coefficients
to be centered normalized complex Gaussian variables \cite{DH2}; however, we did not investigate further such an issue.
%Similarly if the coefficients $f_j$ are complex  i.i.d.~random variables  uniformly distributed in the annulus $\{z\in\CCC:\delta\le|z|\le 1/2\}$,
%then Corollay \ref{biascoappeso} ensures that invariant maximal tori occur with positive probability.
}
\end{rmk}
%%%%%%%%%%%%%%%%%%%%%%%%%%%%%%%%%%%%%%%%%%%%%%%%%%%%%%%%%%%%%%%%%%%%%%%%%% 

%%%%%%%%%%%%%%%%%%%%%%%%%%%%%%%%%%%%%%%%%%%%%%%%%%%%%%%%%%%%%%%%%%%%%%%%%% 
\subsection{Almost-periodic solutions for a full measure set of potentials}\vspace{-.2cm}
%%%%%%%%%%%%%%%%%%%%%%%%%%%%%%%%%%%%%%%%%%%%%%%%%%%%%%%%%%%%%%%%%%%%%%%%%% 

From the perspective of physical applications, one is interested in considering the NLS equation with a fixed potential $V$,
say chosen in $\mathcal G$,  and showing that the set of
initial data giving rise to almost-periodic solutions has asymptotically full measure.
This is a delicate issue which requires a more careful analysis.

Let $\mathfrak G\subset \ol\matU_{\!{1/4}}(\ell^{N,\io}(\RRR))$ denote the set of potentials for which 
there exists a non-empty measurable set of initial data
$\matS_V\subset \ol\matU_{\!{1/2}}(\mathtt{g}(s,\al))$
such that for any $W\in \matS_V$ the solution  $u(x,t)$ to the Cauchy problem \eqref{nls}
satisfies the properties 1, 2 and 3 in Theorem \ref{main}.

%%%%%%%%%%%%%%%%%%%%%%%%%%%%%%%%%%%%%%%%%%%%%%%%%%%%%%%%%%%%%%%%%%%%%%%%%% 
\begin{rmk} \label{added}
\emph{
{By definition, with the notation of Theorem \ref{main},
the set $\mathfrak G$ contains the set $\mathcal G$, and one has $\matS_V = \mathcal T_V$ if $V\in\mathcal G$.
However, if compared with $\mathcal G$, the set $\mathfrak G$ is not only measurable but of full measure.
This can be seen by relying on the scaling properties enlightened in Remark \ref{riscalo}:
the measure of the set
{$\{ (V,W) \in  \ol{\matU}_{\!1/4}(\ell^{N,\io}(\RRR)) \times \ol{\matU}_{\!\rho}(\mathtt g(s,\al) ) : W \in \matS_V\}$}
%$\{ (V,W) \in \ol{\matU}_{\!1/4}(\ell^{N,\io}(\RRR)) \times \ol{\matU}_{\!\rho}(\mathtt g(s,\al) ): W \in \matS_V\}$
tends to 1 as $\rho$ tends to 0, and this excludes the existence of a subset $\mathfrak A \subset \ol\matU_{\!{1/4}}(\ell^{N,\io}(\RRR))$
such that $\mu_{1,1/4}(\mathfrak A)>0$ and $\mu_{2,\rho}(\matS_V \cap \ol{\matU}_{\!\rho}(\mathtt g(s,\al)))=0$ for all $V\in\mathfrak A$.}
}
\end{rmk}
%%%%%%%%%%%%%%%%%%%%%%%%%%%%%%%%%%%%%%%%%%%%%%%%%%%%%%%%%%%%%%%%%%%%%%%%%% 

%Let $\mathfrak G\subset \mathfrak G_0$ be the set of potentials such that $\matS_V$ has asymptotically full measure, namely 
%\begin{equation}
%	\label{clitennestra}
%	\liminf_{\rho\to 0^+} \mu_{2,\rho} (\matS_V \cap \ol\matU_{\rho}(\mathtt{g}(s,\al)))=1\,.
%\end{equation}

{In Section \ref{stab} we prove the following result.}

%%%%%%%%%%%%%%%%%%%%%%%%%%%%%%%%%%%%%%%%%%%%%%%%%%%%%%%%%%%%%%%%%%%%%%%%%% 
\begin{theorem}\label{eureka}
{
Consider the Cauchy problem \eqref{nls}, under the same hypotheses as in Theorem \ref{main}.}
Given  any positive sequence
$\varrho:=\{\rho_k\}_{k\ge1}\in \ell^{1/2(N+1)}(\NNN,\RRR_+)$,
there exists a  measurable set  of potentials 
$ {\mathfrak G}_\varrho\subset \mathfrak G$ with full measure
such that,  for any $V\in {\mathfrak G}_\varrho$, the set $\matS_V$ satisfies 
\begin{equation} \nonumber %\label{clitennestra2}
\liminf_{k\to \infty } \mu_{2,\rho_k} (\matS_V \cap \ol\matU_{\rho_k}(\mathtt{g}(s,\al)))=1 . 
\end{equation}
%
%	For all $\g_0>0$ small enough, 
%Any measurable set $\mathfrak B\subseteq\ol\matU_{\!{1/4}}(\ell^{N,\io}(\RRR))\setminus \mathfrak G$ satisfies
%	\[
%	\mu_{1,1/4}(\mathfrak B) = 0\,,
%	\]
%	namely  the $\mu_{1,1/4}$ external measure of $\mathfrak G$ is equal to one.
%such that for all $V\in \mathfrak G$ there exists a set of initial data $\matS_V\in  \ol\matU_{\!{1/2}}(\mathtt{g}(s,\al))$ with
%\begin{equation}
%		\label{clitennestra}
%	\lim_{\rho\to 0^+} \mu_{2,\rho} (\matS_V \cap \ol\matU_{\rho}(\mathtt{g}(s,\al)))=1
%\end{equation}
%such that for any $W\in \matS_V$, the unique solution $u(x,t)$ to the Cauchy problem \eqref{nls} satisfies items 1,2,3 of Theorem \ref{main}.
%	\red{
%	Equivalently, any measurable set $\mathfrak A \supseteq  \mathfrak G$ satisfies
%	\[
%	\mu_{1,1/4}(\mathfrak A) = 1\,.
%	\]
	%}
\end{theorem}

\begin{rmk} \label{woulldbeusefule}
\emph{
In Theorem \ref{eureka} one would like to consider  the set of potentials for which
$\liminf_{\rho\to 0^+ } \mu_{2,\rho} (\matS_V \cap \ol\matU_{\rho}(\mathtt{g}(s,\al)))=1$.
However, unfortunately, it is not clear to us even whether such a set is measurable.
This is why we restricted ourselves to taking the limit on a sequence. While some summability condition for the sequence $\varrho$
 is needed in our proof, we choose the space $\ell^{1/2(N+1)}$ only to simplify the argument.
}
\end{rmk}
%%%%%%%%%%%%%%%%%%%%%%%%%%%%%%%%%%%%%%%%%%%%%%%%%%%%%%%%%%%%%%%%%%%%%%%%%% 

%%%%%%%%%%%%%%%%%%%%%%%%%%%%%%%%%%%%%%%%%%%%%%%%%%%%%%%%%%%%%%%%%%%%%%%%%% 
\subsection{Lyapunov statistical stability}\vspace{-.2cm}
%%%%%%%%%%%%%%%%%%%%%%%%%%%%%%%%%%%%%%%%%%%%%%%%%%%%%%%%%%%%%%%%%%%%%%%%%% 

An important dynamical consequence of %Theorem \ref{eureka} 
our analysis is the following result regarding the global stability for
a large measure set of solutions, showing that  the elliptic fixed point $\und{u}(x,t)=0$ is statistically Lyapunov stable
w.r.t.~the Gevrey norm in the sense of the following result,
which is proved  in Section \ref{stab}.

%%%%%%%%%%%%%%%%%%%%%%%%%%%%%%%%%%%%%%%%%%%%%%%%%%%%%%%%%%%%%%%%%%%%%%%%%% 
\begin{theorem} \label{proposta}
Consider the Cauchy problem \eqref{nls}, under the same hypotheses as in Theorem \ref{main}.
There exists $\g^*>0$ such that for any $\g\in(0,\g^*)$
there exist $\rho=\rho(\g)>0$ and
{$\calG(\gamma) \subset \gotG$}
such that
\vspace{-.3cm}
\begin{enumerate}
\itemsep-0.2em
\item setting \footnote{Here and henceforth, given any set $A$ in any space $X$, $A^c$ denotes the complementary of $A$ in $X$.}
\[
{\Gamma(\g)}
: =\{(V,W)\in \ol\matU_{1/4}(\ell^{N,\io}(\RRR)) \times \ol\matU_{\rho}(\mathtt{g}(s,\al)) :  V\in \calG(\gamma)\,, W\in \matS_V \} ,
\]
one has $\ol\mu_{1/4,\rho}((\Gamma(\g))^c) \le \g/\g_*$,
\item one has $\mu_{1,1/4}((\calG(\gamma))^c) \le \sqrt{\g/\g_*}$,
\item one has $\mu_{2,\rho} (\matS_V^c \cap \ol\matU_{\rho}(\mathtt{g}(s,\al))
\le \sqrt{\g/\g_*} $ for all $V\in \calG(\g)$,
\item for any initial datum $W\in \matS_V\cap \ol\matU_{\rho}(\mathtt{g}(s,\al))$ the solution $u(x,t)$ to \eqref{nls} is such that
\[
\|\calF(u(\cdot ,t))\|_{s,\al} \le 2\|W\|_{s,\al}  \qquad \forall  t \in\RRR .
\]
%$u(\cdot,t)\in \calF^{-1}(\matU_\epsilon(\mathtt{g}(s,\al)))$
\end{enumerate}
\end{theorem}
%%%%%%%%%%%%%%%%%%%%%%%%%%%%%%%%%%%%%%%%%%%%%%%%%%%%%%%%%%%%%%%%%%%%%%%%%% 

	%%%%%%%%%%%%%%%%%%%%%%%%%%%%%%%%%%%%%%%%%%%%%%%%%%%%%%%%%%%%%%%%%%%%%%%%%% }
%\red{
%%%%%%%%%%%%%%%%%%%%%%%%%%%%%%%%%%%%%%%%%%%%%%%%%%%%%%%%%%%%%%%%%%%%%%%%%%% 
%\begin{coro}[Probabilistic control of Gevrey norms]\label{proposta}
%	Under the same assumptions as in Theorem \ref{main} the following holds:
%	 there exists  a full measure set $\mathfrak G\subset \mathfrak G_0$ such that
%	  for all $V\in \mathfrak G$  there exists $\rho_0>0$ such that for all $\rho<\rho_0$  and for any initial datum $W\in \matS_V\cap \matU_{\rho}(\mathtt{g}(s,\al))$
%	the solution $u(x,t)$ to \eqref{nls} is such that
%	$$\|\calF(u(\cdot ,t))\|_{s,\al} \le 2\|W\|_{s,\al}   $$
%	%$u(\cdot,t)\in \calF^{-1}(\matU_\epsilon(\mathtt{g}(s,\al)))$
%	 for all $t\in\RRR$.
%\end{coro}
%%%%%%%%%%%%%%%%%%%%%%%%%%%%%%%%%%%%%%%%%%%%%%%%%%%%%%%%%%%%%%%%%%%%%%%%%%% 
%}

%%%%%%%%%%%%%%%%%%%%%%%%%%%%%%%%%%%%%%%%%%%%%%%%%%%%%%%%%%%%%%%%%%%%%%%%%% 
\begin{rmk}\label{facaldo}
\emph{
Let $\g^*$ be as in Theorem \ref{proposta}. Set
\[
{\calG}:=\bigcup_{\gamma\in(0,\g^*)}\calG(\gamma).
\]
As a consequence of Theorem \ref{proposta},
for all $V\in\calG$ and all $\epsilon >0$  there exists $\de>0$ such that for any initial datum 
$W\in \matS_V\cap \matU_\de(\mathtt{g}(s,\al))$
the solution $u(x,t)$ to \eqref{nls} is such that
$u(\cdot,t)\in \calF^{-1}(\matU_\epsilon(\mathtt{g}(s,\al)))$ for all $t\in\RRR$.
This yields that the origin is statistically Lyapunov stable -- statistically
since the initial data which do not move away from the origin form a set of large but not full measure.
}
\end{rmk}
%%%%%%%%%%%%%%%%%%%%%%%%%%%%%%%%%%%%%%%%%%%%%%%%%%%%%%%%%%%%%%%%%%%%%%%%%% 

Theorem \ref{proposta} follows from the fact that the almost-periodic solutions in Theorem \ref{main} are such that
the ``linear actions'' $|u_j|^2$ (see \eqref{1.2bis} for the notation) are approximately constant for all times.
In fact, we give a more precise geometric description of the solutions and show that,
as it happens in classical (finite-dimensional) KAM theory,
also the support of any almost-periodic solution in Theorem \ref{main} is a torus,
which is a slight deformation of the ``flat torus'' where all the linear actions are constant.
However, since in our case most tori are infinite-dimensional, some care is needed in defining such objects.

%%%%%%%%%%%%%%%%%%%%%%%%%%%%%%%%%%%%%%%%%%%%%%%%%%%%%%%%%%%%%%%%%%%%%%%%%% 
\subsection{Structure of the invariant tori}\vspace{-.2cm}
%%%%%%%%%%%%%%%%%%%%%%%%%%%%%%%%%%%%%%%%%%%%%%%%%%%%%%%%%%%%%%%%%%%%%%%%%% 

For any $a>0$, the \emph{thickened torus}
\begin{equation} \nonumber %\label{torino}
\TTT^\ZZZ_a:=\Bigl\{\f \in \CCC^\ZZZ:\quad \mbox{Re}\f_j\in \TTT\,,\quad |\mbox{Im}\f_j |< a \,,\quad \forall j\in \ZZZ \Bigr\} ,
\end{equation}
endowed with the metric
\begin{equation} \nonumber
{\rm dist}(\varphi,\varphi'):= \sup_{j\in\ZZZ} \bigl(|\mbox{Re}(\varphi_j-\varphi'_j)\hbox{ mod }{2\pi}|+|\mbox{Im}(\varphi_j-\varphi'_j)|\bigl)\,,
\end{equation}
is a Banach manifold modelled on $\ell^{\io}(\RRR)$.
For any function $\TTT^\ZZZ_a\mapsto \mathtt{g}(s,\al)$ we can define the derivatives w.r.t.~the angles $\f_j$ in the natural way,
and characterize the analytic functions as the limits w.r.t.~the uniform topology of trigonometric polynomials depending on a finite 
number of angles \cite{MP,CGP2}. 
Then, in order to prove that \eqref{nls} admits an almost-periodic solution
with the properties stated in Theorem \ref{main},
actually we prove the following.
Fix $s>0$, $\al\in(0,1)$ and $N\ge 3$ (see Remark \ref{N3}),
and assume $\e$ to be small enough;
then there are a constant $a=a(\e,s,\al,N)$ and, for any $V\in \mathcal G$ and any $W\in \TT_V$,
a vector
{$\om=\om(V,W)\in\RRR^\ZZZ$}
such that the Cauchy problem \eqref{nls} admits a solution $u(x,t)$ of the form
\begin{equation}\label{unnumero}
u(x,t) = \UU(x,\om t) :=\sum_{j\in\ZZZ}e^{\ii j x} \uu_j(\om t) 
\end{equation}
and the map $\mathfrak{i}_{\,\UU} \!:\TTT^\ZZZ_a\to \mathtt{g}(s,\al)$, defined as
\begin{equation} \label{iU}
\mathfrak{i}_{\,\UU}(\f) := \{\uu_j(\f)\}_{j\in\ZZZ} , 
\end{equation}
is analytic. This implies that the solution $u(x,t)$ is almost-periodic in $t$ with frequency $\om$. 

However, the map $\mathfrak{i}_{\,\UU}$ in \eqref{iU} may well not be injective, so that both periodic and 
quasi-periodic solutions may -- and in fact are -- also obtained.
Furthermore, differently from the finite-dimensional case, the image $\mathfrak{i}_{\,\UU}(\TTT^\ZZZ)$ 
may not be a submanifold of $\gsa$.
Recall that,\footnote{
We follow the terminology used by P\"oschel and Trubowitz -- note, however, that sometimes in the literature
a subset with the listed properties is called an \emph{embedded submanifold}.}
given a Banach space $E$, a subset $M\subseteq E$ is a \emph{submanifold of $E$}
if for any $p\in M$ there is an open subset $U\subseteq E$ containing $p$ and
a homeomorphism $\phi \! :U\to V$, with $V$ an open subset of a Banach space $F$ such that
$F= F_h \oplus F_v$ and $\phi(U\cap M) = V \cap F_h$ \cite{potru}.
In other words, if in a neighbourhood of each point $p$ there is a coordinate system such that $M$ is ``locally horizontal''.
Nevertheless, although $\mathfrak{i}_{\,\UU}(\TTT^\ZZZ)$ may fail to be a submanifold for all $W\in \TT_V$,
if we define, for $\de\in(0,1/2)$,
\begin{subequations} \label{mesciarsi}
\begin{align}
\calA_{\de} (\mathtt g(s,\al)) & :=
\biggl\{\mathtt{x} \! \in \! \mathtt g(s,\al) \,:\,  \de\le |\mathtt{x}_j|e^{s\jap{j}^\al} \!\! \le \frac{1}{2} \; \forall\, j\in\ZZZ\biggr\}
\subset \ol\matU_{1/2}(\gsa) ,
\label{mesciarsia} \\
\matA_{V,\de} & := \TT_V\cap\calA_\de(\mathtt g(s,\al)) , \phantom{\ol{\sum}}
\end{align}
\end{subequations}
and take $W\in\matA_{V,\de}$, then the corresponding set $\mathfrak{i}_{\,\UU}(\TTT^\ZZZ)$ turns out to be a submanifold of $\gsa$.
In particular, the corresponding homeomorphism is analytic. The discussion above is summarised in the following result,
which is proved in Section \ref{outro} (recall Remark \ref{N3} as to the condition on $N$).

%Moreover, if we make explicit the dependence on $W$ of the solution $U(x,\om t)$, so as to write $\mathfrak{i}_U=\mathfrak{i}_{U,W}$,
%then there exists a family of infinite-dimensional tori
%\[
%\bigcup_{W \in \TT_{V,\de}}  \mathfrak{i}_{U,W}(\TTT^\ZZZ)
%\]
%and a Lipschitz embedding

%%%%%%%%%%%%%%%%%%%%%%%%%%%%%%%%%%%%%%%%%%%%%%%%%%%%%%%%%%%%%%%%%%%%%%%%%% 
\begin{theorem}\label{biascoappeso}
Fix $s>0$, $\al\in(0,1)$ and $N\ge3$.
For all $\e\in(-\e_*,\e_*)$, with $\e_*$ as in Theorem \ref{main}, there exists $a=a(\e,s,\al,N)$ such that,
for all $V\in\mathcal{G}$ and all $W\in \TT_V $, with $\mathcal{G}$ and $\TT_V$ as in Theorem \ref{main},
there exists an analytic map $\mathfrak{i}_{\,\UU}\!:\TTT^\ZZZ_a\to\gsa$ such that
the unique solution $u(x,t)$ to \eqref{nls} is of the form \eqref{unnumero}
and has the following property:
for all $\de>0$, there exists $\e_\star:=\e_{\star}(s,\al,N,\de)\le\e_*$ such that, for all $\e\in(-\e_\star,\e_\star)$,
all $V\in\mathcal{G}$ and all $W\in \matA_{V,\de}$, 
the set $\mathfrak{i}_{\,\UU}(\TTT^\ZZZ)$
is a submanifold of $\gsa$ analytically homeomorphic to $\TTT^\ZZZ$.
\end{theorem}
%%%%%%%%%%%%%%%%%%%%%%%%%%%%%%%%%%%%%%%%%%%%%%%%%%%%%%%%%%%%%%%%%%%%%%%%%% 

We conclude by discussing more in detail
the structure of the sets described by the hulls
of the almost-periodic solutions,\footnote{See footnote \ref{hull} for the definition of hull.}
as it emerges from the proof of Theorem \ref{biascoappeso} in Section \ref{outro}
(see in particular Remark \ref{crucial2}).
%Recallthat the hull of an almost-periodic function
%$f(t)$ is the closure w.r.t.~the uniform topology of the set $\cup_{\tau\in\RRR}f(t+\tau)$.
In order to prove Theorem \ref{biascoappeso} -- as well as Theorem\ref{main} --
we start by considering the NLS equation
\begin{equation} \label{NLS2}
\ii u_t - u_{xx} + \VV * u + \e |u|^4u =0 , \qquad x\in\TTT ,
\end{equation}
and ignoring for the time being the initial datum in \eqref{nls}.
For $\e=0$, all solutions to \eqref{NLS2} are of the form
\begin{equation} \label{marione}
u_{\rm lin}(x,t):= \sum_{j\in\ZZZ} c_j e^{\ii j x +\ii(j^2+V_j)t}\,,
\end{equation}
that is of the form \eqref{unnumero} with
$\uu_j(\f)=\uu_{{\rm lin},j}(\f):= c_j e^{\ii \f_j}$
and  $\om=\omega_{\rm lin}:=\{j^2+V_j\}_{j\in\ZZZ} $.
Then, we introduce the map
\begin{equation} \label{ilin}
\mathfrak{i}_{\rm lin}(\f) :=\{\uu_{{\rm lin},j}(\f)\}_{j\in\ZZZ} ,
\end{equation}
and set, for any $c\in\gsa$, any $I=\{I_j\}_{j\in\ZZZ}\in\RRR^\ZZZ_+\cap\mathtt{g}(2s,\al)$ and any $\theta\!\in\!\TTT^\ZZZ$,
\begin{equation} \label{c2}
|c|^2 := \{|c_j|^2\}_{j\in\ZZZ} , \qquad
\sqrt{I}:=\{\sqrt{I_j}\,\}_{j\in\ZZZ} , \qquad
c e^{\ii\theta}:=\{c_je^{\ii\theta_j}\}_{j\in\ZZZ}.
\end{equation}
If, for any choice of the \emph{actions} $I$, we define
\begin{equation} \label{TI}
\matT_I:= \bigl\{ c' \in \mathtt g(s,\al) : {|c'|^2= I} \bigr\} =
\bigl\{ c' \in \mathtt g(s,\al) : c' = \sqrt{I} e^{\ii\theta} , \, \theta\in\TTT^\ZZZ \bigr\} ,
\end{equation}
the set
\begin{equation} \nonumber %\label{toropiatto}
\mathfrak{i}_{\rm lin}(\TTT^\ZZZ) := \bigcup_{\f\in\TTT^\ZZZ} \mathfrak{i}_{\rm lin}(\f) = \matT_{|c|^2}, 
\end{equation}
is what is usually called a \emph{flat torus}.
This means that all solutions to the linear equations are supported on
invariant flat tori and we have the stratification
\begin{equation}\label{parafoli}
\gsa = \bigcup_{I\in\RRR^\ZZZ_+\cap\mathtt{g}(2s,\al)}\matT_I = \bigcup_{c\in\gsa} \matT_{|c|^2} .
\end{equation}
Finally, for each fixed {$c\in \gsa$} the linear dynamics on the torus described by the flow %Kronecker flow
$t \mapsto \varphi+ \omega_{\rm lin}t$ is dense on $\matT_{|c|^2}$ with respect to the {weak-$^*$ topology}, 
which turns out to be  the product  topology. 
Thus, by construction,
the image of the hull of the linear solution \eqref{marione} is the torus $\matT_{|c|^2}$ endowed with the product topology. 
It is however convenient to think of the tori  $\matT_{|c|^2}$ as immersions in $\gsa$.

%%%%%%%%%%%%%%%%%%%%%%%%%%%%%%%%%%%%%%%%%%%%%%%%%%%%%%%%%%%%%%%%%%%%%%%%%% 
\begin{rmk} \label{no-si-mani}
\emph{
In general $\matT_{|c|^2}$ is not a submanifold and hence \eqref{parafoli}
does not provide a foliation of $\mathtt g(s,\al)$, because some components of $c$ may be arbritrarily small \cite{potru}.
In order to have a foliation, we need a quantitative lower bound on $|c_j|$ for all $j\in\ZZZ$.
In particular, it is possible to extract from $\mathtt g(s,\al)$ a suitable open and dense subset $\IIIIO$
%\[
%\matP:=
%\bigcup_{\de>0}
%\bigcup_{\de\in(0,1/2)}
%\biggl\{ c  \in \! \mathtt g(s,\al) :  |c_j|e^{s\jap{j}^\al} \ge \de \; \forall\, j\in\ZZZ\biggr\} ,
%\calA_\de(\mathtt g(s,\al))
%\]
which is foliated in flat tori which are invariant for the linear dynamics (see Remark \ref{crucial2} for details).
}
\end{rmk}
%%%%%%%%%%%%%%%%%%%%%%%%%%%%%%%%%%%%%%%%%%%%%%%%%%%%%%%%%%%%%%%%%%%%%%%%%% 

For $\e\neq0$, it is therefore natural to look for a solution to \eqref{NLS2} of the form \eqref{unnumero},
which is identified by a pair $(\mathfrak{i}_{\,\UU},\omega)$, 
with $\mathfrak{i}_{\,\UU}\!:\TTT^\ZZZ_a\to \mathtt{g}(s,\al)$ and $\omega\in \RRR^\ZZZ$,
which continues the unperturbed solution
$(\mathfrak{i}_{\rm lin},\omega_{\rm lin})$, 
%$(\{c_j e^{\ii \f_j}\}_{j\in\ZZZ}, \{j^2+V_j\}_{j\in\ZZZ})$ 
and satisfies
\begin{equation} \label{NLS3}
\ii \omega\cdot\partial_\f \, \uu_j + (j^2 +V_j) \, \uu_j + \e (|\UU|^4  \UU)_j=0 ,  \qquad j \in \ZZZ ,
\end{equation}
where\footnote{Here and henceforth, $\ol{z}$ denotes the complex conjugate of $z\in\CCC$.}
\begin{equation} \label{quintic}
\left( |\UU|^4 \UU \right)_{j} = \!\!\!\!\!
\sum_{\substack{j_1,j_2,j_3,j_4 \in \ZZZ \\ j_1-j_2+j_3 - j_4 +j_5 = j } }
\uu_{j_1} \ol{\uu}_{j_2} \uu_{j_3} \ol{\uu}_{j_4} \uu_{j_5}\,.
\end{equation}

In analogy with the classical (finite-dimensional) KAM theory, one expects that
for all $V\in\mathcal{G}$ there is a large measure Cantor set
$\IIII \subseteq \RRR^\ZZZ_+\cap\ol{\matU}_{1/4}(\mathtt{g}(2s,\al))$, %for some $\rho_0 \in (0,1)$,}
such that for all $I\in\IIII$ the flat torus $\matT_I$ survives the perturbation slightly deformed.
We prove that there is a family of flat tori
\begin{equation} \label{TC}
\matT_{\IIII} := \bigcup_{I\in\IIII} \matT_I ,
\end{equation}
a ball $\ol{\matU}_{\rho_0}(\gsa)$ for some $\rho_0\in(1/2,1)$ and a Lipschitz map $\gotWW\!\!:\ol{\matU}_{\rho_0}(\gsa) \to \gsa$
such that the image of the restriction of $\gotWW$ to each $\matT_I$ is an invariant torus,
that is a set homotopically equivalent to $\matT_I$ and invariant for the nonlinear dynamics.
Such tori, being the image under the map $\gotWW$ of a flat torus, are in principle
%In fact, any invariant torus for the nonlinear system is the image under the map $\gotWW$ of a flat torus of the linear system,
%so that \emph{a fortiori} it is 
not necessarily submanifolds (see Remark \ref{no-si-mani}). In any case, the dynamics on $\gotWW(\matT_I)$ is
conjugated to a linear dynamics on $\matT_I$ and hence, as in the case $\e=0$,
the flow is dense on $\gotWW(\matT_I)$ w.r.t.~the weak-$^*$ topology.

%%%%%%%%%%%%%%%%%%%%%%%%%%%%%%%%%%%%%%%%%%%%%%%%%%%%%%%%%%%%%%%%%%%%%%%%%% 
\begin{rmk} \label{si-no-mani}
\emph{
The low regularity of the map $\gotWW$ is due to the fact that, whereas, writing $c=\sqrt{I}e^{\ii \f}$,
for fixed $I$ the dependence of the map on the angles $\f$ is analytic, on contrast the actions $I$ are restricted to
a Cantor set, and the extension of $\gotWW$ outside such a set is no more than Lipschitz.
}
\end{rmk}
%%%%%%%%%%%%%%%%%%%%%%%%%%%%%%%%%%%%%%%%%%%%%%%%%%%%%%%%%%%%%%%%%%%%%%%%%% 

If we introduce the set
\begin{equation} \label{C0}
\IIIIC := \Bigl\{ c \in\IIIIO :  |c|^2 \in\IIII \Bigr\} , 
\end{equation}
then, for any $c\in\IIIIC$, the set $\gotWW(\matT_{|c|^2})$
turns out to be a submanifold (in fact, an \emph{embedded torus}) and, according to Kuksin and P\"oschel terminology \cite{kupo}, the set
\begin{equation} \label{inutile}
\matE := \bigcup_{c\in\IIIIC} \gotWW(\matT_{|c|^2} )
\end{equation}
is a \emph{Cantor manifold} foliated in tori which are invariant for \eqref{NLS2}.
As a map from $\matE$ to $\gsa$, the embedding $\gotWW$ is only a lipeomorphism, by Remark \ref{si-no-mani}.

%%%%%%%%%%%%%%%%%%%%%%%%%%%%%%%%%%%%%%%%%%%%%%%%%%%%%%%%%%%%%%%%%%%%%%%%%% 
\begin{rmk} \label{remarkalpostodiabuse}
\emph{
Of course, each initial datum $W\in \matE:=\gotWW(\matT_\IIII) \subset \gsa$ gives rise to an almost-periodic solution.
The crucial difference w.r.t.~the finite dimensional case is that $\IIII$ having large measure does not guarantee that
$\matE$ has large measure as well, and in principle the set $\matE$ might not even be measurable.
This is the reason for formulating Theorem \ref{main}, Corollary \ref{facaldo} and Theorem \ref{biascoappeso} in terms of
the initial data rather than the actions.
}
\end{rmk}
%%%%%%%%%%%%%%%%%%%%%%%%%%%%%%%%%%%%%%%%%%%%%%%%%%%%%%%%%%%%%%%%%%%%%%%%%% 

%%%%%%%%%%%%%%%%%%%%%%%%%%%%%%%%%%%%%%%%%%%%%%%%%%%%%%%%%%%%%%%%%%%%%%%%%% 
%%%%%%%%%%%%%%%%%%%%%%%%%%%%%%%%%%%%%%%%%%%%%%%%%%%%%%%%%%%%%%%%%%%%%%%%%% 
\zerarcounters
\section{A quantitative Moser counterterm theorem}\label{secmoser}
%%%%%%%%%%%%%%%%%%%%%%%%%%%%%%%%%%%%%%%%%%%%%%%%%%%%%%%%%%%%%%%%%%%%%%%%%% 
%%%%%%%%%%%%%%%%%%%%%%%%%%%%%%%%%%%%%%%%%%%%%%%%%%%%%%%%%%%%%%%%%%%%%%%%%% 

{As anticipated in Section \ref{intro}, it is more convenient to reformulate \eqref{1.2} as a counterterm problem.
We take $\om$ in the set
\begin{equation}\label{rettangolo}
\gQ := \Big\{ \om\in \RRR^\ZZZ \; :\; |\om_j-j^2| \le 1/2 \Big\} 
\end{equation}
and, instead of the original equation \eqref{1.2}, we study the equation \eqref{1.3}, with $\eta$ a parameter sequence to be suitably fixed.
If we set $\UUU_j(\om t):=u_j(t)$ and define 
\begin{equation} \label{unnumerobis}
\UUU(x,\f) :=\sum_{j\in\ZZZ}e^{\ii j x} U_j(\f) ,
\vspace{-.2cm}
\end{equation}
then \eqref{1.3} becomes
\begin{equation} \label{NLS4}
\ii \omega\cdot\partial_\f \uuu_j + (\omega_j+\eta_j) \uuu_j+ \e (| \UUU|^4 \UUU )_j=0 , \qquad j\in\ZZZ ,
\end{equation}
where $(| \UUU|^4 \UUU )_j$ is meant as in \eqref{quintic}. We call \eqref{NLS4} the \emph{modified equation}
and the sequence $\eta$ the \emph{counterterm} -- by borrowing the terminology used in Quantum Field Theory.}
Then the unknown for \eqref{NLS4} becomes the pair $(\UUU,\h)$,
with $U\!:\TTT \times \TTT^\ZZZ_a\to \calF^{-1}(\mathtt{g}(s,\al))$ and $\h\in \ell^\io(\RRR)$.

%%%%%%%%%%%%%%%%%%%%%%%%%%%%%%%%%%%%%%%%%%%%%%%%%%%%%%%%%%%%%%%%%%%%%%%%%% 
\begin{rmk} \label{Ueta}
\emph{
The modified equation \eqref{NLS4} reduces to \eqref{NLS3} if
\begin{equation} \label{compatibility}
\om_j+\eta_j=j^2+V_j, \qquad j\in\ZZZ .
\end{equation}
In other words, if $(U,\h)$ solves \eqref{NLS4} and $\h=\{\h_j\}_{j\in\ZZZ}$ is such
that \eqref{compatibility} is satisfied, then
the function $\UU(x,\om t)=\UUU(x,\om t)$ solves \eqref{NLS2}.
{We call \eqref{compatibility} the \emph{compatibility condition}.}
}
\end{rmk}
%%%%%%%%%%%%%%%%%%%%%%%%%%%%%%%%%%%%%%%%%%%%%%%%%%%%%%%%%%%%%%%%%%%%%%%%%% 

%%%%%%%%%%%%%%%%%%%%%%%%%%%%%%%%%%%%%%%%%%%%%%%%%%%%%%%%%%%%%%%%%%%%%%%%%% 
\subsection{The set of frequencies}\vspace{-.2cm}
%%%%%%%%%%%%%%%%%%%%%%%%%%%%%%%%%%%%%%%%%%%%%%%%%%%%%%%%%%%%%%%%%%%%%%%%%% 

Then, at first, we neglect the constraint \eqref{compatibility}.
In order to prove that there exists a solution $(\UUU,\h)$ to \eqref{NLS4}
%such that the bounds
%\[
%\sup_{\f\in \TTT^\ZZZ_a}\| U_j(\f)-c_j e^{\ii \f_j}\|_{s,\al} \le C_0 |\e| , \qquad \|\eta\|_{\io} \le C_0 |\e|
%\]
%are verified for a suitable constant $C_0$ and for all $c\in \ol{\matU}_1(\gsa)$ and all 
%satisfying
we need to require $\om\in \gQ$ to satisfy an appropriate Bryuno non-resonance condition, as detailed below.
Set\footnote{\label{11}Here and henceforth, $\|\cdot\|_1$ is the standard $\ell^1$-norm and the subscript $f$ stands for \emph{finite support}.}
\begin{subequations} \label{interi}
\begin{align}
\ZZZ^\ZZZ_f & := \Bigl\{\nu\in\ZZZ^\ZZZ : \|\nu\|_1<\io \Bigr\}, \phantom{\biggr\}}
\label{interia} \\
\ZZZ^\ZZZ_{f,0} & :={\biggl\{ \nu\in\ZZZ^\ZZZ_f :  \sum\limits_{i\in\ZZZ} \nu_i=0 \biggr\},}
\label{interib}
\end{align}
\end{subequations}
and, for $\al>0$, define
\begin{equation}\label{normalfa}
|\nu|_{\al}:=\sum_{i\in\ZZZ} \jap{i}^\al |\nu_i|,
\end{equation}
so that $\|\cdot\|_1=|\cdot|_0$. For any fixed $\al\in(0,1)$ and $\om\in\gQ$, introduce the non-increasing function 
$\be^{(0)}_\om \! : [1,+\io) \to\RRR$ given by
\begin{equation}\label{beta}
{
\beta^{(0)}_\om(x):= \inf \bigl\{ |\om\cdot\nu| : \nu\in\ZZZ^{\ZZZ}_{f,0} , \; 0 <|\nu|_{\al/2}\le x %, \;  \sum\limits_{i\in\ZZZ} \nu_i=0 
\bigr\} ,}
% \inf_{\substack{\nu\in\ZZZ^{\ZZZ}_{f} \vspace{-.1cm}\\ 0 <|\nu|_{\al/2}\le x\\ \sum\limits_{i\in\ZZZ} \nu_i=0}}  |\om\cdot\nu|,
\end{equation}
the \emph{weak Bryuno function}
\begin{equation}\label{juno}
\BB(\om) := \sum_{m\ge1}\frac{1}{2^m}\log \left(\frac{1}{\be^{(0)}_\om(2^m)}\right) ,
\end{equation}
and define
\begin{equation}\label{set}
\gotB^{(0)} := \{\om\in{\gQ}\;:\; \BB(\om)<\io\}.
\end{equation}
We say that $\om\in\gotB^{(0)}$ is a (infinite-dimensional) vector satisfying the \emph{weak Bryuno condition}.

%%%%%%%%%%%%%%%%%%%%%%%%%%%%%%%%%%%%%%%%%%%%%%%%%%%%%%%%%%%%%%%%%%%%%%%%%% 
\begin{rmk} \label{nuovissimormk2}
\emph{
One could replace the sequence $\{2^m\}_{m\in\NNN}$ with an arbitrary diverging increasing sequence
$\{r_m\}_{m\ge0}$ on $[1,+\io)$, with no significant change in the following discussion \cite{CGP}.
}
\end{rmk}
%%%%%%%%%%%%%%%%%%%%%%%%%%%%%%%%%%%%%%%%%%%%%%%%%%%%%%%%%%%%%%%%%%%%%%%%%% 

%%%%%%%%%%%%%%%%%%%%%%%%%%%%%%%%%%%%%%%%%%%%%%%%%%%%%%%%%%%%%%%%%%%%%%%%%% 
\begin{rmk} \label{nuovissimormk}
\emph{
The reason why the function $\BB$ in \eqref{juno} is called the weak Bryuno function
and the vectors $\om\in\gotB^{(0)}$ are said to satisfy the weak Bryuno condition
is that only vectors $\nu\in\ZZZ^\ZZZ_f$
such that
\begin{equation} \label{sumnui0}
\sum_{i\in\ZZZ} \nu_i=0 
\end{equation}
are considered when computing $\beta^{(0)}_\om(x)$ in \eqref{beta}.
If we drop the constraint \eqref{sumnui0} and allow vectors such that \eqref{juno}, with $\beta^{(0)}_\om(x)$ being replaced with 
\[
\beta_\om(x)  := \inf_{\substack{\nu\in\ZZZ^{\ZZZ}_{f} \vspace{-.1cm}\\ 0 <|\nu|_{\al/2}\le x }}  |\om\cdot\nu| ,
\]
is finite, we recover the Bryuno vectors which naturally extend to the infinite-dimensional context \cite{CGP}
the usual Bryuno vectors of the finite-dimensional case \cite{Bry}.
}
\end{rmk}
%%%%%%%%%%%%%%%%%%%%%%%%%%%%%%%%%%%%%%%%%%%%%%%%%%%%%%%%%%%%%%%%%%%%%%%%%% 

%%%%%%%%%%%%%%%%%%%%%%%%%%%%%%%%%%%%%%%%%%%%%%%%%%%%%%%%%%%%%%%%%%%%%%%%%% 
\begin{rmk} \label{nuovormk}
\emph{
The set of \emph{$(\gamma,\tau)$-weak Diophantine vectors}, i.e.~the set
\begin{equation} \nonumber %\label{diofantinoBIS}
\gD^{(0)}(\g,\tau) :=\biggl\{ \omega\in {\gQ} : |\omega\cdot \nu|> \g
\prod_{i\in \ZZZ}\frac{1}{(1+ \jap{i}^{2}|\nu_i|^{2})^\tau} \quad \forall \nu\in \ZZZ^{\ZZZ}_{f,0} \setminus\!\{0\} %, \quad \sum_{i\in\ZZZ} \nu_i=0 
\biggr\} ,
\end{equation}	
is a subset of $\gotB^{(0)}$ \cite[Lemma B.3]{CGP}
and, with the notation in Remark \ref{pesiemisure}, for $\tau>1/2$
there exists a constant $C(\tau)$ such that $\mu_0(\gQ\setminus \gD^{(0)}(\g,\tau)) \le C(\tau)\,\g$
(see \cite{Bjfa,yuan} or Lemma B.5 in \cite{CGP}).
}
\end{rmk}
%%%%%%%%%%%%%%%%%%%%%%%%%%%%%%%%%%%%%%%%%%%%%%%%%%%%%%%%%%%%%%%%%%%%%%%%%% 

%%%%%%%%%%%%%%%%%%%%%%%%%%%%%%%%%%%%%%%%%%%%%%%%%%%%%%%%%%%%%%%%%%%%%%%%%% 
\subsection{Almost-periodic solutions to the modified equation}\vspace{-.2cm}
%%%%%%%%%%%%%%%%%%%%%%%%%%%%%%%%%%%%%%%%%%%%%%%%%%%%%%%%%%%%%%%%%%%%%%%%%% 

Then, we look for a solution $(U,\h)$, with $U(x,\f)=U(x,\f;c,\om,\e)$, to the modified equation \eqref{NLS4} of the form
\begin{equation}\label{perp}
\begin{aligned}
\UUU(x,\f;c,\om,\e) \! & \;= 
\UU_0(x,\f ;c)
+ \UUU_\perp (x,\f; c, \om,\e) , \phantom{\sum_{\nu\in\ZZZ}} \\
\UU_0(x,\f;c)
& := \sum_{j\in\ZZZ} c_j e^{\ii j x+\ii \f_j} , \\
\UUU_\perp (x,\f; c, \om,\e) & := \sum_{j\in\ZZZ} \sum_{\substack{\nu\in\ZZZ^{\ZZZ}_{f} \vspace{-.1cm} \\ \nu\ne \gote_j}} \uuu_{j,\nu}(c,\om,\e) \, e^{\ii jx + \ii\nu\cdot\f} ,
\end{aligned}
\end{equation}
with $c\in \mathtt g({s,\alpha})$ and $\om\in\gotB^{(0)}$, and the coefficients $\uuu_{j,\nu}(c,\om,\e)$
of the function $\UUU_\perp(x,\f; c, \om,\e)$ to be determined.

It is convenient to rewrite \eqref{perp} as
\[
\begin{aligned}
\UUU(x,\f;c,\om,\e) :=\sum_{j\in\ZZZ} \UUU_j(\f;c,\om,\e)e^{\ii j x}= \sum_{j\in\ZZZ} \sum_{\nu\in\ZZZ^{\ZZZ}_{f}} \uuu_{j,\nu} (c,\om,\e)e^{\ii j x+\nu\cdot\f} ,
\end{aligned}
\]
where we have set $\UUU_{j,\gote_j}(c,\om,\e)=c_j$,
if $\gote_j\in\ZZZ^\ZZZ_f$ denotes the integer vector with components $(\gote_j)_i=\de_{ij}$, with $\de_{ij}$ being the Kroneker delta.
%for all $j\in\ZZZ$,
%\[
%U_{j,\nu}(c,\om,\e)= \begin{cases}
%c_j , & \qquad \nu=\gote_j , \\
%u_{j,\nu}(c,\om,\e), & \qquad \nu\neq\gote_j .
%\end{cases}
%\]

Then, %assuming that the coefficients $\uuu_{j,\nu}=\uuu_{j,\nu}(c,\om,\e)$ have an appropriate decay and writing
the modified equation \eqref{NLS4} becomes
\begin{equation} \label{nls-modified}
\begin{array}{rl}
\left( -\om\cdot\nu + \om_j + \eta_j \right) \uuu_{j,\nu} + \e \left( |\UUU|^4 \UUU \right)_{j,\nu} = 0 ,  & \qquad j \in \ZZZ, \; \nu \neq \gote_j , \\
\eta_j  c_j + \e \left( |\UUU|^4 \UUU \right)_{j,\gote_j} = 0 , & \qquad j \in \ZZZ , \phantom{\displaystyle{\int}}
\end{array}
\end{equation}
where we have shortened $\h_j=\h_j(c,\om,\e)$ and $\UUU_{j,\nu}=\UUU_{j,\nu}(c,\om,\e)$, and set
\[
\left( |\UUU|^4 \UUU \right)_{j,\nu} =
\sum_{\substack{j_1,j_2,j_3,j_4,j_5\in\ZZZ \\ j_1-j_2+j_3 - j_4 +j_5 = j}}
\sum_{\substack{\nu_1,\nu_2,\nu_3,\nu_4,\nu_5\in\ZZZ^\ZZZ_f  \\ \nu_1 -\nu_2 + \nu_3-\nu_4 +\nu_5 =\nu }}
\uuu_{j_1,\nu_1} \ol{\uuu}_{j_2,\nu_2} \uuu_{j_3,\nu_3} \ol{\uuu}_{j_4,\nu_4} \uuu_{j_5,\nu_5} .
\]
The following result rephrases Theorem 2.13 in \cite{CGP}.

%%%%%%%%%%%%%%%%%%%%%%%%%%%%%%%%%%%%%%%%%%%%%%%%%%%%%%%%%%%%%%%%%%%%%%%%%% 
\begin{theorem} \label{moser}
Fix $s>0$, $\al\in(0,1)$ and $\om\in\gotB^{(0)}$. There exists  $\e_0=\e_0(s,\al,\om)$ such that
for all $\e\in(-\e_0,\e_0)$ and all $c\in \ol{\matU}_1(\mathtt{g}(s,\al))$ there exist two sequences
\[
\{ \uuu_{j,\nu} (c,\om,\e)\}_{(j,\nu) \in \ZZZ \times \ZZZ_f^\ZZZ } , 
\qquad
\h(c,\om,\e) = \{\h_j(c,\om,\e)\}_{j\in \ZZZ} , 
\]
satisfying the following properties:
\vspace{-.2cm}
\begin{enumerate}
\itemsep0em
\item the corresponding $(\UUU,\h)$ solves the modified equation \eqref{nls-modified},
\item there exist $s_0'=s_0'(\e) \in (0,s)$ and, for all $s_2\in(s_0',s)$, a constant $C=C(s,s_2,\al,\om)$ such that
\vspace{-.2cm}
\begin{equation} \label{moser-bounds}
\sup_{j\in\ZZZ} \left| \h_j(c,\om,\e) \right|  < C |\e| , 
\qquad
\sup_{j\in\ZZZ} 
\sup_{\substack{\nu\in\ZZZ^{\ZZZ}_{f} \setminus\{\gote_j}\}} 
\left| \uuu_{j,\nu}(c,\om,\e) \right| e^{s_1 |\nu|_\al}
e^{s_2\jap{j}^\al} \le C|\e| ,
\vspace{-.2cm}
\end{equation}
with $s_1:=s-s_2$,
\item the sequence $\eta(c,\om,\e)$ is real analytic in $|c|^2:=\{|c_j|^2\}_{j\in\ZZZ}$, 
\item the function
$U_\perp (x,\f; c, \om,\e)$  in \eqref{perp}-- and hence $U(x,\f;c,\om,\e)$ -- is separately analytic in $c,\bar c$,
\item
there exists $a=\aaa(\e,s,\al,\om)$ such that, setting
\begin{equation} \label{Uj}
\mathfrak{i}_{\UUU}(\f):=\{\uuu_j(\f;c,\om,\e)\}_{j\in\ZZZ} , \qquad
\UUU_j(\f;c,\om,\e) := \sum_{{\nu\in\ZZZ^\ZZZ_f }} \uuu_{j,\nu}(c,\om,\e) e^{\ii\nu\cdot \varphi},
\vspace{-.2cm}
\end{equation}
the map $\mathfrak{i}_{\UUU} \!: \TTT^\ZZZ_{a} \to \gsa$ is analytic.
%\vspace{-.2cm}
\end{enumerate}
\end{theorem}
%%%%%%%%%%%%%%%%%%%%%%%%%%%%%%%%%%%%%%%%%%%%%%%%%%%%%%%%%%%%%%%%%%%%%%%%%% 

%%%%%%%%%%%%%%%%%%%%%%%%%%%%%%%%%%%%%%%%%%%%%%%%%%%%%%%%%%%%%%%%%%%%%%%%%% 
\begin{rmk} \label{separate}
\emph{
In Theorem \ref{moser}, separate analyticity is meant as follows \cite[Definition 2.11]{CGP}.
Given two complex Banach spaces $(Z,\| \cdot \|_Z)$ and $(Y, \|\cdot \|_Y)$,
a function $F\!:Z\to Y$ is said to be \emph{separately analytic} in $z,\bar z$ if there is an analytic function
${F}_{\rm ext}\!:Z\times Z\to Y$ such that $F(z)={F}_{\rm ext}(z,\ol{z})$ for all $z\in Z$.
}
\end{rmk}
%%%%%%%%%%%%%%%%%%%%%%%%%%%%%%%%%%%%%%%%%%%%%%%%%%%%%%%%%%%%%%%%%%%%%%%%%% 

%%%%%%%%%%%%%%%%%%%%%%%%%%%%%%%%%%%%%%%%%%%%%%%%%%%%%%%%%%%%%%%%%%%%%%%%%% 
%\begin{rmk} \label{separateb}
%\emph{
%\blue{
%The dependence of $s_0'(\e)$ on $\e$ is such that we can say, equivalently, that for all $s_2\in(0,s)$ there exists $\e(s_2)$ small enough
%such that the bounds \eqref{moser-bounds} hold provided we have $|\e|\le \e(s_2)$ (see \cite{CGP} for details).
%In particular we can fix $s_2\ge s/2$ and choose $\e(s_2)$ accordingly. Note that we have always $s_1,s_2\in(0,s)$.
%(mi sembrava che servisse dopo, se no si pu\`o togliere)}
%}
%\end{rmk}
%%%%%%%%%%%%%%%%%%%%%%%%%%%%%%%%%%%%%%%%%%%%%%%%%%%%%%%%%%%%%%%%%%%%%%%%%% 

Both $\h(c,\om,\e)$  and $\UUU_\perp(x,\f;c,\om,\e)$ may be explicitly computed as absolutely convergent series in $\e$, i.e.~for 
all $c\in \overline{\matU}_1(\mathtt{g}(s,\al))$ and all $\om\in\gotB^{(0)}$ and for $\e$ small enough, we may expand
\begin{subequations} \label{formale+fout}
\begin{align}
\UUU_\perp(x,\f;c,\om,\e) & \!=\! 
\sum_{k\ge1}\e^k \sum_{j\in\ZZZ} \! \UUU_{j}^{(k)}(\f;c,\om) \, e^{\ii j x} , \qquad
% = \sum_{k\ge1}\e^k \sum_{j\in\ZZZ} \sum_{\substack{\nu\in\ZZZ^{\ZZZ}_{f} \vspace{-.1cm} \\ \nu \neq \gote_j}}
% \uuu_{j,\nu}^{(k)}(c,\om) \, e^{\ii (j x+\nu\cdot\f)} ,
\UUU_{j}^{(k)}(\f;c,\om) \!=\! \!\! \sum_{\substack{\nu\in\ZZZ^{\ZZZ}_{f} \vspace{-.1cm} \\ \nu \neq \gote_j}}
\!\! \uuu_{j,\nu}^{(k)}(c,\om) \, e^{\ii \nu\cdot\f} ,
\label{fout} \\
\h_j (c,\om,\e) & = \sum_{k\ge1} \e^k \h_j^{(k)} (c,\om) 
\label{formale}
\end{align}
\end{subequations}
where the coefficients $\uuu_{j,\nu}^{(k)}(c,\om)$ and $\h_j^{(k)} (c,\om)$ admit a graphical representation in terms of trees \cite[Section 4]{CGP}.

The bounds in item 2 of Theorem \ref{moser} are obtained by studying the coefficients $\uuu_{j,\nu}^{(k)}(c,\om)$ and $\h_{j,\nu}^{(k)}(c,\om)$
in \eqref{formale+fout} and proving that, for any $k\ge 1$, any $j\in\ZZZ$ and any $\nu\in\ZZZ^\ZZZ_f\setminus\{\gote_j\}$, one has
\cite[Proposition 7.83]{CGP}
\begin{equation} \label{moser-bounds-k}
\bigl| \h_j^{(k)}(c,\om) \bigr|  < D_0^k , 
\qquad
\bigl| \uuu^{(k)}_{j,\nu}(c,\om) \bigr| e^{s_1 |\nu|_\al}
e^{s_2\jap{j}^\al} \le D_0^k  ,
\end{equation}
for some positive constant $D_0=D_0(s,\al,\om)$, which goes to infinity as $s$ goes to $0$.
{By introducing the bounds \eqref{moser-bounds-k} in \eqref{fout} and
recalling the notation \eqref{c2} and the definitions \eqref{perp} and \ref{Uj}, we obtain the bounds
\begin{equation} \label{Ujbound}
\left\| \mathfrak{i}_{\UUU} (\f) - c e^{\ii \f} \right\|_{s,\al} \le D_1 |\e| ,
\end{equation}
for some positive constant $D_1=D_1(s,\al,\om)$.}

%%%%%%%%%%%%%%%%%%%%%%%%%%%%%%%%%%%%%%%%%%%%%%%%%%%%%%%%%%%%%%%%%%%%%%%%%% 
\begin{rmk}\label{troppo}
\emph{
By looking at formula $(4.27)$ of \cite{CGP},
{it is easy to prove inductively}
that  $\uuu^{(k)}_{j,\nu}(c,\om)$ is a homogenous polynomial in 
$c,\bar c$ of degree $4k+1$, while $\h_j^{(k)}(c,\om)$ is a homogeneous polynomial in $|c|^2$ of degree $2k$.
}
\end{rmk}
\subsection{Properties of the almost-periodic solutions}\vspace{-.2cm}
%%%%%%%%%%%%%%%%%%%%%%%%%%%%%%%%%%%%%%%%%%%%%%%%%%%%%%%%%%%%%%%%%%%%%%%%%% 

From the symmetry property of the modified equation \eqref{NLS4} we may deduce that
if a solution of the form \eqref{unnumerobis} exists, then it satisfies the following crucial symmetries:
\begin{enumerate}[topsep=0ex]
\itemsep0em
\item \emph{Momentum conservation}: for any $\theta\!\in\!\TTT^\ZZZ$, %if we set $c e^{\ii\theta}:=\{c_je^{\ii\theta_j}\}_{j\in\ZZZ}$,  
we have \cite[Remark 2.17]{CGP}
\begin{equation}\label{chinottoneri}
\UUU(x,\f;ce^{\ii\theta},\om,\e)  = \UUU(x,\f + \theta;c,\om,\e) ,\qquad
\h(c,\om,\e)=\h(ce^{\ii\theta},\om,\e).
\end{equation}
\item \emph{Translation-covariance}:
for any $x\in\RRR$, we have \cite[Remark 2.17]{CGP}
\begin{equation}\label{chinottonerib}
\UUU(x,\f;c,\om,\e)  = \UUU(0,\f + \varsigma x ;c,\om,\e) = U(0,\f ;c e^{\ii \varsigma x} ,\om,\e) ,
\end{equation}
{with $\varsigma = \{j \}_{j\in \ZZZ} $ (see Remark \ref{pesiemisure})},
i.e.~the solutions we obtain are travelling waves.
\item \emph{Gauge-covariance}: for any $\la\in\RRR$ we have \cite[Remark 2.9]{CGP}
\begin{equation}\label{chinottofaschifo}
\UUU(x,\f;c,\om,\e) =e^{-\ii \la} \UUU(x,\f;ce^{\ii\la},\om,\e) .
\end{equation}
\item{\emph{Scaling-covariance}}: for any $\rho>0$ we have \cite[Remark 2.18]{CGP}
\begin{equation}\label{moltoschifo}
\UUU(x,\f;c,\om,\e) =\rho^{-1} \UUU(x,\f;\rho c,\om,\rho^{-4}\e) ,\qquad
\h(c,\om,\e) = \h(\rho c,\om,\rho^{-4}\e).
\end{equation}
\end{enumerate}

%%%%%%%%%%%%%%%%%%%%%%%%%%%%%%%%%%%%%%%%%%%%%%%%%%%%%%%%%%%%%%%%%%%%%%%%%% 
\begin{rmk}\label{verifica}
\emph{
The symmetries listed above can also be verified \emph{a posteriori} from the
explicit construction of the solution to the modified equation \eqref{NLS4} in terms of trees \cite{CGP}.
The scaling covariance is the symmetry referred to in Remark \ref{riscalo}.
}
\end{rmk}
%%%%%%%%%%%%%%%%%%%%%%%%%%%%%%%%%%%%%%%%%%%%%%%%%%%%%%%%%%%%%%%%%%%%%%%%%% 

%%%%%%%%%%%%%%%%%%%%%%%%%%%%%%%%%%%%%%%%%%%%%%%%%%%%%%%%%%%%%%%%%%%%%%%%%% 
\begin{rmk}\label{cervellospento}
\emph{
By using \eqref{chinottoneri}, we immediately deduce that, if $c_j=0$ for some $j\in\ZZZ$, then $\del_{\f_j}\UUU(x,\f;c,\om,\e)=0$.
}
\end{rmk}
%%%%%%%%%%%%%%%%%%%%%%%%%%%%%%%%%%%%%%%%%%%%%%%%%%%%%%%%%%%%%%%%%%%%%%%%%% 

%%%%%%%%%%%%%%%%%%%%%%%%%%%%%%%%%%%%%%%%%%%%%%%%%%%%%%%%%%%%%%%%%%%%%%%%%% 
\subsection{Asymptotic expansion of the counterterms}
%%%%%%%%%%%%%%%%%%%%%%%%%%%%%%%%%%%%%%%%%%%%%%%%%%%%%%%%%%%%%%%%%%%%%%%%%% 

If we assume the frequencies to have a suitable special dependence on $j$,
then the counterterm turns out to have the same dependence as well. 
More precisely, the subset of frequencies in $\gQ$ we consider is as follows. 
Define
\begin{equation} \nonumber %\label{bienne}
\matW_N:= [-1/4,1/4]^{N-1}\times \ol{\matU}_{\!{1/2}}(\ell^{N,\io}(\RRR)) ,
\end{equation}
{and introduce in $\matW_N$ the norm
\begin{equation} \label{stanormetta}
\|\ze\|_{m,\io} := \max\{ \|\kappa\|_{\io} , \|\x\|_{m,\io} \} .
\end{equation}
For}
any $\ze\in\matW_N$, write $\ze=(\ka,\xi)$, with $\ka=(\ka_0,\ka_2,\ldots,\ka_{N-1})\in[-1/4,1/4]^{N-1}$
and $\xi\in \ol{\matU}_{\!{1/2}}(\ell^{N,\io}(\RRR))$.
For all $j\in\ZZZ$, define also 
\begin{equation}\label{espome}
\om_j(\ze) := \om_j(\ka,\xi)=
\begin{cases}
\ka_0+\xi_0,&\qquad j=0, \\
\displaystyle{ j^2 +\kappa_0+ \sum_{q=2}^{N-1} \frac{\kappa_q}{j^q} +\xi_j} \, ,& \qquad j \neq 0 ,
\end{cases}
\end{equation}
and set
\begin{equation} \nonumber %\label{KN}
\matK_N:=\{ \ze\in \matW_N:\,\omega(\ze)\in \gotB^{(0)}\} ,
\end{equation}
with $\gotB^{(0)}$ defined in \eqref{set}. Then, the following result holds \cite[Proposition 2.21]{CGP}.
	
%%%%%%%%%%%%%%%%%%%%%%%%%%%%%%%%%%%%%%%%%%%%%%%%%%%%%%%%%%%%%%%%%%%%%%%%%% 
\begin{prop} \label{prop2.21}
Fix $s>0$, $\al\in(0,1)$ and $N\ge 3$. Let $\e_0$ be as in Theorem \ref{moser}.
For any $\ze\in \matK_N$ there exists ${\e_1=\e_1(s,\al,\ze,N)\in(0,\e_0)}$ such that,
for all $\e\in(- \e_1,\e_1)$  and all $c\in\ol{\matU}_1(\mathtt{g}(s,\al))$,
there are a vector $\gota(c,\ze,\e)= (\gota_0(c,\ze,\e),\gota_2(c,\ze,\e),\ldots,\gota_{N-1}(c,\ze,\e))\in\RRR^{N-1}$
and a sequence $\gotr(c,\ze,\e)= \{\gotr_j(c,\ze,\e)\}_{j\in\ZZZ} \in \ell^{N,\io}(\RRR)$, both
real analytic in $|c|^2$, such that
%depending on $c$ only through $|c|^2$, such that 
%
\begin{equation} \label{nemmenounnome!}
\h_j(c,\om(\ze),\e) = 
\begin{cases}
\displaystyle{\gota_0(c,\ze,\e)+ \sum_{q=2}^{N-1} \frac{\gota_q(c,\ze,\e)}{j^q} + {\gotr_j(c,\ze,\e)},}& \qquad j \neq 0 ,\\
\gota_0(c,\ze,\e) + \gotr_0(c,\ze,\e), \phantom{\displaystyle{\int}} & \qquad j=0 . \\
\end{cases}
\end{equation}
Moreover, there is a positive constant
{$C_0=C_0(s,\al,\ze,N)$,}
independent of both $\e$ and $c$, such that
\begin{equation} \nonumber %\label{tosse}
%\sup_{q=0,2,\ldots,N-1}|\gota_q(c,\ze,\e)| \le C_0|\e|,
\|\gota(c,\ze,\e)\|_{\io} \le C_0|\e|,
\qquad \|\gotr(c,\ze,\e)\|_{N,\io}\le C_0| \e| .
\end{equation}
Finally,
{$\gota(c,\ze,\e)$ and  $\gotr(c,\ze,\e)$} are analytic in $\e$ and their Taylor expansion satisfies the same properties as $\eta$.
{In particular,}
$\gota^{(k)}$ and $\gotr^{(k)}$ are homogeneous polynomials in $|c|^2$ of degree $2k$.
\end{prop}
%%%%%%%%%%%%%%%%%%%%%%%%%%%%%%%%%%%%%%%%%%%%%%%%%%%%%%%%%%%%%%%%%%%%%%%%%% 

%%%%%%%%%%%%%%%%%%%%%%%%%%%%%%%%%%%%%%%%%%%%%%%%%%%%%%%%%%%%%%%%%%%%%%%%%% 
%%%%%%%%%%%%%%%%%%%%%%%%%%%%%%%%%%%%%%%%%%%%%%%%%%%%%%%%%%%%%%%%%%%%%%%%%% 
\zerarcounters
\section{Lipschitz continuity and uniform bounds}
\label{uni}
%%%%%%%%%%%%%%%%%%%%%%%%%%%%%%%%%%%%%%%%%%%%%%%%%%%%%%%%%%%%%%%%%%%%%%%%%% 
%%%%%%%%%%%%%%%%%%%%%%%%%%%%%%%%%%%%%%%%%%%%%%%%%%%%%%%%%%%%%%%%%%%%%%%%%% 

Theorem \ref{moser} is an abstract Moser-like counterterm theorem which ensures the existence of
a solution $(U,\eta)$ to the modified equation \eqref{NLS4},
%the solution $U(x,\f)$ and of the counterterm $\eta$
for all $c\in \matU_1(\mathtt{g}(s,\alpha))$ and all $\om\in\gotB^{(0)}$.
The radius of convergence $\e_0$ is uniform in $c$, for $c\in\ol{\matU}_1(\mathtt{g}(s,\al))$,
but depends strongly on the frequency $\om$.
The same happens for the radius of convergence $\e_1$ as a function of $\zeta$ in Proposition \ref{prop2.21}.
However, in order to apply the result to the NLS equation \eqref{nls}, we need
to  solve the compatibility equation \eqref{compatibility} and express both the frequency
$\om$ and the counterterm $\h$ in terms of the potential $V$. This is an implicit
function problem and, if we aim to solve it, we need Lipschitz regularity
on some open set of the frequencies where one has uniform bounds on the radius of convergence in $\e$.

To study the Lipschitz continuity of the functions we are interested in, it is useful to introduce some notation.
If $X$ and $Y$ are normed spaces endowed with the norms $\|\cdot \|_X$ and $\|\cdot\|_Y$, respectively,
given any non-empty $D\subseteq X$ and any function $f\!:D\to Y$, we define the seminorm
\begin{equation}\label{brrrr2}
\|f\|_{{\rm Lip}(D\subseteq X,Y)} := \sup_{\substack{\mathtt x,\mathtt x'\in D \\ \mathtt x\neq \mathtt x'}}
\frac{\|f( \mathtt x)-f( \mathtt x')\|_Y}{\| \mathtt x- \mathtt x' \|_X} .
\end{equation}
Below, we use repeatedly the notation \eqref{brrrr2}, contextualized to the normed space we consider each time.

%%%%%%%%%%%%%%%%%%%%%%%%%%%%%%%%%%%%%%%%%%%%%%%%%%%%%%%%%%%%%%%%%%%%%%%%%% 
\subsection{The set of good parameters}\vspace{-.2cm}
%%%%%%%%%%%%%%%%%%%%%%%%%%%%%%%%%%%%%%%%%%%%%%%%%%%%%%%%%%%%%%%%%%%%%%%%%% 

In order to restrict $\ze$ to a  subset $\matK_N(\g)\subseteq\matK_N$
in which the radius of convergence $\e_1$ is bounded uniformly in $\ze$, we proceed as follows.
Fix $\g>0$ and $\tau>1/2$ and set
\begin{equation} \label{notte}
	 \be^*(x,\g,\tau) := \g \!\!\!\! 
	\inf_{\substack{\nu\in\ZZZ^{\ZZZ}_{f} \vspace{-.1cm}\\ 0 <|\nu|_{\al/2}\le x \vspace{-.1cm}\\ {\sum_{i\in\ZZZ} \nu_i = 0}}} 
	%\inf_{\substack{\nu\in\ZZZ^{\ZZZ}_{f} \vspace{-.1cm}\\ 0 <|\nu|_{\al/2}\le x}}
	\prod_{i\in\ZZZ} \frac{1}{(1+\jap{i}^2|\nu_i|^2)^\tau} .
	%\qquad \BB(r,\g) := \sum_{m\ge1}\frac{1}{r_{m-1}}\log\frac{1}{\be^*(r_m,\g)} , 
\end{equation}
One can verify \cite[Lemma B.3]{CGP} that
\begin{equation} \label{adessoce}
	\BB(\g,\tau) := \sum_{m\ge1}\frac{1}{2^{m}}\log\frac{1}{\be^*(2^m,\g,\tau)}<\io.
\end{equation}

Define
\begin{equation} \label{Bbeta1}
\gotB(\gamma,\tau)   = \bigl\{ \om\in\gotB^{(0)} : \be_\om^{(0)}(2^m) \ge \be^*(2^m,\g,\tau) \quad \forall m\ge0 \bigr\} .
\end{equation}
%
%%%%%%%%%%%%%%%%%%%%%%%%%%%%%%%%%%%%%%%%%%%%%%%%%%%%%%%%%%%%%%%%%%%%%%%%%% 
\begin{rmk}\label{diofa}
\emph{
For any $\g,\tau>0$ one has $\gD^{(0)}(\g,\tau)\subseteq \gotB(\gamma,\tau)$ \cite[Appendix B]{CGP}.
}
\end{rmk}
%%%%%%%%%%%%%%%%%%%%%%%%%%%%%%%%%%%%%%%%%%%%%%%%%%%%%%%%%%%%%%%%%%%%%%%%%% 

%%%%%%%%%%%%%%%%%%%%%%%%%%%%%%%%%%%%%%%%%%%%%%%%%%%%%%%%%%%%%%%%%%%%%%%%%% 
\begin{rmk} \label{fuoridalnulla}
\emph{
In the following we need a stronger condition on $\tau$ than the condition $\tau>1/2$ considered above.
In fact, we have to require at least  $\tau=N+1$
(see Section \ref{misura}, in particular Proposition \ref{44} and Remark \ref{tau=N+1}).
Hence, from now on, we fix $\tau=N+1$.
}
\end{rmk}
%%%%%%%%%%%%%%%%%%%%%%%%%%%%%%%%%%%%%%%%%%%%%%%%%%%%%%%%%%%%%%%%%%%%%%%%%% 

%For $N\ge3$ (see Remark \ref{N3}), %and $\tau=N+1$. %and $\be=\{\be_m\}_{m\in\NNN}$ such that \eqref{pedanteria} holds. Then,
The appropriate subset $\matK_N(\g)$ is defined as
\begin{equation} \label{bgamma} 
\matK_N(\g)  :=\{ \ze \in \matW_N : \omega(\ze) \in \gotB(\gamma,N+1) \}  .
%(\gamma,N+1)
\end{equation}
Following \cite{CGP}, we call $\matK_N(\g)$ the set of \emph{good parameters}.

%%%%%%%%%%%%%%%%%%%%%%%%%%%%%%%%%%%%%%%%%%%%%%%%%%%%%%%%%%%%%%%%%%%%%%%%%% 
\begin{rmk}\label{piugen}
\emph{
We defined the sequence $\{\beta^*(2^m,\gamma,\tau)\}_{m\ge0}$ so that the following holds:
\begin{enumerate}[topsep=0ex]
\itemsep0em
\item the vectors in $\gotB(\gamma,\tau)$ satisfy the weak Bryuno condition,
\item the set $\gotB(\gamma,\tau)$ contains the set of $(\gamma,\tau)$-weak Diophantine vectors (see Remark \ref{nuovormk})
\item for $\zeta\in\matK_N(\gamma)$ we have a uniform bound on the small divisors.
\end{enumerate}
In fact,
instead of the sequence $\{\beta^*(2^m,\gamma,N+1)\}_{m\ge0}$,
we may consider a different sequence as soon as the properties above are satisfied.
Precisely, given any sequence $\be=\{\be_m\}_{m\in\NNN}$ such that 
\begin{subequations} \label{pedanteria}
\begin{align}
& \be^*(2^m,\g,N+1)  \ge \be_m,\qquad m\ge0, 
\label{pedanteriaa} \\
& \BB_\be := \sum_{m\ge1}\frac{1}{2^{m-1}}\log\frac{1}{\be_m} < \io ,
\label{pedanteriab}
\end{align}
\end{subequations}
if we define the set
\begin{equation} \label{Bbeta2}
\gotB_\be   = \bigl\{ \om\in\gotB^{(0)} : \be_\om^{(0)}(2^m) \ge \be_m \quad \forall m\ge0 \bigr\} ,
\end{equation}
by construction we have
$\gotB(\g,N+1) \subset \gotB_\be$.
Then, instead of the set \eqref{bgamma}, we may consider the larger set 
\[
\matK_{N,\be}  :=\{ \ze \in \matW_N : \omega(\ze) \in \gotB_\be \}  ,
%(\gamma,N+1)
\]
and this allows us to consider almost-periodic solutions with frequencies
which are not necessarily Diophantine and still admit uniform estimates~\cite[Proposition 3]{CGP2}.
However, here we prefer to use the sets $\gotB(\g,N+1)$ and $\matK_N(\gamma)$ in \eqref{bgamma}
to keep the same notation as in \cite{CGP}.
}
\end{rmk}
%%%%%%%%%%%%%%%%%%%%%%%%%%%%%%%%%%%%%%%%%%%%%%%%%%%%%%%%%%%%%%%%%%%%%%%%%% 

%%%%%%%%%%%%%%%%%%%%%%%%%%%%%%%%%%%%%%%%%%%%%%%%%%%%%%%%%%%%%%%%%%%%%%%%%% 
%\begin{defi}[\textbf{Good parameters}]
%\label{unibr}
%Let $\gotB_\be$ be defined in \eqref{Bbeta}.
%Let $\gotB(\g,\tau)  \subset\gotB^{(0)}$ denote the set of weak Bryuno vectors $\om$ with the property that
%$\be^{(0)}_\om(2^m)\ge \be^*(2^m,\g,\tau)$ for all $m\ge0$.
%For $N\ge 0$, 
%\end{defi}
%%%%%%%%%%%%%%%%%%%%%%%%%%%%%%%%%%%%%%%%%%%%%%%%%%%%%%%%%%%%%%%%%%%%%%%%%% 

%%%%%%%%%%%%%%%%%%%%%%%%%%%%%%%%%%%%%%%%%%%%%%%%%%%%%%%%%%%%%%%%%%%%%%%%%% 
\subsection{Lipschitz extensions}\vspace{-.2cm}
%%%%%%%%%%%%%%%%%%%%%%%%%%%%%%%%%%%%%%%%%%%%%%%%%%%%%%%%%%%%%%%%%%%%%%%%%% 

Given any subset $\matU\subseteq\matW_N$, with $N\ge 3$, and any map $f\!:\matU\to E$, 
with $(E,\|\cdot\|_E)$ some Banach space, we define the family of Lipschitz norms of $f$ as
\begin{equation}\label{lipnorm}
|f|_{\matU,m, E}^{{\rm{Lip}}} := \sup_{\ze\in \matU} \| f(\ze) \|_E +  \|f\|_{{\rm Lip}(\matU \subseteq \matW_N, E)} ,
%+ \sup_{\substack{ \ze,\ze'\in \matU \\ \ze\neq\ze'}} \frac{\|f(\ze) -f(\ze') \|_E}{\|\ze-\ze'\|_{m,\io}}
\qquad m=0,\ldots,N . 
\end{equation}
%
%\blue{For future convenience, we write, for $m\ne 0$,
%
%\begin{equation} \nonumber %\label{torta}
%|f|_{m,E}^{{\rm{Lip}}} := |f|_{\matW_N,m,E}^{{\rm{Lip}}} ,
%\end{equation}
%}

%%%%%%%%%%%%%%%%%%%%%%%%%%%%%%%%%%%%%%%%%%%%%%%%%%%%%%%%%%%%%%%%%%%%%%%%%%% 
%\begin{rmk} \label{serve?}
%\emph{
%\blue{
%The norm $|\cdot|_{m,E}^{{\rm{Lip}}}$ in \eqref{torta} is decreasing in $m$.}
%}
%\end{rmk}
%%%%%%%%%%%%%%%%%%%%%%%%%%%%%%%%%%%%%%%%%%%%%%%%%%%%%%%%%%%%%%%%%%%%%%%%%%% 

Define
\begin{equation}\label{ugot}
\gotu_{j}(c,\ze,\e)
:= \sum_{k\ge1}\e^k \UUU_{j}^{(k)}(0;c,\om(\ze)) ,
\qquad 
\gotu(c,\ze,\e)
:=\{\gotu_j(c,\ze,\e)\}_{j\in\ZZZ} , 
\end{equation}
so that we can write
\[
U(x,0) = \sum_{j\in\ZZZ} \bigl( c_j + \gotu_j(c,\ze,\e) \bigr) e^{ijx} ,
\]
and requiring that $U(x,0)=\WW(x)$, that is the initial datum in \eqref{nls}, reads
\begin{equation} \label{WcU}
W_j=c_j + \gotu_j(c,\ze,\e) , \qquad j\in\ZZZ .
\end{equation}
 
The vector $\gota(c,\ze,\e)$ and the sequence $\gotr(c,\ze,\e)$ appearing in \eqref{nemmenounnome!},
as well as the sequence $\gotu(c,\ze,\e)$ in \eqref{ugot}, satisfy suitable estimates which allow us
to extend them to the whole $\matW_N$.
This is the content of the two following results.

%%%%%%%%%%%%%%%%%%%%%%%%%%%%%%%%%%%%%%%%%%%%%%%%%%%%%%%%%%%%%%%%%%%%%%%%%% 
\begin{prop} \label{fabrizio}
Fix $N\ge 3$ and $\g\in(0,1)$.
For any $s>0$ and any $\al\in(0,1)$ %and all $s_1 \ge 0$ and $s_2>0$ such that $s_1+s_2=s$
there exist $\e_2=\e_2(s,\al,N,\g)>0$ and a positive constant $C_1=C_1(s,\al,N,\g)$ such that,
for all  $\e\in(-\e_2, \e_2)$ and all $c\in \ol{\matU}_1(\mathtt{g}(s,\al))$, one has
\begin{subequations} \label{dolore!}
\begin{align}
|\h(c,\om(\cdot),\e)|^{{\rm{Lip}}}_{\matK_N(\g),0, \ell^\io(\RRR)} & \le C_1 |\e| , 
\label{dolore!a} \\
%|\Uperp(c,\om(\cdot),\e) |^{{\rm{Lip}}}_{\matK_N(\g),\mathtt{W}(s_1,s_2,\al)}  & \le C |\e| , 
%\label{dolore!b} \\
%|\gota_q(c,\cdot,\e)|^{{\rm{Lip}}}_{\matK_N(\g),0,\RRR} & \le C_1 |\e|  , \qquad q=0,2,\ldots,N-1, 
|\gota(c,\cdot,\e)|^{{\rm{Lip}}}_{\matK_N(\g),0,\RRR^{N-1}} & \le C_1 |\e|  , 
\label{muffa1} \\
|\gotr(c,\cdot,\e)|^{{\rm{Lip}}}_{\matK_N(\g),N,\ell^{N,\io}(\RRR)} & \le C_1 |\e|  ,
\label{muffa2} 
\\
|\gotr(c,\cdot,\e)|^{{\rm{Lip}}}_{\matK_N(\g),0,\ell^{2,\io}(\RRR)}  & \le C_1 |\e|  ,
\label{muffa3} \\
|\gotu(c,\cdot,\e)|^{{\rm{Lip}}}_{\matK_N(\g),0,{\mathtt g}(s,\al)}  & \le C_1 |\e| .
\label{dolorissimo}
\end{align}
\end{subequations}
\end{prop}
%%%%%%%%%%%%%%%%%%%%%%%%%%%%%%%%%%%%%%%%%%%%%%%%%%%%%%%%%%%%%%%%%%%%%%%%%% 

%%%%%%%%%%%%%%%%%%%%%%%%%%%%%%%%%%%%%%%%%%%%%%%%%%%%%%%%%%%%%%%%%%%%%%%%%% 
\proof
The bounds \eqref{dolore!a} to \eqref{muffa3} are proved in \cite{CGP} (see Proposition 2.27 therein).
The same argument leads to the proof of the last bound \eqref{dolorissimo} as well.
Indeed, both $\gota(c,\ze,\e)$ and $\gotr(c,\ze,\e)$ are limits of quantities which can be written as sums over trees \cite[Section 8]{CGP},
and the bounds are obtained by estimating the value of each tree separately; analogously, 
each $\gotu_{j}(c,\ze,\e)$ admits a tree representation which can be dealt with in the same way \cite[Lemma 9.4]{CGP},
so that, by reasoning as in the proof of Lemma A.2 in \cite{CGP}, we obtain that, for some positive constant $D_1$, 
\begin{equation}\label{dolore!b}
\sup_{\f\in \TTT_a} |\UUU_{j}^{(k)}(\f;c,\om(\cdot)) |^{{\rm{Lip}}}_{\DgN,0,\RRR}   \le  D_1^k e^{-s\jap{j}^\al},
\end{equation}
with $\UUU_{j}^{(k)}(\f;c,\om)$ defined in \eqref{fout}, and hence
\begin{equation} \nonumber %\label{forsesi}
\sum_{k\ge1}|\e|^k\sup_{\f\in \TTT_a}\sup_{j\in\ZZZ} | \UUU_{j}^{(k)}(\f;c,\om(\cdot)) |^{{\rm{Lip}}}_{\DgN,0,\RRR} e^{s\jap{j}^\al}  \le D_2 |\e| ,
\end{equation}
for a suitable positive constant $D_2$,
provided $\e$ is small enough. This implies that $\gotu(c,\ze,\e)\in \matU_\rho (\mathtt g(s,\al))$ for some $\rho=O(\e)$
and hence the bound \eqref{dolorissimo} follows,
possibly with a different constant $C_1$ with respect the other bounds. 
Thus, by suitably redefining the constant $C_1$, all the bounds are found to hold with the same value of $C_1$.
\qed
%%%%%%%%%%%%%%%%%%%%%%%%%%%%%%%%%%%%%%%%%%%%%%%%%%%%%%%%%%%%%%%%%%%%%%%%%% 

%%%%%%%%%%%%%%%%%%%%%%%%%%%%%%%%%%%%%%%%%%%%%%%%%%%%%%%%%%%%%%%%%%%%%%%%%% 
\begin{prop} \label{piergiorgio}
Fix $N\ge 3$ and $\g\in(0,1)$. For $s>0$ and $\al\in(0,1)$,
let $\e_2$ be as in Proposition \ref{fabrizio}.
Then, for all  $\e\in(-\e_2, \e_2)$ and all $c\in \ol{\matU}_1(\mathtt{g}(s,\al))$, there exists Lipschitz extensions
\[
\begin{aligned}
A(c,\ze,\e) & = (A_0(c,\ze,\e), A_2(c,\ze,\e),\ldots,A_{N-1}(c,\ze,\e)) , \\
R(c,\ze,\e) & = \{ R_j(c,\ze,\e) \}_{j\in\ZZZ} , \\
\gotU(c,\ze,\e) 
& = \{ \gotU_j(c,\ze,\e) \}_{j\in\ZZZ}  
\end{aligned}
\]
of the functions
$\gota(c,\ze,\e)$, $\gotr(c,\ze,\e)$, $\gotu(c,\ze,\e)$, respectively,  to the whole $\matW_N$,
continuous w.r.t.~the product topology and satisfying the same bounds as in \eqref{dolore!}, with $\matW_N$ replacing $\matK_N(\g)$.
In particular, the counterterm \eqref{nemmenounnome!} can be extended to the function
\begin{equation} \label{nemmenounnome!ext}
\h_j(c,\om(\ze),\e) = 
\begin{cases}
\displaystyle{A_0(c,\ze,\e)+ \sum_{q=2}^{N-1} \frac{A_q(c,\ze,\e)}{j^q} + {R_j(c,\ze,\e)},}& \qquad j \neq 0 , \\
A_0(c,\ze,\e) + R_0(c,\ze,\e), \phantom{\displaystyle{\int}} & \qquad j=0 , \\
\end{cases}
\end{equation}
defined for all $\ze\in\matW_N$.
\end{prop}
%%%%%%%%%%%%%%%%%%%%%%%%%%%%%%%%%%%%%%%%%%%%%%%%%%%%%%%%%%%%%%%%%%%%%%%%%% 

%%%%%%%%%%%%%%%%%%%%%%%%%%%%%%%%%%%%%%%%%%%%%%%%%%%%%%%%%%%%%%%%%%%%%%%%%% 
\proof
The extensions of the functions $\gota(c,\ze,\e)$ and $\gotr(c,\ze,\e)$ 
are proved to exist and satisfy the bounds \eqref{dolore!a} to \eqref{muffa3} in
%ref.~\cite[Corollary 9.8 and formula (9.29)]{CGP}.
\cite{CGP} (see Corollary 10.8 and formula (10.29) therein).
Similarly, applying the McShane Theorem also to the function $u^{(k)}_{j,\nu}(c,\om(\cdot))$ and reasoning in the same way,
we find that the function $\gotu(c,\ze,\e)$ can be extended as well to a function $\gotU(c,\ze,\e)$ which is defined on the whole $\matW_N$,
where it admits the same bound as $\gotu(c,\ze,\e)$ in $\matK_N(\g)$.
\qed
%%%%%%%%%%%%%%%%%%%%%%%%%%%%%%%%%%%%%%%%%%%%%%%%%%%%%%%%%%%%%%%%%%%%%%%%%% 

%%%%%%%%%%%%%%%%%%%%%%%%%%%%%%%%%%%%%%%%%%%%%%%%%%%%%%%%%%%%%%%%%%%%%%%%%% 
\begin{rmk} \label{extensionbytrees}
\emph{
In the proof given in \cite{CGP} that the functions $\gota(c,\ze,\e)$ and $\gotr(c,\ze,\e)$ can be extended,
one of the main features is that
each component $\gota_q(c,\ze,\e)$, for $q=0,2,\ldots,N-1$ and each component $\gotr_j(c,\ze,\e)$, for $j\in\ZZZ$, 
is expressed as an absolutely convergent series where the summands are graphically represented as renormalized trees \cite[Section 7]{CGP},
and their values equal the values of the trees, which in turn are functions of a finite number of variables.
Then, all the components are extended by extending each single tree value. This allows us 
to obtain extensions $A(c,\ze,\e)$ and $R(c,\ze,\e)$, which, by construction, are continuous w.r.t.~the product topology in all their variables
and Lipschitz-continuous in the variables $\ze$. The same comment applies to the extension $\gotU(c,\ze,\e)$ of the function $\gotu(c,\ze,\e)$.
}
\end{rmk}
%%%%%%%%%%%%%%%%%%%%%%%%%%%%%%%%%%%%%%%%%%%%%%%%%%%%%%%%%%%%%%%%%%%%%%%%%% 

Finally, the following result  is  proved by relying once more on \cite{CGP}.

%%%%%%%%%%%%%%%%%%%%%%%%%%%%%%%%%%%%%%%%%%%%%%%%%%%%%%%%%%%%%%%%%%%%%%%%%% 
\begin{prop} \label{nuovaprop}
For all $(c,\ze,\e)\in \matU_{1}(\mathtt g{(s,\alpha)})\times \matW_N\times (-{\e_2},{\e_2})$
the extended functions $A(c,\ze,\e)$, $R(c,\ze,\e)$ and $\gotU(c,\ze,\e)$ in Proposition \ref{piergiorgio}
are separately analytic in $c,\bar c\in \ol{\matU}_1(\mathtt g(s,\al))$.
In particular, both $A(c,\ze,\e)$ and $R(c,\ze,\e)$ are real analytic in $|c|^2$.
\end{prop}
%%%%%%%%%%%%%%%%%%%%%%%%%%%%%%%%%%%%%%%%%%%%%%%%%%%%%%%%%%%%%%%%%%%%%%%%%% 

%%%%%%%%%%%%%%%%%%%%%%%%%%%%%%%%%%%%%%%%%%%%%%%%%%%%%%%%%%%%%%%%%%%%%%%%%% 
\prova
Using the notation in \cite{CGP} (see Appendix A.2, specifically formula (A.2) therein) and setting
 \[
 \ZZZ^{\ZZZ}_{+,f}:=\Bigl\{ a \in \ZZZ^\ZZZ_f \,:\, a_j\ge0\ \forall\,j\in\ZZZ \Bigr\}, \qquad j(a,b):= \sum_{i\in\ZZZ} i(a_i-b_i) ,
 \]
we may write
\begin{equation} \nonumber %\label{ecchep}
\uuu_{j,\nu}(c,\om(\ze),\e) = \!\!\!\!\!\!\!\! \sum_{ \substack{a,b\in \ZZZ^\ZZZ_{+,f}\\ a-b = \nu , \, j(a,b)=j}} \!\!\!\!\!\!\!\!\!
\uuu^{(a,b)}_{j,\nu}(\om(\ze),\e) \, c^a \bar c^b ,
\qquad  c^a \bar c^b:=\prod_{i\in\ZZZ} c_i^{a_i}\bar c_i^{b_i}\,,
\end{equation}
where the coefficients $\uuu^{(a,b)}_{j,\nu}(\om(\ze),\e)$ satisfy the bounds \cite[formula (A.3)]{CGP}
\begin{equation} \label{A3}
\sup_{\substack{j\in\ZZZ \\ \nu\in\ZZZ^\io_f}} e^{s_1|\nu|_\al} e^{(s-s_1)\jap{j}^\al}
\!\!\!\!\!\!\!\!\!
\sum_{ \substack{a,b\in \ZZZ^\ZZZ_{+,f}\\ a-b = \nu , \, j(a,b)=j}} \!\!\!\!\!\!\!\!\!
|\uuu^{(a,b)}_{j,\nu}(\om(\ze),\e)|
\sup_{c,\ol{c}\in \matU_1(\mathtt{g}(s,\al))}|c^a \bar c^b |<\io .
\end{equation}
Therefore, we obtain
\begin{equation} \nonumber %\label{addiritturina}
\gotu_j(c,\ze,\e) = \sum_{\nu\in\ZZZ^\ZZZ_f\setminus\{\gote_j\}}
\sum_{ \substack{a,b\in \ZZZ^\ZZZ_{+,f}\\ a-b = \nu , \, j(a,b)=j}} 
\uuu^{(a,b)}_{j,\nu}(\om(\ze),\e) \, c^a \bar c^b
\end{equation}
and hence the sequence $\gotu(c,\ze,\e)$ is
%normally
separately analytic in $c,\ol{c}$ for $c\in\matU_1(\mathtt g(s,\al))$.
A similar result holds also for $\gota(c,\ze,\e)$ and $\gotr(c,\ze,\e)$, namely we may expand
\begin{subequations}\label{addirittura}
\begin{align}
\gota_q(c,\ze,\e) & = \sum_{a\in\ZZZ^\ZZZ_{+,f}} \gota^{(a)}_q(\ze,\e)|c|^{2a}, \qquad q=0,2,\ldots,N-1 ,
\label{addiritturaa} \\
\gotr_j(c,\ze,\e) & = \sum_{a\in\ZZZ^\ZZZ_{+,f}} \gotr^{(a)}_j(\ze,\e)|c|^{2a} , \qquad j \in \ZZZ ,
\label{addiritturab}
\end{align}
\end{subequations}
with (see \eqref{c2} for the notation $|c|^2$)
\[
|c|^{2a}:= \prod_{i\in\ZZZ} |c_i|^{2a_i}
\]
and the coefficients $\gota^{(a)}(\ze,\e)$ and $\gotr^{(a)}(\ze,\e)$ satisfying bounds analogous to \eqref{A3}.
This shows that also the functions $\gota(c,\ze,\e)$ and $\gotr(c,\ze,\e)$ 
are separately analytic in $c,\bar c$.

Since the extension procedure is performed for each single tree value contributing to the coefficients $u_{j,\nu}^{(a,b)}(\om(\ze),\e)$,
$\gota^{(a)}_q(\ze,\e)$ and $\gotr_j^{(a)}(\ze,\e)$ (see Remark \ref{extensionbytrees}),
we conclude that the extended functions $A(c,\ze,\e)$, $R(c,\ze,\e)$ and $\gotU(c,\ze,\e)$
inherit the separate analyticity property of the functions $\gota(c,\ze,\e)$, $\gotr(c,\ze,\e)$ and $\gotu(c,\ze,\e)$.
\qed
%%%%%%%%%%%%%%%%%%%%%%%%%%%%%%%%%%%%%%%%%%%%%%%%%%%%%%%%%%%%%%%%%%%%%%%%%% 

%%%%%%%%%%%%%%%%%%%%%%%%%%%%%%%%%%%%%%%%%%%%%%%%%%%%%%%%%%%%%%%%%%%%%%%%%% 
\begin{rmk}\label{ancorascale}
\emph{
By using  Remark \ref{troppo} and Proposition \ref{prop2.21}, if we restrict
{all the involevd functions}
to $c\in\ol\matU_\rho(\gsa)$,
{with $\rho\in(0,1]$,}
then the estimates \eqref{dolore!a} to \eqref{muffa3} hold with $C_1$ replaced with $C_1\rho^4$ 
and the estimate \eqref{dolorissimo} holds with $C_1$ replaced with $C_1\rho^5$.
The same argument applies to the extensions $A,R,\gotU$.
}
\end{rmk}
%%%%%%%%%%%%%%%%%%%%%%%%%%%%%%%%%%%%%%%%%%%%%%%%%%%%%%%%%%%%%%%%%%%%%%%%%% 

%%%%%%%%%%%%%%%%%%%%%%%%%%%%%%%%%%%%%%%%%%%%%%%%%%%%%%%%%%%%%%%%%%%%%%%%%% 
Let us now discuss more in detail the regularity with respect  to $c$.
%%%%%%%%%%%%%%%%%%%%%%%%%%%%%%%%%%%%%%%%%%%%%%%%%%%%%%%%%%%%%%%%%%%%%%%%%% 

%%%%%%%%%%%%%%%%%%%%%%%%%%%%%%%%%%%%%%%%%%%%%%%%%%%%%%%%%%%%%%%%%%%%%%%%%% 
\begin{coro} \label{lip-c}
For any $\rho_0\in(0,1)$ and for all $(c,\ze,\e)\in
{\ol{\matU}_{\rho_0}(\mathtt g{(s,\alpha)})}
\times \matW_N\times (-{\e_2},{\e_2})$,
the extended functions $A(c,\ze,\e)$, $R(c,\ze,\e)$ and $\gotU(c,\ze,\e)$ in Proposition \ref{piergiorgio}
are Lipschitz continuous in
{$c \in \ol{\matU}_{\rho_0}(\mathtt g(s,\al))$}.
\end{coro}
%%%%%%%%%%%%%%%%%%%%%%%%%%%%%%%%%%%%%%%%%%%%%%%%%%%%%%%%%%%%%%%%%%%%%%%%%% 

%%%%%%%%%%%%%%%%%%%%%%%%%%%%%%%%%%%%%%%%%%%%%%%%%%%%%%%%%%%%%%%%%%%%%%%%%% 
\proof
Lipschitz continuity in $c$ is ensured by the separate analyticity in $c,\bar c$.
Indeed, the bounds on the Lipschitz norms in the variables $c$ follow from the Cauchy estimates: if
$(E,\|\cdot\|_E)$ is a Banach space and
 a function $f\!:\ol{\matU}_{\!1}(\mathtt{g}(s,\al))\to E$ is separately analytic in $c,\bar c$, then,
for a suitable constant $C$ depending on $\rho_0$,
{we have \cite{Muj}}
\begin{equation} \label{constantC}
{
|f|_{{\rm Lip}(\ol{\matU}_{\rho_0}(\mathtt g(s,\al)) \subseteq \mathtt g(s,\al),E)} } 
= \!\!\!
\sup_{\substack{c, c'\in \ol{\matU}_{\!\rho_0}(\mathtt g(s,\al)) \\ c \neq c'}}
\frac{\|f(c) -f(c') \|_E}{\|c- c'\|_{s,\al}}\le C \sup_{{c\in \ol{\matU}_{\!1}(\mathtt g(s,\al)) }}\|f(c)\|_E , 
\end{equation}
which immediately implies that $f$ is Lipschitz continuous in $c$.
\qed
%%%%%%%%%%%%%%%%%%%%%%%%%%%%%%%%%%%%%%%%%%%%%%%%%%%%%%%%%%%%%%%%%%%%%%%%%% 

%%%%%%%%%%%%%%%%%%%%%%%%%%%%%%%%%%%%%%%%%%%%%%%%%%%%%%%%%%%%%%%%%%%%%%%%%% 
\begin{rmk} \label{2/3}
\emph{
In the following we require $\rho_0$ to be greater than $1/2$,
for instance $\rho_0=2/3$ (which gives $C=3$ in \eqref{constantC}).
In fact, the dependence on $c$ of the extended function is much more than Lipschitz; however, we do not need more regularity than that.
}
\end{rmk}
%%%%%%%%%%%%%%%%%%%%%%%%%%%%%%%%%%%%%%%%%%%%%%%%%%%%%%%%%%%%%%%%%%%%%%%%%% 

Since the functions $A(c,\ze,\e)$, $R(c,\ze,\e)$ and $\gotU(c,\ze,\e)$ are Lipschitz continuous in $\ze \in \matW_N$ as well \cite{CGP},
it is convenient is to introduce a suitable Lipschitz norm in order to take into account both variables 
$c \in \ol{\matU}_{\!\rho_0}(\mathtt g(s,\al))$ and $\ze\in \matW_N$.

To this aim,
%recalling \eqref{bienne},
we set
\begin{equation}\label{asilo}
\matY_0:= \ol{\matU}_{\!\!\rho_0}(\mathtt g(s,\al))\times\matW_N%[-1/4,1/4]^{N-1}\times \ol{\matU}_{\!1/2}(\ell^{N,\io})
\end{equation}
and, for any $m=0,\ldots,N$, in $\matY_0$ we introduce the norm 
\begin{equation}\label{stanormazza}
\norm (c,\ze) \norm_{m}:= \max \bigl\{\|c\|_{s,\al},
%\|\kappa\|_{\io},
\|\ze\|_{m,\io} \bigr\} ,
\end{equation}
with $\|\ze\|_{m,\io}$ defined in \eqref{stanormetta},
{and call $\matY_0(m)$ the normed space $\matY_0$ endowed with the norm $\norm \cdot\norm_{m}$.}
Then, given
{any subset $\matU\subset \matY_0(m)$ and any function $f\!: \matU\to E$,}
where $(E,\|\cdot\|_E)$ is some Banach space, similarly to \eqref{lipnorm} we define
\begin{equation} \label{lipnorm2}
|f|_{\matU,m,E}^{{\rm{Lip}}} := \sup_{(c,\ze) \in \matU} \| f(c,\ze) \|_E \; + 
{\|f\|_{{\rm Lip} (\matU \subseteq \matY_0(m),E)} .}
%\sup_{\substack{(c,\ze),(c',\ze') \in \matY_0 \\ (c,\ze) \neq (c',\ze')}}  \frac{\|f(c,\ze) -f(c',\ze') \|_E}{\norm (c,\ze)- (c',\ze')\norm_{m}}\,  .
\end{equation}

By collecting together the results above, we conclude that, for $\e$ small enough, there are
a constant $C$ and a map
\begin{align} \nonumber
\Psi \!:
\ol{\matU}_{\!\!\rho_0} (\mathtt g(s,\al))
\times\matW_N & \to 
{\matU}_{\rho_\e}(\mathtt g(s,\al))\times [-\rho_\e,\rho_\e]^{N-1}\times {\matU}_{\rho_\e}(\ell^{N,\io}(\RRR))
 \subset \matY_0 \\
 (c,\ze) & \mapsto \Psi(c,\ze):=(\gotU(c,\ze,\e), A(c,\ze,\e),R(c,\ze,\e)) , \nonumber
\end{align}
with $\rho_\e=C|\e|$, such that $\Psi$ is continuous w.r.t.~the product topology, is
Lispschitz continuous, and verifies the bounds
\begin{subequations} \label{doloreesteso}
	\begin{align}
	\label{dolore!exta} 
	|\gotU(\cdot,\cdot,\e)|^{{\rm{Lip}}}_{
	\matY_0,0,\mathtt g(s,\al)} & \le C_1 |\e|,\\
	%|\Uperp(c,\om(\cdot),\e) |^{{\rm{Lip}}}_{\DgN,\mathtt{W}(s_1,s_2,\al)}  & \le C |\e| ,  %\label{dolore!b} \\
	|A(\cdot,\cdot,\e)|^{{\rm{Lip}}}_{
	\matY_0,0,\RRR^{N-1}} & \le C_1 |\e|  , %\qquad q=0,2,\ldots,N-1, 
	\label{muffaext1} \\
	|R(\cdot,\cdot,\e)|^{{\rm{Lip}}}_{
	\matY_0,N,\ell^{N,\io}(\RRR)} & \le C_1 |\e|  \,,
	\label{muffaext2}\\
		|R(\cdot,\cdot,\e)|^{{\rm{Lip}}}_{
		\matY_0,0,\ell^{2,\io}(\RRR)} & \le C_1 |\e|  \, ,
	\label{muffaext3}
	\end{align}
\end{subequations}
with $C_1=C_1(s,\al,N,\g)$ as in Proposition \ref{fabrizio}.

%%%%%%%%%%%%%%%%%%%%%%%%%%%%%%%%%%%%%%%%%%%%%%%%%%%%%%%%%%%%%%%%%%%%%%%%%%
%%%%%%%%%%%%%%%%%%%%%%%%%%%%%%%%%%%%%%%%%%%%%%%%%%%%%%%%%%%%%%%%%%%%%%%%%% 
\zerarcounters
\section{Almost-periodic solutions:~proof of Theorem \ref{main}, part 1}
\label{eppoi}
%%%%%%%%%%%%%%%%%%%%%%%%%%%%%%%%%%%%%%%%%%%%%%%%%%%%%%%%%%%%%%%%%%%%%%%%%% 
%%%%%%%%%%%%%%%%%%%%%%%%%%%%%%%%%%%%%%%%%%%%%%%%%%%%%%%%%%%%%%%%%%%%%%%%%% 

We are now ready to study the existence and regularity of the solutions to \eqref{nls}.  
Recalling that $\ze=(\kappa,\xi)$ and collecting together \eqref{espome}, \eqref{nemmenounnome!}, \eqref{compatibility} and \eqref{WcU},
 we have to solve the set of equations
\begin{subequations} \label{system}
\begin{align}
c + \gotU(c,\kappa,\xi,\e) = W ,
\label{systema} \\
\kappa + A(c,\kappa,\xi,\e)  = 0 , 
\label{systemb} \\ %\qquad q=0, 2,\ldots, N-1 , \\ %\xi_0 + R_0(c,\ze,\e) = V_ 0 , \\
\x + R(c,\kappa,\xi,\e) = V , %\neq 0 , 
\label{systemc}
\end{align} %\qquad j \in \ZZZ 
\end{subequations}
in order to obtain $c=c(V,W,\e)$, $\kappa=\kappa(V,W,\e)$ and $\xi=\xi(V,W,\e)$.    Having done that,
for all $(V,W)$ satisfying $\ze(V,W,\e):=(\ka(V,W,\e),\x(V,W,\e)) \in\matK_N(\g)$,
 we define
\begin{equation}\label{noci}
\UU(x,\f ) %\om(\ze(V,W,\e))t) 
:= \UUU(x,\f; %\om(\ze(V,W,\e))t;
c(V,W,\e),\om(\ze(V,W,\e)),\e) ,
\end{equation}
where $U$ is the solution of the form \eqref{perp} of the modified equation \eqref{NLS4}  whose existence is ensured by Theorem \ref{moser}. 
Then, the function
\begin{equation} \label{lasoluzione!}
u(x,t) =  \UU(x,\om(\ze(V,W,\e))t)
\end{equation}
is an almost-periodic solution to \eqref{nls} of the form \eqref{unnumero}.
Indeed, the first equation in \eqref{system} yields $u(x,0)=\WW(x)$,
while the remaining two ensure that, for all $j\in\ZZZ$, one has
\[
\begin{aligned}
j^2 + V_j & = \omega_j(\ze(V,W,\e)) + \h_j(\ze(V,W,\e)) \\
& =  j^2 + \ka_0 (c,\ze(V,W,\e),\e) + A_0 (c,\ze(V,W,\e),\e) \\
& + \sum_{q=2}^{N-1}
\frac{\ka_q(c,\ze(V,W,\e),\e) + A_q(c,\ze(V,W,\e),\e)}{j^q}+ \xi_j + R_j(c,\ze(V,W,\e),\e),
\end{aligned}
\]
{where the term with the sum is missing for $j=0$,}
so that the compatibility condition \eqref{compatibility} is satisfied and,
if $\ze(V,W,\e)\in\matK_N(\g)$, then \eqref{nls-modified} coincides with \eqref{NLS3} and hence with \eqref{NLS2}.

%%%%%%%%%%%%%%%%%%%%%%%%%%%%%%%%%%%%%%%%%%%%%%%%%%%%%%%%%%%%%%%%%%%%%%%%%% 
\begin{rmk} \label{BCGPvsCGP}
\emph{
With respect to the implicit function problem considered in \cite{CGP}, where the compatibility condition \eqref{compatibility},
for fixed $c\in\ol{\matU}_1(\gsa)$, was ensured by suitably choosing $\ze$ as a function of the parameter $V$,
here we introduce $W$ as a further parameter and aim at choosing both $\ze$ and $c$ as functions of the parameters $V$ and $W$.
}
\end{rmk}
%%%%%%%%%%%%%%%%%%%%%%%%%%%%%%%%%%%%%%%%%%%%%%%%%%%%%%%%%%%%%%%%%%%%%%%%%% 

%%%%%%%%%%%%%%%%%%%%%%%%%%%%%%%%%%%%%%%%%%%%%%%%%%%%%%%%%%%%%%%%%%%%%%%%%% 
\subsection{The implicit function problem, part 1}\vspace{-.2cm}
%%%%%%%%%%%%%%%%%%%%%%%%%%%%%%%%%%%%%%%%%%%%%%%%%%%%%%%%%%%%%%%%%%%%%%%%%% 

Thus, we are left with the problem of solving \eqref{system}.
First, we use \eqref{muffaext1} to solve
the second equation in \eqref{system},
i.e.~the finite-dimensional fixed point  problem
\begin{equation} \label{system2}
\kappa+ A(c,\kappa,\xi,\e)  = 0\,,
\end{equation}
where $c$, $\xi$ and $\e$ are parameters.
From now on
{until the end of the section},
we consider $\e$ fixed once and for all, so we drop it from notation,
and we look for a solution\footnote{{Here and henceforth, for convenience, we shift the order of the two variables $\xi$ and $c$.
The reason for doing that is clarity, since we aim to introduce a correspondence between $(\xi,c)$ and $(V,W)$,
still maintaining the order of the two variables $V$ and $W$ considered so far.}}
{$\kappa=\kappa(\xi,c)$}
to \eqref{system2}.
{Set
\begin{equation} \label{V1}
\matV_0 :=   \ol{\matU}_{\!1/2}(\ell^{N,\io}(\RRR))  \times \ol{\matU}_{\!\rho_0}(\mathtt g(s,\al))
\end{equation}
call $\matV_0(m)$ the Banach space $\matV_0$ endowed with the norm 
\begin{equation}\label{stanormazzabis}
\norm (\xi,c) \norm_{m}:= \max \bigl\{ \|\xi\|_{m,\io} , \|c\|_{s,\al} \bigr\} 
\end{equation}
and, for any subset $\matU \subseteq \matV_0(m)$ and any function $f\!:\matU \to E$, with $E$ being some normed space,
define, by analogy with \eqref{lipnorm2},}
\begin{equation} \label{lipnormservemamancava}
|f|_{\matU,m,E}^{{\rm{Lip}}} := \sup_{(\xi,c) \in \matU} \| f(\xi,c) \|_E \; + 
{\|f\|_{{\rm Lip} (\matU \subseteq \matV_0(m),E)} } .
\end{equation}
Then, by Corollary \ref{lip-c}, for any $\rho_0\in(0,1)$, there exists a solution $\kappa$ to \eqref{system2}, defined on
$\matV_0$, which is continuous w.r.t.~the product topology and satisfies the bound
\begin{equation} \label{kappa-lip}
|\kappa|^{{\rm{Lip}}}_{\matV_0,0,\RRR^{N-1}}  \le 2 C_1 |\e| ,
\end{equation}
with $C_1=C_1(s,\al,N,\g)$ as in \eqref{doloreesteso}.

%%%%%%%%%%%%%%%%%%%%%%%%%%%%%%%%%%%%%%%%%%%%%%%%%%%%%%%%%%%%%%%%%%%%%%%%%% 
\begin{rmk} \label{rmkV1}
\emph{
The space $\matV_0$ defined in \eqref{V1} is closed w.r.t.~the product topology.
}
\end{rmk}
%%%%%%%%%%%%%%%%%%%%%%%%%%%%%%%%%%%%%%%%%%%%%%%%%%%%%%%%%%%%%%%%%%%%%%%%%% 

%%%%%%%%%%%%%%%%%%%%%%%%%%%%%%%%%%%%%%%%%%%%%%%%%%%%%%%%%%%%%%%%%%%%%%%%%% 
\subsection{The implicit function problem, part 2}\vspace{-.2cm}
%%%%%%%%%%%%%%%%%%%%%%%%%%%%%%%%%%%%%%%%%%%%%%%%%%%%%%%%%%%%%%%%%%%%%%%%%% 

Now, we pass to the first and third equations in \eqref{system}.
Introduce, for notation convenience, the Banach spaces $(X_k,\|\cdot\|_{X_k})$, with
\begin{equation}\label{xk}
{X_k := \ell^{k,\io}(\RRR) \times \mathtt g(s,\al) , }
\qquad k \in \NNN ,
\end{equation}
and,
{for any function $F=(F_1,F_2) \in X_k$,}
define the norm
\begin{equation}\label{xknorm}
\| (F_1,F_2) \|_{X_k} := \max \bigl\{ \|F_1\|_{s,\al} , \|F_2\|_{k,\io} \bigr\} .
\end{equation}
{Moreover, given any subset $D\subset X_N$  and any function $F\!:D\to X_N$, set
\begin{equation} \label{brrrr} 
|F|_{D,N}^{{\rm{Lip}}} := \sup_{\mathtt x \in D} \| F( \mathtt x ) \|_{X_N}
+ 
\|F\|_{{\rm Lip}(D\subset X_0,X_2)} 
%\sup_{\substack{\mathtt x , \mathtt x' \in D \\ \mathtt x \neq \mathtt x'}} \frac{\|F (\mathtt x)  -F(\mathtt x')\|_{X_2}}{\|\mathtt x- \mathtt x'\|_{X_0}}
+
 \|F\|_{{\rm Lip}(D \subset X_N,\ell^{N,\io}(\RRR))}  , 
%\sup_{\substack{\mathtt x , \mathtt x' \in D \\ \mathtt x \neq \mathtt x'}} \frac{\|F_2(\mathtt x) -F_2(\mathtt x') \|_{\ell^{N,\io}(\RRR)}}{\|\mathtt x-\mathtt x'\|_{X_N}} ,
\end{equation}
where we have used the inclusion relations $X_{k+1} \subset X_{k}$ for all $k\in\NNN$.}

Next, we introduce the map $
G= %\{G_j\}_{j\in\ZZZ} = 
(G_1,G_2) % = \{(G_{1j},G_{2j})\}_{j\in\ZZZ}
\!:\matV_0 \to X_N$, %\mathtt g(s,\al )\times \ell^{N,\io}(\RRR)$,
by setting
\begin{equation} \nonumber %\label{mapG}
G_{1}(\xi,c):=-\gotU(c,\kappa(\xi,c),\xi,\e)\,,\qquad G_{2}(\xi,c):=-R(c,\kappa(\xi,c),\xi,\e) , %,\quad j\in \ZZZ\,.
\end{equation}
{with $\gotU(c,\kappa,\xi,\e)$ and $R(c,\kappa,\xi,\e)$ being the extensions introduced in Proposition \ref{piergiorgio}.}
By construction the map $G$ is well-defined and continuous w.r.t.~the product topology, and
\[
\sup_{(\xi,c) \in \matV_0} \| G(\xi,c)  \|_{X_N} \le C_1|\e| ,
\]
with $C_1$ as in \eqref{doloreesteso}. Moreover, $G$ satisfies the bound
%:
\begin{equation}\label{ciliegia}
|G|_{\matV_0,N}^{{\rm{Lip}}} \le 
C_1' |\e| .
\end{equation}
for a suitable constant 
{$C_1' =C_1'(s,\al,N,\g,\rho_0) \ge C_1$}
depending also on $\rho_0$ (we refer to the proof of Corollary \ref{lip-c} for a similar argument).

Using the notation above, we rewrite the first and third equations in \eqref{system} as
\begin{equation} \label{systemboh}
(\xi,c) - G(\xi,c)= (V,W) .
\end{equation}
If we fix $\rho_0$ at a value larger than $1/2$, say $\rho_0=2/3$ (see Remark \ref{2/3}), and set
\begin{equation} \label{V0}
\matV := \ol{\matU}_{\!1/4}(\ell^{N,\io}(\RRR)) \times \ol{\matU}_{\!1/2}(\mathtt g(s,\al)) \,\subset\, \matV_0 ,
\end{equation}
we look for a solution  to \eqref{systemboh}
defined on $\matV$
of the form
\begin{equation} \label{prefisso}
(\xi,c)=(V,W) + \Delta(V,W)  
%= (W + \Phi(V,W), V +\xi(V,W)).
= (V + \Delta_1(V,W), W +\Delta_2(V,W)) .
\end{equation}
Note that requiring that \eqref{systemboh} admits a solution of the form \eqref{prefisso}
is the same as requiring that the function $\Delta\!:\matV \to X_N$ solves the fixed point equation
\begin{equation} \label{fisso}
\Delta(V,W)= G \big((V,W)+ \Delta(V,W)\big) .
\end{equation}

{To proceed along the lines outlined above we need a preliminary abstract result,
for which it is useful to rely once more on the notation \eqref{brrrr2}.
%it is convenient to introduce some further notation. If $X$ and $Y$ are normed spaces endowed with the norms $\|\cdot \|_X$ and $\|\cdot\|_Y$,
%respectively, given any non-empty $D\subseteq X$ and any function $f\!:D\to Y$, we define
%
%\begin{equation}\label
%\|f\|_{{\rm Lip}(D\subseteq X,Y)} := \sup_{\substack{\mathtt x,\mathtt x'\in D \\ \mathtt x\neq \mathtt x'}}
%\frac{\|f( \mathtt x)-f( \mathtt x')\|_Y}{\| \mathtt x- \mathtt x' \|_X} .
%\end{equation}
%
The following result is easily proved.\footnote{Here and below we are not assuming that $\|f\|_{{\rm Lip}(D\subseteq X,Y)}$ is finite.} }

%%%%%%%%%%%%%%%%%%%%%%%%%%%%%%%%%%%%%%%%%%%%%%%%%%%%%%%%%%%%%%%%%%%%%%%%%% 
\begin{lemma}\label{lippa}
Let $X$, $Y$ and $Z$ be normed spaces, and let $D\subseteq X$ and $E\subseteq Y$ be given.
For any $f\!:D\to E$ and $g\!:E\to Z$, one has
\begin{equation} \nonumber %\label{cippa}
%\red{\|g\circ f\|^{\rm Lip}_{D,0,Z}}
\|g\circ f\|_{{\rm Lip}(D\subseteq X,Z)}
\leq
%\red{\|f\|^{\rm Lip}_{D,0,Y}}
\|f\|_{{\rm Lip} (D\subseteq X,Y)}
%\red{\|g\|^{\rm Lip}_{E,0,Z}}
\|g\|_{{\rm Lip} (E\subseteq Y,Z)} ,
\end{equation}
where the Lipschitz seminorms are defined according to \eqref{brrrr2}.
\end{lemma}
%%%%%%%%%%%%%%%%%%%%%%%%%%%%%%%%%%%%%%%%%%%%%%%%%%%%%%%%%%%%%%%%%%%%%%%%%% 

%%%%%%%%%%%%%%%%%%%%%%%%%%%%%%%%%%%%%%%%%%%%%%%%%%%%%%%%%%%%%%%%%%%%%%%%%% 
%\begin{rmk} \label{brrrr-rmk}
%\emph{
%By comparing the notation in \eqref{brrrr2} with \eqref{brrrr}, we may write
%
%\begin{equation}\label{cornovaglia}
%|F|_{D,N}^{{\rm{Lip}}} = \sup_{\mathtt x \in D} | F( \mathtt x ) |_{X_N}
%+\|F \|_{{\rm Lip}(D\subseteq X_0,X_2)} +\|F \|_{{\rm Lip}(D\subseteq X_N,\ell^{\io,N}(\RRR))} .
%\end{equation}
%
%}
%\end{rmk}
%%%%%%%%%%%%%%%%%%%%%%%%%%%%%%%%%%%%%%%%%%%%%%%%%%%%%%%%%%%%%%%%%%%%%%%%%% 

Coming back to \eqref{fisso}, we prove the following result.

%%%%%%%%%%%%%%%%%%%%%%%%%%%%%%%%%%%%%%%%%%%%%%%%%%%%%%%%%%%%%%%%%%%%%%%%%% 
\begin{lemma}\label{contraggo}
Fix $N\ge 3$ and $\g\in(0,1)$.
For any $s>0$
{and any $\al\in(0,1)$}, %and all $s_1 \ge 0$ and $s_2>0$ such that $s_1+s_2=s$
%Using the same notations as in Proposition  \ref{fabrizio} fix $N\ge 3$, $\g\in(0,1)$, $s>0$ and $G$ be given by \eqref{mapG}.
there exists $\e_3=
{\e_3(s,\al,N,\g)}
\in (0,\e_2]$, with $\e_2$ as in Proposition  \ref{fabrizio},
such that for any $\e\in  (-{\e_3},{\e_3})$ there exists a map
$\Delta\!:\matV\to X_N$
which solves \eqref{fisso}, is continuous w.r.t.~the product topology and is bounded as
\[
|\Delta|_{\matV,N}^{{\rm{Lip}}} \le C_2 |\e| ,
\]
for a suitable positive constant
{$C_2=C_2(s,\al,N,\g)$.}
\end{lemma}
%%%%%%%%%%%%%%%%%%%%%%%%%%%%%%%%%%%%%%%%%%%%%%%%%%%%%%%%%%%%%%%%%%%%%%%%%% 

%%%%%%%%%%%%%%%%%%%%%%%%%%%%%%%%%%%%%%%%%%%%%%%%%%%%%%%%%%%%%%%%%%%%%%%%%% 
\proof
The proof is based on a fixed point argument for Lipschitz-continuous functions \cite{BMP2}.
Let $(\calmF,\|\cdot\|_{\calmF})$ be the Banach space\footnote{Note that $\matV$ is a closed set also w.r.t. the product topology.
By Tychonoff's theorem $\matV$ is also compact.}
of the functions ${\rm F}\!:\matV \to X_N$ continuous with respect to the product topology, with
\[
\|{\rm F}\|_{\calmF} := \sup_{\mathtt x\in\matV}  \|{\rm F}(\mathtt x) \|_{X_N} ,
\]
where $\mathtt x$ is short for $(V,W)$.
Fix $\rho:=C_1'|\e|$,
with $C_1'$ the constant appearing in \eqref{ciliegia}.
For $\e$ small enough,
if $\| {\rm F} \|_\calmF \le \rho$, then
\begin{equation} \label{cipster}
({\rm Id}+{\rm F})(\matV)\subset \matV_0.
\end{equation}
The bound \eqref{ciliegia} implies that the operator 
$G\circ({\rm Id}+{\rm F})\!: \ol{\matU}_{\!\rho}(\calmF) \longrightarrow \ol{\matU}_{\!\rho}(\calmF)$
%defined as $\Phi({\rm F}):=G\circ({\rm Id}+{\rm F})$,
 is a contraction.
Therefore, there exists a unique function ${\rm F}=\Delta$ which satisfies \eqref{fisso}, and such a function
is continuous w.r.t.~to the product topology~\cite[Lemma B.1]{BMP2}.

Next, we deduce the Lipschitz bounds,
%by a \blue{ bootstrap argument,}
using that
$\Delta$ solves \eqref{fisso},
that is
$\Delta=G\circ({\rm Id}+\Delta)$.
We apply Lemma \ref{lippa}, with %\footnote{Note that, %$\Phi=g\circ f$ and, by \eqref{cipster}, $f\!:\matV_0\to \matV_1$.}
\[
X=Y=X_0, \qquad Z=X_2, \qquad
D = \matV,
\qquad E=\matV_0,
\qquad f={\rm Id}+\Delta , \qquad g=G .
\]
Using that by \eqref{ciliegia} we bound
\[
\|G\|_{{\rm Lip}(\matV_0\subseteq X_0,X_2)}
\le |G|_{\matV_0,N}^{{\rm{Lip}}} \le
C_1'|\e| ,
\]
we get
\begin{equation} \label{scozia}
\begin{aligned} 
\| \Delta \|_{{\rm Lip}(
\matV\subseteq X_0,
X_2)}
& \le C_1'|\e| \| {\rm Id} + \Delta \|_{{\rm Lip}(
\matV\subseteq X_0,
X_0)} \\
& \le C_1'|\e|\bigl(1+ \| \Delta \|_{{\rm Lip}(
\matV\subseteq X_0,
X_0)}\bigr) \\
& \le C_1'|\e|\bigl(1+ \| \Delta \|_{{\rm Lip}(
\matV\subseteq  X_0,
X_2)}\bigr) .
\end{aligned}
\end{equation}
By writing
\[
\Delta (\mathtt x) - \Delta (\mathtt x') = \left( G(\mathtt x + \Delta (\mathtt x)) - G(\mathtt x + \Delta(\mathtt x')) \right) +
\left( G(\mathtt x + \Delta (\mathtt x')) - G(\mathtt x' + \Delta (\mathtt x')) \right) ,
\]
we deduce immediately that $\| \Delta \|_{{\rm Lip}(\matV\subseteq X_0,X_2)}$ is finite
since
\[
\frac{\|\Delta (\mathtt x) - \Delta (\mathtt x') \|_{X_2}}{\|\mathtt x - \mathtt x'\|_{X_0}}
\left( 1- \|G\|_{{\rm Lip}(\matV_0\subseteq X_0,X_2)} \right) \le
\|G\|_{{\rm Lip}(\matV_0\subseteq X_0,X_2)} .
\]
Then \eqref{scozia} implies that
\begin{equation} \label{irlanda}
\|  \Delta \|_{{\rm Lip}(\matV\subseteq X_0,X_2)}
\leq 2 C_1'|\e| ,
\end{equation}
as soon as $|\e| \le 1/2C_1'$.
Thus, \eqref{irlanda} holds for $\e$ small enough.
%\blue{
%In order to prove \eqref{irlanda}, let us assume by contradiction that $\|F \|_{{\rm Lip}(\matV_0\subseteq X_0,X_2)} > 11 C_1|\e|$.
%This means that there exist $\mathtt x, \mathtt x'\in \matV_0$, with $\mathtt x \neq \mathtt x'$,
%such that, choosing as $D$ the set made by such two points, i.e.~$D=\{\mathtt x, \mathtt x'\}$, we find
%\[
%\|F \|_{{\rm Lip}(D\subseteq X_0,X_2)} = \red{
%\frac{|F(x)-F(x')|_{X_2}}{|x-x|_{X_0}} }
%> 11 C_1|\e|,
%\]
%which leads to a contradiction with \eqref{scozia} for $\e$ small enough.
%}
Analogously, we prove that
\begin{equation} \nonumber %\label{irlanda2}
\| \Delta \|_{{\rm Lip}(\matV_0\subseteq X_N,X_N)} \le 2\,  C_1'|\e| , 
\end{equation}
for $\e$ small enough. Finally, the bound
\begin{equation} \nonumber %\label{irlanda3}
\sup_{(V,W) \in \matV_0} \| \Delta( W,V) \|_{X_N} \le C_1' |\e|
\end{equation}
follows trivially from the fact that $\Delta \in \ol{\matU}_{\!\rho}(\matF)$,
with $\rho=C_1'|\e|$.
This complete the proof and gives $C_2=5C_1'$.
\qed
%%%%%%%%%%%%%%%%%%%%%%%%%%%%%%%%%%%%%%%%%%%%%%%%%%%%%%%%%%%%%%%%%%%%%%%%%% 

%%%%%%%%%%%%%%%%%%%%%%%%%%%%%%%%%%%%%%%%%%%%%%%%%%%%%%%%%%%%%%%%%%%%%%%%%%
\begin{rmk} \label{biL}
\emph{
The bound \eqref{ciliegia} on $G$, together with the bound on $\Delta$ in Lemma \ref{contraggo}, yields that the map
$\uno-G\!:\matV_0 \to \ell^{N,\io}(\RRR) \times \mathtt g(s,\al)$ is bi-Lipschitz -- that is Lipschitz and invertible,
with inverse $\uno+\Delta\!:\matV \to \matV_0$ %\ell^{N,\io}(\RRR) \times \mathtt g(s,\al)$
which is also Lipschitz -- and that, for both maps,
the Lipschitz constant is $1+O(\e)$. As pointed out in Remark \ref{BCGPvsCGP}, differently from \cite{CGP}, such a map
allows us to write the initial data (as well as the potentials) in terms of the coefficients $c$ which identify the linear solutions
(as well as of the parameters $\xi$).
}
\end{rmk}
%%%%%%%%%%%%%%%%%%%%%%%%%%%%%%%%%%%%%%%%%%%%%%%%%%%%%%%%%%%%%%%%%%%%%%%%%%

%%%%%%%%%%%%%%%%%%%%%%%%%%%%%%%%%%%%%%%%%%%%%%%%%%%%%%%%%%%%%%%%%%%%%%%%%%
\noindent {\it Proof of Theorem \ref{main}, part 1.}
Collecting together the results proved so far, we have found that, fixed $\g\in (0,1)$ and $\e$ small enough, setting
\begin{equation}\label{paiodi}
\Gamma_N(\g):=\{(V,W)\in \matV \, : \, \ze(V,W)\in \matK_N(\g) \},
\end{equation}
then for any $(V,W)\in \Gamma_N(\gamma)$, the function \eqref{noci} is an almost periodic solution to the Cauchy problem \eqref{nls}.
\qed

%%%%%%%%%%%%%%%%%%%%%%%%%%%%%%%%%%%%%%%%%%%%%%%%%%%%%%%%%%%%%%%%%%%%%%%%%%

%%%%%%%%%%%%%%%%%%%%%%%%%%%%%%%%%%%%%%%%%%%%%%%%%%%%%%%%%%%%%%%%%%%%%%%%%%
\begin{rmk} \label{good}
\emph{
We call the set $\Gamma_N(\g)$ in \eqref{paiodi} the set of \emph{good parameters} in the space of initial data and potentials.
By definition, if $(V,W)$ is a good parameter in $\matV$, then $\ze=\ze(V,W)$ is a good parameter in $\matW_N$.
}
\end{rmk}
%%%%%%%%%%%%%%%%%%%%%%%%%%%%%%%%%%%%%%%%%%%%%%%%%%%%%%%%%%%%%%%%%%%%%%%%%%

%%%%%%%%%%%%%%%%%%%%%%%%%%%%%%%%%%%%%%%%%%%%%%%%%%%%%%%%%%%%%%%%%%%%%%%%%%
\begin{rmk} \label{tildekappabis}
\emph{
Combining Lemma \ref{contraggo} with the definition \eqref{espome} and setting (see \eqref{prefisso})
\[
{\tilde{\kappa}(V,W):=\kappa(\xi(V,W),c(V,W)) , \qquad \xi(V,W) := \Delta_1(V,W) ,}
\]
we obtain from \eqref{espome} and \eqref{prefisso}
\begin{equation} \label{raccapriccio}
\om_j(\ze(V,W)) 
=\begin{cases}
\tilde{\kappa}_0(V,W) + V_0 + \Xi_0(V,W),
& \qquad j=0,\\
\displaystyle{j^2 + \tilde{\kappa}_0(V,W)+\sum_{\pow=2}^{N-1}\frac{\tilde{\kappa}_\pow(V,W)}{j^\pow} + V_j + \Xi_j(V,W),}
& \qquad j\ne0 .
\end{cases}
\end{equation}
The bound \eqref{kappa-lip}, together with \eqref{cipster}, and Lemma \ref{contraggo} impliy that
\[
|\tilde\kappa|^{{\rm{Lip}}}_{\matV,0,\RRR^{N-1}} \le |\kappa|^{{\rm{Lip}}}_{\matV_0,0,\RRR^{N-1}} \le 2 C_1 |\e| ,
\qquad 
\|\xi(V,W)\|_{N,\infty} \le 2 C_1|\e|  .
\]
}
\end{rmk}
%%%%%%%%%%%%%%%%%%%%%%%%%%%%%%%%%%%%%%%%%%%%%%%%%%%%%%%%%%%%%%%%%%%%%%%%%%

%%%%%%%%%%%%%%%%%%%%%%%%%%%%%%%%%%%%%%%%%%%%%%%%%%%%%%%%%%%%%%%%%%%%%%%%%%
\begin{rmk} \label{tildekappa}
\emph{
{
Set, for $\rho\in(0,\rho_0)$,
\begin{equation}\label{vurho}
{\matV_\rho:=  \ol\matU_{\!1/4}(\ell^{N,\infty}(\RRR)) \times \ol\matU_{\!\rho}(\gsa) .}
\end{equation}
By construction one has $\matV_{1/2}=\matV$, as defined in \eqref{V0}, and, for any $\rho\in(0,1/2]$,}
from Proposition \ref{fabrizio} and Remark \ref{ancorascale} we deduce that, for all $(V,W)\in \matV_\rho$,
\begin{subequations} \label{stimepalla}
\begin{align}
|\tilde\kappa_q(V,W)| & \le 2 C_1|\e|\rho^4 , \qquad q=0,2,\dots, N-1 ,
\label{stimepallaa} \\  
\|\xi(V,W)\|_{N,\infty} & \le 2 C_1|\e|\rho^4 .
\label{stimepallab} 
\end{align}
\end{subequations}
Finally, for all $(V,W)\in \matV_\rho$ and all  $\f\in \TTT^\ZZZ$, the function $\calU$ defined in \eqref{noci} is such that
\begin{equation}\label{sposto}
%\begin{aligned}
%& 
\| \calF( 
%{\calU(\cdot,\f;c(V,W),\e))} 
{\calU(\cdot,\f))} 
- W e^{\ii \f}  \|_{s,\al}
%\\ & \qquad \qquad 
\le \|\calF(
{\calU(\cdot,\f))} 
%{\calU(\cdot,\f;c(V,W),\e))} 
-  c(V,W) e^{\ii \f}  \| _{s,\al}  + {\|\Delta_2 \| _{s,\al} }
\le 4 C_1\e \rho^5 ,
%\end{aligned}
\end{equation}
{where we have used the notation \eqref{c2}, for both $W e^{\ii \f}$ and $c(V,W) e^{\ii \f}$,
and the bounds \eqref{Uj} together with the identity \eqref{noci}.}
}
\end{rmk}
%%%%%%%%%%%%%%%%%%%%%%%%%%%%%%%%%%%%%%%%%%%%%%%%%%%%%%%%%%%%%%%%%%%%%%%%%%

%%%%%%%%%%%%%%%%%%%%%%%%%%%%%%%%%%%%%%%%%%%%%%%%%%%%%%%%%%%%%%%%%%%%%%%%%% 
%%%%%%%%%%%%%%%%%%%%%%%%%%%%%%%%%%%%%%%%%%%%%%%%%%%%%%%%%%%%%%%%%%%%%%%%%% 
\zerarcounters
\section{Measure estimates:~proof of Theorem \ref{main}, part 2}
\label{misura}
%%%%%%%%%%%%%%%%%%%%%%%%%%%%%%%%%%%%%%%%%%%%%%%%%%%%%%%%%%%%%%%%%%%%%%%%%% 
%%%%%%%%%%%%%%%%%%%%%%%%%%%%%%%%%%%%%%%%%%%%%%%%%%%%%%%%%%%%%%%%%%%%%%%%%% 

To conclude the proof of Theorem \ref{main} it remains to show that
the set of good parameters $\Gamma_N(\gamma)$ defined in \eqref{paiodi} has large measure for $\g$ small enough.

Let us define
\begin{equation} \label{raccapriccio2}
\ol\om_j(V,W) \!:=\!
\om_j(\ze(V,W)) \!-\! \tilde{\ka}_0(V,W) \!=\!\!
\begin{cases}
V_0 + \Xi_0(V,W),
& \quad j=0,\\
\displaystyle{j^2 + \!\! \sum_{\pow=2}^{N-1}\frac{\tilde{\kappa}_\pow(V,W)}{j^\pow} + V_j +  \Xi_j(V,W), }
& \quad j\ne0 .
\end{cases}
\end{equation}

A consequence of
%mass conservation 
the gauge-invariance \cite[Remark 2.9]{CGP} is the following result.

%%%%%%%%%%%%%%%%%%%%%%%%%%%%%%%%%%%%%%%%%%%%%%%%%%%%%%%%%%%%%%%%%%%%%%%%%%
\begin{lemma}\label{canali}
The set of good parameters \eqref{paiodi} may be defined in terms of the function $\ol\om$ as
\[
{\Gamma_N(\g)	= \bigl\{ (V,W)\in \matV : \ol\om(V,W)\in \gotB(\gamma,N+1) \bigr\} } ,
\]
with $\gotB(\gamma,\tau)$ defined in \eqref{Bbeta1}.
\end{lemma}
%%%%%%%%%%%%%%%%%%%%%%%%%%%%%%%%%%%%%%%%%%%%%%%%%%%%%%%%%%%%%%%%%%%%%%%%%% 

%%%%%%%%%%%%%%%%%%%%%%%%%%%%%%%%%%%%%%%%%%%%%%%%%%%%%%%%%%%%%%%%%%%%%%%%%% 
\prova
Recall the definition of $\matY_0$ and $\matV$ in \eqref{asilo} and in \eqref{V0}, respectively:
then, $(V,W)\in \matV$ ensures that $(c(V,W),\ze(V,W))\in \matY_0$ and,
due to the bound \eqref{kappa-lip} and to Lemma \ref{contraggo}, 
both $\om(\ze(V,W))$ and $\ol\om(V,W)$ belong to $\gQ$ and
$\om(\ze(V,W))\in\gotB_\be$
%(\g,N+1)$
if
$\om=\om(\ze(V,W))$ is such that
$
\be^{(0)}_\om(2^m)\ge \be_m
$
%^*(2^m,\g,\tau)
%\qquad \mbox{
for all 
%}
$m\ge0$.
For any $\nu\in \ZZZ_{f,0}^\ZZZ$, from the definition \eqref{interib} we deduce that
\[
\om(\ze(V,W))\cdot\nu- \ol\om(V,W)\cdot\nu= \tilde{\ka}_0(V,W) \sum_{i\in\ZZZ} \nu_i = 0 ,
\]
so that the thesis immediately follows.
\qed
%%%%%%%%%%%%%%%%%%%%%%%%%%%%%%%%%%%%%%%%%%%%%%%%%%%%%%%%%%%%%%%%%%%%%%%%%% 

%\vspace{.3cm}

%%%%%%%%%%%%%%%%%%%%%%%%%%%%%%%%%%%%%%%%%%%%%%%%%%%%%%%%%%%%%%%%%%%%%%%%%% 
\subsection{Measure of the set of good parameters}\vspace{-.2cm}
%%%%%%%%%%%%%%%%%%%%%%%%%%%%%%%%%%%%%%%%%%%%%%%%%%%%%%%%%%%%%%%%%%%%%%%%%% 

For all $\nu\in \ZZZ^{\ZZZ}_{f,0} \setminus\{0\}$ and all $\de>0$, set 
\begin{equation} \label{Rnude}
\mathfrak R_\nu(\de):=\bigl\{(V,W)\in \matV : |\ol\om(V,W)\cdot \nu|\le \de \bigr\} .
\end{equation}
The set $\mathfrak R_\nu(\de)$ is measurable w.r.t.~the measure induced by the product measure on $\matV$
since it is the preimage of  a closed set
under the map $(V,W) \mapsto \ol\om(V,W)\cdot\nu$,
which, by Lemma \ref{contraggo}, is continuous in the product topology. 
Similarly, for each $W\in \ol{\matU}_{\!{1/2}}(\mathtt g(s,\al))$, the section
\[
\mathfrak R_\nu(W,\de) := \Bigl\{V\in\ol{\matU}_{\!{1/4}}(\ell^{N,\io}(\RRR)) : |\ol\om(V,W)\cdot \nu|\le \de \Bigr\}
\]
is measurable w.r.t.~the measure induced by the product measure on $\ol{\matU}_{\!{1/4}}(\ell^{N,\io}(\RRR))$.
Our aim is to find an upper bound on the measure  of $\mathfrak R_\nu(W,\de)$ and consequently of $\mathfrak R_\nu(\de)$.
This is provided by the following result.

%%%%%%%%%%%%%%%%%%%%%%%%%%%%%%%%%%%%%%%%%%%%%%%%%%%%%%%%%%%%%%%%%%%%%%%%%% 
\begin{lemma}\label{misuretta}
For all $W\in \ol{\matU}_{\!{1/2}}(\mathtt g(s,\al))$ and all $\nu\in \ZZZ^{\ZZZ}_{f,0} \setminus\{0\}$
there exists $\de_0>0$ such that for all $\de\in(0,\de_0)$
\begin{equation} \label{upper}
\ol\mu_{1/4,1/2}(\mathfrak R_\nu(\de)) = 
\mu_{2,1/2} \left( \mu_{1,1/4} \left( \mathfrak R_\nu(\cdot,\de) \right) \right) 
\le \mu_{1,1/4} \left( \mathfrak R_\nu(W,\de) \right)
\le C \jap{i_0(\nu)}^N \de ,
\end{equation}
{where $C$ is an absolute constant}, and %positive constant independent of $N$, and 
\begin{equation} \label{i0max}
i_0(\nu):= \max\bigl\{i\in \ZZZ:   |\nu_i|=\|\nu\|_{\io}\bigr\} .
\end{equation}
\end{lemma}
%%%%%%%%%%%%%%%%%%%%%%%%%%%%%%%%%%%%%%%%%%%%%%%%%%%%%%%%%%%%%%%%%%%%%%%%%% 

%%%%%%%%%%%%%%%%%%%%%%%%%%%%%%%%%%%%%%%%%%%%%%%%%%%%%%%%%%%%%%%%%%%%%%%%%% 
\proof
For all $\nu\in\ZZZ^{\ZZZ}_{f,0}$ we have
\[
\ol\omega(V,W)\cdot \nu = \sum_{j\in \ZZZ} j^2 \nu_j + V\cdot\nu
+ \xi(V,W)\cdot\nu +  \sum_{q=2}^{N-1} \tilde\kappa_q(V,W) \sum_{j\in \ZZZ\setminus\{0\}} 
\frac{\nu_j}{j^q}.
\]
Fix $\nu\neq0$ and write $V= (V_{i_0},V')$, with $i_0=i_0(\nu)$ defined in \eqref{i0max} and $V'=\{V_j\}_{j\in\ZZZ\setminus\{i_0\}}$.

For $h\neq0$ small enough we obtain, from \eqref{raccapriccio2},
for some positive
{constant $C_3=C_3(s,\al,N)$,}
\begin{equation} \nonumber %\label{C(N)}
\begin{aligned} 
& \frac{|\ol \om((V_{i_0}+h,V'),W)\cdot \nu - \ol\om((V_{i_0},V'),W) \cdot \nu |}{h} \\
& \qquad \ge | \nu_{i_0}| 
\Biggl( 1 - |\Delta|^{{\rm{Lip}}}_{\matV_0,0,\ell^{2,\io}(\RRR)}\sum_{j\in\ZZZ}\jap{j}^{-2} 
- \sum_{q=2}^{N-1} |\tilde\kappa_q|^{\text{Lip}}_{\matV_0,0,\RRR}  \sum_{j\in \ZZZ} \jap{j}^{-q} \Biggr) \\
& \qquad \ge  |\nu_{i_0} |(1- C_3 |\e| ) \ge \frac12,
\end{aligned}
\end{equation}
{where we have used the bound on the norm $|\Delta|_{\matV,N}^{{\rm{Lip}}}$ in Lemma \ref{contraggo}
to estimate the norm $|\Delta|^{{\rm{Lip}}}_{\matV_0,0,\ell^{2,\io}(\RRR)}$.}

Reasoning as in \cite{CGP}
(which in turn follows closely \cite{BMP3}),
the latter bound allows us to conclude that, for $\de$ small enough, the set
$\mathfrak R_\nu(W,\de)$ is contained in the normal domain
\[
\EEE :=\Biggl\{ V\in \ell^{N,\io}(\RRR):  a(V') \le V_{i_0} \le b(V') \,, \; 
\sup_{\substack{j\in\ZZZ \\  j\ne i_0}}|V_j| \jap{j}^N \le \frac{1}{4} \Biggr\} ,
\]
where $a$ and $b$ are two  functions, continuous w.r.t.~the product topology, such that 
\[
b(V') |i_0|^{N} \le \frac{1}{4} ,
\qquad  
0 \le b(V')- a (V') \le 2\de , \qquad \forall V'\in \ol{\matU}_{\!1/4}(\ell^{N,\io}(\RRR)).
\]
Thus, we may bound
\[
\mu_{1,1/4}({\EEE}) \le 8 \jap{i_0}^N \de 
\]
and hence, by Fubini's theorem, we obtain
\[
\ol\mu_{1/4,1/2}(\mathfrak R_\nu(\de)) = 
\mu_{2,1/2} \left( \mu_{1,1/4} \left( \mathfrak R_\nu(\cdot,\de) \right) \right) \le
\mu_{1,1/4} \left( \mathfrak R_\nu(\cdot,\de) \right)  \le
\mu_{1,1/4}({\EEE}) \le 8 \jap{i_0}^N \de. 
\]
This yields the bound \ref{upper}, with $C=8$.
\qed
%%%%%%%%%%%%%%%%%%%%%%%%%%%%%%%%%%%%%%%%%%%%%%%%%%%%%%%%%%%%%%%%%%%%%%%%%%

%%%%%%%%%%%%%%%%%%%%%%%%%%%%%%%%%%%%%%%%%%%%%%%%%%%%%%%%%%%%%%%%%%%%%%%%%%
\begin{rmk} \label{rmkC(N)}
\emph{
We stress that
{the bounds of \cite{CGP} provide in the proof of Lemma \ref{misuretta}
a constant $C_3(s,\al,N)$  such that $C_3(s,\al,N) \to \io$ as $N\to \io$.}
This prevents us from obtaining uniformity in $N$ and trivially extending our result to convolution potentials of class $C^{\io}$.
}
\end{rmk}
%%%%%%%%%%%%%%%%%%%%%%%%%%%%%%%%%%%%%%%%%%%%%%%%%%%%%%%%%%%%%%%%%%%%%%%%%%

%%%%%%%%%%%%%%%%%%%%%%%%%%%%%%%%%%%%%%%%%%%%%%%%%%%%%%%%%%%%%%%%%%%%%%%%%%
\begin{prop} \label{44}
Let $\Gamma_N(\g)$ be defined as in Lemma \ref{canali}.
There is an absolute constant $C_*$ such that, for any $\g>0$,
{the set $\matV \setminus\Gamma_N(\g)$ has measure ${\ol\mu_{1/4,1/2}} (\matV \setminus\Gamma_N(\g))\le C_*\g$.}
\end{prop}
%%%%%%%%%%%%%%%%%%%%%%%%%%%%%%%%%%%%%%%%%%%%%%%%%%%%%%%%%%%%%%%%%%%%%%%%%%

%%%%%%%%%%%%%%%%%%%%%%%%%%%%%%%%%%%%%%%%%%%%%%%%%%%%%%%%%%%%%%%%%%%%%%%%%%
\prova
Fix $\de=\g \de_{\nu}^{N+1}$ in \eqref{Rnude}, with
\begin{equation} \label{deltanu}
\de_{\nu}:= \prod_{i\in \ZZZ} (1+\jap{i}^2\nu_i^2)^{-1} .
\end{equation}
Then we have,
{by Lemma \ref{misuretta},} 
%\cite[Lemma B.5]{CGP}
%
\begin{equation} \label{final-measure}
%\meas
{\ol\mu_{1/4,1/2}}
(\mathfrak R_\nu(\g \de_\nu^{N+1}))\le \g \de_\nu 
\end{equation}
and hence
\begin{equation} \nonumber
%\meas
{\ol\mu_{1/4,1/2}} (\matV \setminus\Gamma_N(\g))\le
\sum_{\nu\in\ZZZ^\ZZZ_{f,0}}
%\meas
{\ol\mu_{1/4,1/2}}
(\mathfrak R_\nu(\g \de_\nu^{N+1}))\le \g 
\sum_{ \nu\in\ZZZ^\ZZZ_{f}} 
\prod_{i\in \ZZZ} 
\frac{1}{1+\jap{i}^2\nu_i^2}<C_*\g ,
\end{equation}
where the proof of the last inequality is due to Bourgain \cite{Bjfa}
and $C_*$ is a suitable constant independent of $\g$.
\qed
%%%%%%%%%%%%%%%%%%%%%%%%%%%%%%%%%%%%%%%%%%%%%%%%%%%%%%%%%%%%%%%%%%%%%%%%%%

%%%%%%%%%%%%%%%%%%%%%%%%%%%%%%%%%%%%%%%%%%%%%%%%%%%%%%%%%%%%%%%%%%%%%%%%%%
\begin{rmk} \label{tau=N+1}
\emph{
We need $\de_\nu$ as in \eqref{deltanu} in order to obtain an upper bound on the measure of the set
$\matV_0\setminus\Gamma_N(\g)$ proportional to $\g$. On the other hand,
requiring that $(V,W) \notin {\mathfrak R}_{\nu}(\g\de_{\nu}^{N+1})$ for all $\nu\in\ZZZ^\ZZZ_{f,0} \setminus\{0\}$,
with ${\mathfrak R}_{\nu}(\de)$ defined in \eqref{Rnude}, implies that the frequency $\om(\ze(V,W))$
belongs to the set $\gD^{(0)}(\g,\tau)$ introduced in Remark \ref{nuovormk} with $\tau=N+1$.
This explains why we have fixed $\tau=N+1$ in Section \ref{uni} (see Remark \ref{fuoridalnulla}).
}
\end{rmk}
%%%%%%%%%%%%%%%%%%%%%%%%%%%%%%%%%%%%%%%%%%%%%%%%%%%%%%%%%%%%%%%%%%%%%%%%%%

%%%%%%%%%%%%%%%%%%%%%%%%%%%%%%%%%%%%%%%%%%%%%%%%%%%%%%%%%%%%%%%%%%%%%%%%%%
\begin{rmk} \label{misuraomega0}
\emph{
{If, instead of the frequency $\ol\om(W,V)$, we consider $\omega^{(0)}(V) :=\{j^2+V_j\}_{j\in \ZZZ}$ and define
\[
\Gamma^{(0)}_N(\g) = \bigl\{ (V \in \ol\matU_{\!1/4}(\ell^{N,\infty}(\RRR))  : \om^{(0)}(V)\in \gotB(\g,N+1) \bigr\} ,
\]
then, reasoning along the same lines as in the proof of Proposition \ref{44} and possibly redefining the constant $C_*$,
we find that $\mu_{1,1/4}(\ol\matU_{\!1/4}(\ell^{N,\infty}(\RRR)) \setminus \Gamma^{(0)}_N(\g)) \le C_*\g$ as well.}
}
\end{rmk}
%%%%%%%%%%%%%%%%%%%%%%%%%%%%%%%%%%%%%%%%%%%%%%%%%%%%%%%%%%%%%%%%%%%%%%%%%%

%%%%%%%%%%%%%%%%%%%%%%%%%%%%%%%%%%%%%%%%%%%%%%%%%%%%%%%%%%%%%%%%%%%%%%%%%%
\subsection{Measure estimates in the space of initial data}\vspace{-.2cm}
%%%%%%%%%%%%%%%%%%%%%%%%%%%%%%%%%%%%%%%%%%%%%%%%%%%%%%%%%%%%%%%%%%%%%%%%%%

Now, we have all the ingredients to conclude our analysis and prove the measure estimates stated in Section \ref{setup}.

\vspace{.3cm}

%%%%%%%%%%%%%%%%%%%%%%%%%%%%%%%%%%%%%%%%%%%%%%%%%%%%%%%%%%%%%%%%%%%%%%%%%%
\noindent {\it Proof of Theorem \ref{main}, part 2.}
Since the set $\Gamma_N(\g)$ is measurable, together with its sections
respect to the $V$ and $W$ variables, a direct application of Fubini's Theorem ensures that,
for $\g$ and $\e$ small enough,  there exists a positive measure set
$\mathcal{G}\in \ol{\matU}_{\!1/4}(\ell^{N,\io}(\RRR))$ such that  for all $V\in \mathcal{G}$
the section
\begin{equation} \label{TV}
\mathcal T_V := 
\Gamma_N(\g)\vert_V
= \Bigl\{ W \in \ol{\matU}_{\!1/2}(\mathtt g(s,\al)) : (V,W) \in \Gamma_N(\g) \Bigr\}
\end{equation}
has measure proportional to $1- O(\g)$ and hence close to 1 for $\g$ small.
By construction, the condition that
$V\in \mathcal{G}$ and, correspondingly,
$W\in \mathcal T_V$ ensures that $(V,W)\in \Gamma_N(\g)$ and hence
{$\ze(V,W)\in \matK_N(\g)$.}
Set
\begin{equation} \label{g0eps}
{\g_0(\e) = \g_0(\e,s,\al,N):= %2 
\inf\{\g>0\;:\;  |\e|\le \e_3(s,\al,N,\g)\} ,}
\end{equation}
with $\e_3(s,\al,N,s,\g)$ as in Lemma \ref{contraggo}.
Using the results of Section 2.4 in \cite{CGP}, one sees that $\gamma_0(\e)\to 0$ as $\e\to 0$,
which implies that $\mathcal G$ and $\TT_V$ have asymptotically full measure.
This completes the proof of Theorem \ref{main}.
\qed
%%%%%%%%%%%%%%%%%%%%%%%%%%%%%%%%%%%%%%%%%%%%%%%%%%%%%%%%%%%%%%%%%%%%%%%%%%

%%%%%%%%%%%%%%%%%%%%%%%%%%%%%%%%%%%%%%%%%%%%%%%%%%%%%%%%%%%%%%%%%%%%%%%%%% 
%%%%%%%%%%%%%%%%%%%%%%%%%%%%%%%%%%%%%%%%%%%%%%%%%%%%%%%%%%%%%%%%%%%%%%%%%% 
\zerarcounters
\section{Reuslts for a full measure set of potentials}
\label{stab}
%%%%%%%%%%%%%%%%%%%%%%%%%%%%%%%%%%%%%%%%%%%%%%%%%%%%%%%%%%%%%%%%%%%%%%%%%% 
%%%%%%%%%%%%%%%%%%%%%%%%%%%%%%%%%%%%%%%%%%%%%%%%%%%%%%%%%%%%%%%%%%%%%%%%%% 

We pass now to the proof of Theorem \ref{eureka} and of the stability result stated in Theorem \ref{proposta}.
We start with a scaling argument. For any
{$\rho\in(0,1/2]$,}
let $\matV_\rho$ be defined as in \eqref{vurho}.
Consider \eqref{nls}, with both $\e$ and $(V,W)\in \matV_\rho$ fixed.
As seen in Remark \ref{riscalo}, equation \eqref{nls} has an almost-periodic solution
{$u$  if  and only if the rescaled equation \eqref{nls}
\begin{equation} \nonumber %\label{nls1}
\begin{cases}
\ii u_t - u_{xx} + \VV * u +(2 \rho)^4 \e |u|^4u =0 , \qquad x\in\TTT, \\
u(x,0) = \WW'(x) ,
\end{cases}
\end{equation}
has an  almost-periodic solution $u' \!=\! (2\rho)^{-1} u$ with initial  datum $\WW' \!=\! (2\rho)^{-1}\WW$.
{For any $\g_1>0$ let $\rho_1 \in(0,1/2)$ be such that}
\begin{equation}\label{boh??}
\gamma_1 >  \gamma_0((2\rho_1)^4\e),
\end{equation}
with $\g_0(\e)$ as in \eqref{g0eps}.
The scaling properties mentioned in Remark \ref{riscalo} ensure that %, at fixed $\e$,
we can find an almost-periodic solution for all $(V,W)\in \matD(\g_1)$, with
\begin{equation} \label{digamma}
\matD(\g_1) := \{(V,W)\in \matV_{\rho_1} : \ol\om(V,(2\rho_1)^{-1}W)\in \gD^{(0)}(\g_1,N+1) \}  ,
\end{equation}
where $\matV_\rho$ is given by \eqref{vurho} and $\gD^{(0)}(\g,\tau)$ is
the set of $(\g,\tau)$-weak Diophantine vectors defined in Remark \ref{nuovormk}.
{Note that, in \eqref{digamma}, for $\mu_0(\gD^{(0)}(\g_1,N+1))$ to be positive we need $\g_1$
to be such that $K_1 \g_1<1$, with $K_1:=C(N+1)$, if $C(\tau)$ is as defined in Remark \ref{nuovormk}.}

%%%%%%%%%%%%%%%%%%%%%%%%%%%%%%%%%%%%%%%%%%%%%%%%%%%%%%%%%%%%%%%%%%%%%%%%%%%
\begin{rmk} \label{rho0rho1}
\emph{
The condition $\rho_1 < 1/2$ ensures that $\rho_1 < \rho_0$, provided that $\rho_0$, as introduced in Corollary \ref{lip-c},
is chosen greater than 1/2 according to Remark \ref{2/3}.
}
\end{rmk}
%%%%%%%%%%%%%%%%%%%%%%%%%%%%%%%%%%%%%%%%%%%%%%%%%%%%%%%%%%%%%%%%%%%%%%%%%%%

%%%%%%%%%%%%%%%%%%%%%%%%%%%%%%%%%%%%%%%%%%%%%%%%%%%%%%%%%%%%%%%%%%%%%%%%%%%
\subsection{Abundance of almost-periodic solutions}\vspace{-.2cm}
%%%%%%%%%%%%%%%%%%%%%%%%%%%%%%%%%%%%%%%%%%%%%%%%%%%%%%%%%%%%%%%%%%%%%%%%%%%

We first prove that for a full measure set of potentials in $\ol{\matU}_{1/4}(\ell^{N,\io}(\RRR))$
the measure of the set of initial data which give rise to almost-periodic solution has asymptotically full measure.

\medskip

\noindent
{\it Proof of Theorem \ref{eureka}.}
From here on, we assume $\e$ to be fixed once and for all.
Fix $\g_1 \in(0,\g_*)$, with $\g_*$ to be determined,
and define $\rho_1$ according to \eqref{boh??}.
For any $\rho\in(0,\rho_1)$, define $\matD_{\rho}(\gamma_1) :=\matV_\rho\cap \matD(\gamma_1)$,
{with $\matD(\gamma_1)$ as in \eqref{digamma}.}
In the rest of the proof we call $C$ any constant depending on $\g_1$, but not on $\e$ and $\rho$.

Set $\omega^{(0)}(V) :=\{j^2+V_j\}_{j\in \ZZZ}$.
By \eqref{raccapriccio2} and \eqref{stimepalla}, there is a constant $C$ such that
\begin{equation} \label{Cg1}
\| \bar\om( V,\rho_0^{-1}W)  - \omega^{(0)}(V)\|_\infty \le C \rho^4 |\e|
\end{equation}
for all $(V,W)\in \matV_\rho$. Define 
\begin{equation}\label{diaframma}
\mathfrak D^{(0)}(\g_1) := \bigl\{ V\in \ol{\matU}_{\!1/4}(\ell^{N,\io}(\RRR)) : 
\om^{(0)}(V) \in \gD^{(0)}(2\g_1,N+1) 
% |\om^{(0)}(V)\cdot \nu | \ge \gamma_1 \delta_\nu^{N+1}  \quad \forall \nu\in\ZZZ^\ZZZ_{f,0}
\bigr\}  
\end{equation}
%
%\red{with $\de_\nu$ defined in \eqref{deltanu},}
and, for $V\in \mathfrak D^{(0)}(\gamma_1)$,
\[
\matS_V (\g_1) := \bigl\{ W\in \ol\matU_{\rho_1}(\gsa) : (V,W)\in \matD({\gamma_1)} \bigr\} . 
\]
Set also $\mathfrak D^{(0)}_{\rho}(\g_1) := \ol\matU_\rho(\gsa)\times  \mathfrak D^{(0)}(\g_1)$
and $\matD^{(0)}_{\rho}(\g_1) :=\matD_{\rho}(\g_1) \cap \mathfrak D^{(0)}_{\rho}(\g_1)$.
{Note that, to ensure $\mu_0(\gD^{(0)}(2\g_1,N+1))$ to be positive,
we need $2 K_1\g_1<1$, with $K_1$ as defined after \eqref{digamma}.}
{Using the measure estimate in Remark \ref{misuraomega0} and the fact that $\gD^{(0)}(\g_1,\tau) \subset \gotB(\g_1,\tau)$
and hence $\gD^{(0)}(\g_1,N+1) \subset \Gamma^{(0)}_N(\g)$
we may estimate, for a suitable constant $K_2$,
\begin{equation}\label{seriali}
\mu_{1,1/4}(\mathfrak D^{(0)}(\g_1)) \ge 1 - K_2 \g_1,
\qquad  \ol\mu_{1/4,\rho} ( \mathfrak D^{(0)}_{\rho}(\g_1)) \ge 1 - K_2 \g_1 ,
\end{equation}
provided $\g_1$ is such that $K_2\g_1<1$.}

Moreover, for all $(V,W)\in \mathfrak D^{(0)}_{\rho}(\g_1)$,
we have, trivially,\footnote{Recall that $\|\cdot\|_1$ denotes the $\ell^1$-norm (see footnote \ref{11}).}
\[
|\bar\om( V,\rho_0^{-1}W) \cdot \nu|\ge
|\omega^{(0)}(V)\cdot \nu |-  C \rho^4 |\e| \, \|\nu\|_1 \ge 
2 \g_1 \de_\nu^{N+1} - C \rho^4 |\e| \, \|\nu\|_1 \ge \g_1 \de_\nu^{N+1} %\ge \g_1  \delta_\nu^{N+1}
\]
for all $\nu\in\ZZZ^\ZZZ_{f,0}$ such that $\gamma_1 \delta^{N+1}_\nu \ge C |\e| \rho^4 \|\nu\|_1$, with $C$ as in \eqref{Cg1}. Set
\[
{\mathfrak R}_{\rho,\nu}(\de):=  \mathfrak D^{(0)}_{\rho}(\g_1) \cap \mathfrak R_\nu(\delta)
\]
with $ \mathfrak R_\nu(\delta)$ defined as in \eqref{Rnude}.
Analogously to \eqref{final-measure}, we find
\[
\ol\mu_{1/4,\rho} ({\mathfrak R}_{\rho,\nu}(\g_1\delta_\nu^{N+1}) ) \le C \g_1 \de_{\nu}\,,
\]
for some  constant $C$.
Thus, using that the sum
\[
\sum_{\nu\in\ZZZ^\ZZZ_{f}} \|\nu\|_1^a \de_{\nu}^{1-a} , \qquad a \in [0,1] ,
\]
is finite for $a<1/3$, as it is straightforward to check, and choosing, say, $a=1/4$, we obtain
\begin{equation}\label{questaqui}
\begin{aligned}
\ol\mu_{1/4,\rho} (\mathfrak D^{(0)}_{\rho}(\g_1) \setminus \matD^{(0)}_{\rho}(\g_1))
& \le \sum_{\nu\in\ZZZ^\ZZZ_{f,0}} \ol\mu_{1/4,\rho} ({\mathfrak R}_{\rho,\nu}(\g_1\delta_\nu^{N+1}) ) \\
& \le  C {\g_1}
\!\!\!\!\!\!\!\!\!\!\!
\sum_{\substack{\nu\in\ZZZ^\ZZZ_{f,0} \\ \g_1 \delta^{N+1}_\nu < C |\e| \rho^4\|\nu\|_1}}
\!\!\!\!\!\!\!\!\!\!\!
\delta_\nu 
%\le  C \g_1^{5/6} \rho^{2/3} \sum_{\nu\in\ZZZ^\ZZZ_{f} }\delta_\nu^{2/3} \le C \rho^{2/3} ,
\le  C \g_1 (\gamma^{-1/4}|\e|^{1/4} \rho)^{\frac{1}{N+1}}\sum_{\nu\in\ZZZ^\ZZZ_{f} } \|\nu\|_1^{1/4} \delta_\nu^{3/4} \\
& \le C_*(\e,\g_1) \, \rho^{1/N+1} ,
\end{aligned}
\end{equation}
with $C_*(\e,\g_1)=C |\e|^{1/(4N+4)}\g_1^{(4N+3)/(4N+4)} $, for some other constant $C$,
{provided $\g_1$ is such that $C_*(\e,\g_1)  \rho^{1/(N+1)} < 1$.}

Next, we consider a positive non-increasing sequence $\varrho=\{\rho_k\}_{k\ge1}\in \ell^{1/2(N+1)}(\NNN,\RRR_+)$  namely such that 
\[
\sum_{k \ge 1} \rho_k^{1/2(N+1)} < \io
\]
and fix $N_1=N_1(\g_1) \in \NNN$ so that
\begin{equation}\label{enneuno}
\sum_{k>N_1} \rho_k^{1/2(N+1)} < \g_1 .
\end{equation}
%
%Consider the sequence of measurable sets  $\mathfrak D^{(0)}_{\rho_k ,\gamma_1}$ and $ \matD_{\rho_k, \gamma_1} $, with $k>N_1$.
Define
\begin{equation} \label{Brhogamma}
\mathfrak B(\rho,\gamma_1) := \left\{ V\in \mathfrak D^{(0)}(\gamma_1) :
\mu_{2,\rho}(\matS_V^c (\g_1)\cap \ol\matU_{\rho}(\gsa))\ge 2 C_*(\e,\g_1) \,  \rho^{1/2(N+1)} \right\} ,
\end{equation}
{and restrict further $\g_1$ by requiring $2C_*(\e,\g_1) \,  \rho^{1/(N+1)} < 1$.}
%and set $\mathfrak B_k(\varrho,\g_1)=\mathfrak B(\rho_k,\gamma_1)$.
By Fubini's theorem, for each $k>N_1$ the set 
$\mathfrak B(\rho_k,\g_1)$
%$\mathfrak B_k(\varrho,\g_1)$
is measurable and 
its measure $\mu_{1,1/4}(\mathfrak B(\rho_k,\g_1)$ can not exceed $\rho_k^{1/2(N+1)}$,
otherwise the estimate \eqref{questaqui} would lead to a contradiction. By \eqref{enneuno}, this yields that
\begin{equation}\label{misuragrande}
\mu_{1,1/4}(\mathfrak D^{(0)}(\g_1)\setminus(\cup_{k>N_1} \mathfrak B(\rho_k,\g_1)) \ge 1 - \g_1 
\end{equation}
and, by construction, for all $V\in \mathfrak D^{(0)}(\g_1) \setminus(\cup_{k>N_1} \mathfrak B(\rho_k,\g_1))$
and for all $k>N_1$, one has
\[
\liminf_{k\to\io} \mu_{2,\rho_k}(\matS_V^c\cap \ol\matU_{\rho_k}(\gsa)) = \lim_{k\to\io} 2 C_*(\e,\g_1) \rho_k ^{1/2(N+1)} = 0 .
\]
{By collecting together all the conditions imposed on $\g_1$, we end up  requiring $\g_1<\g_*$ with
%	, with
%
\begin{equation} \label{g*}
\frac{1}{\g_*} =  2K_1 + \max\{K_1,K_2, (2 C |\e|)^{1/(4N+3)}, 1 \} ,
\end{equation}
where the first summand is due to the fact that we are restricting the analysis to potentials $V \in \ol{\matU}_{\!1/4}(\ell^{N,\io}(\RRR)) \cap
\mathfrak D^{(0)}(\g_1)$. This fixes the value of $\g_*$ introduced at the beginning.
}

Finally, consider an arbitrary non-increasing sequence $\{\gamma_k\}_{k\ge1}$,
with $\g_1$ as above and such that $\gamma_k \to 0^+$ as $k\to\io$,
and define
%$N_k$ by requiring \eqref{enneuno} to be satisfied with $\g_1$ replaced with $\g_k$.
%If we set
\[
\mathfrak G_{\varrho} :=\bigcup_{k\ge1} \biggl(\mathfrak D^{(0)}(\g_k)\setminus \bigcup_{h> N_1(\g_k)}\mathfrak B(\rho_h,\g_k)\biggr) .
\]
Assume that $\mu_{1,1/4}(\mathfrak G_{\varrho}) = 1 - \de$, with $\de>0$. Then, as soon as $\g_k < \de$,  we find,
by using the bound \eqref{misuragrande},
\[
\begin{aligned}
1 - \g_k & \le \mu_{1,1/4}\biggl(\mathfrak D^{(0)}(\g_k)\setminus \bigcup_{h> N_1(\g_k)}\mathfrak B(\rho_h,\g_k )\biggr) \\ &\le
\mu_{1,1/4} \biggl( \; \bigcup_{k \ge 1} \biggl(\mathfrak D^{(0)}(\g_k)\setminus \bigcup_{h> N_1(\g_k)}\mathfrak B(\rho_h,\g_k )\biggr) \biggr)
= 1 - \de ,
\end{aligned}
\]
which leads to a contradiction.
Then the assertion follows.
\qed
%%%%%%%%%%%%%%%%%%%%%%%%%%%%%%%%%%%%%%%%%%%%%%%%%%%%%%%%%%%%%%%%%%%%%%%%%% 

%%%%%%%%%%%%%%%%%%%%%%%%%%%%%%%%%%%%%%%%%%%%%%%%%%%%%%%%%%%%%%%%%%%%%%%%%% 
\subsection{Lyapunov statistical stability of the origin}\vspace{-.2cm}
%%%%%%%%%%%%%%%%%%%%%%%%%%%%%%%%%%%%%%%%%%%%%%%%%%%%%%%%%%%%%%%%%%%%%%%%%% 

Next, we use the results above in order to prove that the elliptic equilibrium point at the origin is Lyapunov statistically stable.

\medskip

%%%%%%%%%%%%%%%%%%%%%%%%%%%%%%%%%%%%%%%%%%%%%%%%%%%%%%%%%%%%%%%%%%%%%%%%%% 
\noindent\emph{Proof of Theorem \ref{proposta}}.
Fix $\g\in(0,\g_*)$, with $\g_*$ as in \eqref{g*}, and $\rho$ such that $\g > \g_0((2\rho)^4\e)$.
{Define $\matD(\g)$ as in \eqref{digamma} and $\mathfrak D^{(0)}(\g)$ as in \eqref{diaframma},
with $\gamma$ and $\rho$ replacing $\gamma_1$ and $\rho_1$, respectively,
and set $\matD_{\rho}(\g) :=\matV_\rho\cap \matD(\g)$,
$\mathfrak D^{(0)}_{\rho}(\g) :=  \ol\matU_\rho(\gsa) \times  \mathfrak D^{(0)}(\g)$
and $\matD^{(0)}_{\rho}(\g) :=\matD_{\rho}(\g) \cap \mathfrak D^{(0)}_{\rho}(\g)$.}
{Then, by construction, $\matD(\g) \subset\mathfrak G$, so that}
we have an almost-periodic solution for any $(V,W)\in \matD(\gamma)$.
{Moreover, using the same notation as in the proof
of Theorem \ref{eureka}, again with $\gamma$ and $\rho$ instead of $\gamma_1$ and $\rho_1$, we find
\[
\begin{aligned}
\ol\mu_{1/4,\rho}(\matD_{\rho}(\gamma)) & \ge 
\ol\mu_{1/4,\rho}(\matD_{\rho}(\gamma) \cap  \mathfrak D^{(0)}(\g)) \\
& \ge \ol\mu_{1/4,\rho}(\mathfrak D^{(0)}(\g)) - \ol\mu_{1/4,\rho}(\mathfrak D^{(0)}(\g)) \setminus (\matD_{\rho}^{(0)}(\gamma)) \\
& \ge 1 - 2K_1 \g - C_*(\e) \g^{1/4} \rho \ge 1 - \frac{\g}{\g_*} ,
\end{aligned}
\]
where we have used the bounds \eqref{seriali} and \eqref{questaqui}.}
By Fubini's theorem, there exists a measurable set of potentials
{$\calG(\gamma)\subset\ol\matU_{1/4}(\ell^{N,\io}(\RRR))$  of measure less than $1- \sqrt{\g/\g_*}$ such that  for 
$V\in\calG(\gamma)$ the set $\matS_V \cap \ol\matU_{\rho}(\ell^{N,\io}(\RRR))$ is measurable and of
measure less than $1-\sqrt{\g/\g_*}$.}

By construction, any $W\in \matS_V$ gives rise to an almost-periodic solution,
{which, according to \eqref{lasoluzione!}, is given by {$u(x,t) =  \UU(x,\om(\ze(V,W,\e))t)$, with $\UU(x,\f)$ of the fom \eqref{noci}.}
By  \eqref{sposto}, for all $\f\in\TTT^{\ZZZ}$ one has
\begin{equation} \nonumber
%\| \mathfrak{i}_{\calU(\f)}\|_{s,\al} \le \| W\|_{s,\al} + O(\e)\| W\|_{s,\al}^5.
\| \calF( {\calU(\cdot,\f))} \|_{s,\al} \le \| W e^{\ii \f}  \|_{s,\al} +  4 C_1\e \rho^5 .
\end{equation}
Then,  we deduce the same estimate for $\calF(u(\cdot,t))$ for all $t\in\RRR$.
Thus the assertion follows provided that $\rho$ -- and hence $\| W\|_{s,\al}$ -- is small enough.}
\qed
%%%%%%%%%%%%%%%%%%%%%%%%%%%%%%%%%%%%%%%%%%%%%%%%%%%%%%%%%%%%%%%%%%%%%%%%%% 

%%%%%%%%%%%%%%%%%%%%%%%%%%%%%%%%%%%%%%%%%%%%%%%%%%%%%%%%%%%%%%%%%%%%%%%%%% 
%%%%%%%%%%%%%%%%%%%%%%%%%%%%%%%%%%%%%%%%%%%%%%%%%%%%%%%%%%%%%%%%%%%%%%%%%% 
\zerarcounters
\section{Invariant tori} \label{outro}
%%%%%%%%%%%%%%%%%%%%%%%%%%%%%%%%%%%%%%%%%%%%%%%%%%%%%%%%%%%%%%%%%%%%%%%%%% 
%%%%%%%%%%%%%%%%%%%%%%%%%%%%%%%%%%%%%%%%%%%%%%%%%%%%%%%%%%%%%%%%%%%%%%%%%% 

For any fixed $\ze\in \matW_N$ and $\e\in(-\e_*,\e_*)$,
the Lipschitz extensions $c\mapsto A(c,\ze,\e)$ and $c \mapsto R(c,\ze,\e)$ appearing in Propositions \ref{piergiorgio} and \ref{nuovaprop}
are constant on any set $\matT_I$ (see \eqref{TI}),
%\[\matT_I:= \{ c \in \mathtt g(s,\al) : {|c|^2= I} \} ,\] 
{with $I=|c|^2\in\RRR_+^\ZZZ \cap {\matU}_{1/4}(\mathtt g(2s,\al))$.}
If we consider only the last two equations in \eqref{system}
and  proceed as in the proof of Lemma 10.11 of \cite{CGP}, 
then, for all
${(V,c)}
\in
 \ol{\matU}_{\!1/4}(\ell^{N,\io}(\RRR))\times \ol{\matU}_{\!\rho_0}(\mathtt{g}(s,\al))
$,
%\red{with $\rho_0$ as in \eqref{V1},}
we obtain a solution $\ze=\widehat\ze(V,c,\e)=(\widehat\kappa(V,c,\e),\widehat\xi(V,c,\e))$ to \eqref{systemb} and \eqref{systemc}
which is still constant on any set $\matT_I$, i.e.~we have
\[
\begin{aligned}
\widehat\kappa(V,c,\e) + A(c,\widehat\kappa(V,c,\e),\widehat\xi(V,c,\e),\e)  & = 0 , \\
\widehat\xi(V,c,\e) + R(c,\widehat\kappa(V,c,\e),\widehat\xi(V,c,\e),\e) & = V ,
\end{aligned}
\]
and $\widehat\ze(V,c,\e)=\widehat\ze(ce^{\ii \theta},V,\e)$ for any $\theta\in\TTT^\ZZZ$.
Moreover,
if we set
\[
\Lambda_N(\g) := \Bigl\{ (V,c) \in \ol{\matU}_{\!1/4}(\ell^{N,\io}(\RRR)) \times \ol{\matU}_{\!\rho_0}(\mathtt{g}(s,\al))  :
\widehat\ze(V,c,\e) \in \matK_N(\g) \Bigr\}
\]
then for all $(V,c)\in \Lambda_N(\g)$,
since the compatibility condition \eqref{compatibility} is satisfied, the function
\begin{equation}\label{glidounnome}
\UU(x,\om(\widehat\ze(V,c,\e)) t) :=\UUU(x,\om(\widehat\ze(V,c,\e))t;c,\om(\widehat\ze(V,c,\e)),\e)
\end{equation}
is a solution to \eqref{NLS2} of the form \eqref{unnumero}, as highlighted in Remark \ref{Ueta}.

%%%%%%%%%%%%%%%%%%%%%%%%%%%%%%%%%%%%%%%%%%%%%%%%%%%%%%%%%%%%%%%%%%%%%%%%%% 
\subsection{Bi-Lipschitz map at fixed potential}\vspace{-.2cm}
%%%%%%%%%%%%%%%%%%%%%%%%%%%%%%%%%%%%%%%%%%%%%%%%%%%%%%%%%%%%%%%%%%%%%%%%%% 

For any fixed $V\in\ol{\matU}_{\!1/4}(\ell^{N,\io}(\RRR))$ and $\e$ small enough,  consider the map
\begin{equation}\label{WWW}
\gotWW(c) = \gotWW(c;V,\e) :=  c + \gotU(c,\widehat{\ze}(V,c,\e),\e) ,
\end{equation}
where $ \gotU(c,\ze,\e)$ is the Lipschitz extension of the function $ \gotu(c,\ze,\e)$ in \eqref{ugot},
according to Proposition \ref{piergiorgio}, so that
\[
\gotWW(c)=\{U_j(0,c;\om(\widehat\ze(c,V,\e),\e)\}_{j\in\ZZZ} \qquad \forall \; (V,c)\in \Lambda_N(\g) .
\]
As a consequence of the bound \eqref{dolore!exta}, the map
$\gotWW\!\!:\ol{\matU}_{\!\rho_0}(\mathtt{g}(s,\al)) \to \mathtt{g}(s,\al)$
is a lipeomorphism from the space of linear solutions to the space of initial data.
In agreement with \eqref{WWW}, if $\gotWW^{-1}\!\!:\gotWW(\ol{\matU}_{\!\rho_0}(\mathtt{g}(s,\al)) \to \ol{\matU}_{\!\rho_0}(\mathtt{g}(s,\al))$
 denotes its inverse, we write $\gotWW^{-1}(W)=\gotWW^{-1}(W;V,\e)$.}

%As a direct consequence of the bound \eqref{dolore!exta}, if we set
%%
%\begin{equation}\label{mappalli}
%\gotW(V,c,\e)  : =  c + \gotU(c,\widehat{\ze}(V,c,\e),\e) ,
%\end{equation}
%%
%where $ \gotU(c,\ze,\e)$ is the Lipschitz extension of the function $ \gotu(c,\ze,\e)$ in \eqref{ugot},
%according to Proposition \ref{piergiorgio},
%then, for any fixed $V\in\ol{\matU}_{\!1/4}(\ell^{N,\io}(\RRR))$ and $\e$ small enough, the map
%{$\gotWW \!:\ol{\matU}_{\!\rho_0}(\mathtt{g}(s,\al)) \to \mathtt{g}(s,\al)$, obtained by setting, at fixed values of $V$ and $\e$,
%%
%\begin{equation} \label{WWW}
%\gotWW(c) = \gotWW(c;V,\e) = \WWW(V,c,\e) , 
%\end{equation}
%%
%is a lipeomorphism from the space of  linear solutions to the one of initial data.}
%For future convenience let us introduce the inverse
%
%\begin{align} \label{W-1}
%\gotWW^{-1} \!:
%\ol{\matU}_{\!1/2}(\mathtt{g}(s,\al))
%& \to \ol{\matU}_{\!\rho_0}(\mathtt{g}(s,\al)) \\
%W & \mapsto \WWW^{-1}(V,W,\e) , \nonumber 
%\end{align}

%For $(V,c)\in \Lambda_N(\g)$, the map \eqref{WWW} establishes 
%\blue{ bijection between the set of the parameters $c$
%and the set of the initial data $W\in\TT_V$ of the solutions to the NLS equation \eqref{NLS2} of the form \eqref{unnumero};}

%%%%%%%%%%%%%%%%%%%%%%%%%%%%%%%%%%%%%%%%%%%%%%%%%%%%%%%%%%%%%%%%%%%%%%%%%% 
\begin{rmk} \label{2/3b}
\emph{
By construction, one has $(V,c)\in \Lambda_N(\g)$ if and only if $(V,W)\in \Gamma_N(\g)$, with 
$\Gamma_N(\g)$ defined in \eqref{paiodi}.
For any fixed $V\in\ol{\matU}_{\!1/4}(\ell^{N,\io}(\RRR))$ and $\e\in(-\e_*,\e_*)$,
consider any $c\in \ol{\matU}_{\!\rho_0}(\mathtt{g}(s,\al))$ such that
$(V,c)\in \Lambda_N(\g)$:
%\ol{\matU}_{\!1/2}(\mathtt{g}(s,\al)) \times \ell^{N,\io}(\RRR)$} such that $\widehat\ze(V,c,\e) \in \DgN $, 
by \eqref{chinottoneri} one has
\[
\begin{aligned}
\gotWW(c e^{\ii\theta}) & = \{ \UUU_j(0;c e^{\ii\theta} ,\om(\widehat\ze(c e^{\ii\theta} ,V,\e),\e)\}_{j\in\ZZZ} =
\{ \UUU_j(0;c e^{\ii\theta} ,\om(\widehat\ze(c ,V,\e),\e)\}_{j\in\ZZZ} \\
& =
\{ \UUU_j(\theta;c,\om(\widehat\ze(V,c,\e),\e)\}_{j\in\ZZZ} = \{ \UU_j(\theta) \}_{j\in\ZZZ} = \mathfrak{i}_{\,\UU}(\theta) ,
\end{aligned}
\]
where the notation \eqref{iU} has been used, and hence the set
\[
\MM(c)=\MM(V,c,\e):= %\WWW( \matT_{|c|^2},V,\e) := 
\bigcup_{\theta\in\TTT^\ZZZ} 
{\gotWW(ce^{\ii\theta})} = \mathfrak{i}_{\,\UU}(\TTT^\ZZZ)
\]
is invariant for \eqref{NLS4} and
the dynamics is conjugated to the translation $\theta \mapsto \theta + \om(\widehat\ze(V,c,\e)) t$.
In particular, if $|c_j|=|c'_j|$ for all $j\in\ZZZ$, then  $\MM(c) = \MM(c')$.
}
\end{rmk}
%%%%%%%%%%%%%%%%%%%%%%%%%%%%%%%%%%%%%%%%%%%%%%%%%%%%%%%%%%%%%%%%%%%%%%%%%%
 
For {
$c\in\ol{\matU}_{\!\rho_0}(\mathtt{g}(s,\al))$},
set $\DD_c:=\{j\in\ZZZ : c_j\ne 0\}$. Moreover, define the maps
$\psi\!:\TTT^{\DD_c} \to \TTT^\ZZZ$ and
$\CCCCC_c \!: \TTT^{\DD_c} \to \ol{\matU}_{\!\rho_0}(\mathtt{g}(s,\al))$
by setting
\begin{subequations}
\begin{align}
%\TTT^{\DD_c}\ni \theta \mapsto \phi (\theta) \in \TTT^\ZZZ,\qquad
\psi_j(\f) & = 
\begin{cases}
\f_j, & \quad j\in \DD_c , \\
0, & \quad j \in \ZZZ \setminus \DD_c, 
\end{cases}
\label{ovetto} \\
\CCCCC_c(\f) & = c \, e^{\ii\psi(\f)} . \phantom{\sum^n}
\label{finpiu}
\end{align}
\end{subequations}
%

%%%%%%%%%%%%%%%%%%%%%%%%%%%%%%%%%%%%%%%%%%%%%%%%%%%%%%%%%%%%%%%%%%%%%%%%%%
\begin{lemma}\label{mottooscrivo}
Fix $V\in\ol{\matU}_{\!1/4}(\ell^{N,\io}(\RRR))$ and $\e\in(-\e_*,\e_*)$. For any $c\in \ol{\matU}_{\!\rho_0}(\mathtt{g}(s,\al))$, the map
\begin{equation}\label{trentenni}
{\TTT^{\DD_c} \ni \f \mapsto \gotWW(\CCCCC_c(\f)) }
\end{equation}
is injective.
\end{lemma}
%%%%%%%%%%%%%%%%%%%%%%%%%%%%%%%%%%%%%%%%%%%%%%%%%%%%%%%%%%%%%%%%%%%%%%%%%%

%%%%%%%%%%%%%%%%%%%%%%%%%%%%%%%%%%%%%%%%%%%%%%%%%%%%%%%%%%%%%%%%%%%%%%%%%%
\prova
The map in \eqref{trentenni} is the composition of the two maps ${\gotWW \!: \ol{\matU}_{\!\rho_0}(\mathtt{g}(s,\al)) \to \gsa}$
and $\CCCCC_c \!: \TTT^{\DD_c} \to \ol{\matU}_{\!\rho_0}(\mathtt{g}(s,\al))$,
both of which are injective.
\qed
%%%%%%%%%%%%%%%%%%%%%%%%%%%%%%%%%%%%%%%%%%%%%%%%%%%%%%%%%%%%%%%%%%%%%%%%%%

%%%%%%%%%%%%%%%%%%%%%%%%%%%%%%%%%%%%%%%%%%%%%%%%%%%%%%%%%%%%%%%%%%%%%%%%%%
\begin{rmk} \label{altroimk}
\emph{
Since the map {$\gotWW$ in \eqref{WWW}} is a lipeomorphism,
the topological properties of the invariant set $\MM(c)$ depend on the properties of the map $\CCCCC_c$.
}
\end{rmk}
%%%%%%%%%%%%%%%%%%%%%%%%%%%%%%%%%%%%%%%%%%%%%%%%%%%%%%%%%%%%%%%%%%%%%%%%%%

%%%%%%%%%%%%%%%%%%%%%%%%%%%%%%%%%%%%%%%%%%%%%%%%%%%%%%%%%%%%%%%%%%%%%%%%%%
\begin{prop} \label{giausato}
Let $c\in
{
\ol{\matU}_{\!\rho_0}(\mathtt{g}(s,\al))}
$ be such that
\[
\delta(c):=\inf_{j\in \DD_c} |c_j|e^{s\jap{j}^\al} >0.
\]
Then, the map $\CCCCC_c\!:\TTT^{\DD_c}\to \matT_{|c|^2}$ is a diffeomorphism and the set $\matT_{|c|^2}$ is a
differentiable manifold embedded in $\mathtt{g}(s,\al)$.
\end{prop}
%%%%%%%%%%%%%%%%%%%%%%%%%%%%%%%%%%%%%%%%%%%%%%%%%%%%%%%%%%%%%%%%%%%%%%%%%%

%%%%%%%%%%%%%%%%%%%%%%%%%%%%%%%%%%%%%%%%%%%%%%%%%%%%%%%%%%%%%%%%%%%%%%%%%%
\prova
Set
\[
\mathtt{g}(s,\al;c):=\Bigl\{\tx =\{\tx _j\}_{j\in \ZZZ\setminus \DD_c}\in \ell^\io(\ZZZ\setminus \DD_c,\CCC) 
:  \|\tx \|_{s,\alpha,c}:= \sup_{j\in \ZZZ\setminus \DD_c} |\tx_j| e^{s\jap{j}^\alpha} <\io\Bigr\} ,
\]
and define
\[
\begin{aligned}
\matC_c&:= [1/2, +\io)^{\DD_c}\ \times \TTT^{\DD_c}\times \mathtt{g}(s,\al;c) \,,\\
 \matA_c&:= \bigl\{ u\in \gsa: |u_j|e^{s\jap{j}^\al}\ge \delta(c)/2 \quad \forall j\in \DD_c \bigr\} .
\end{aligned}
\]
Define also the map
\[
\Psi:\matC_c\to \matA_c,\qquad 
(r,\f,z)\mapsto \Psi(r,\f,z):= 
\begin{cases}
r_j c_j e^{\ii \f_j}, \qquad &j\in \DD_c, \\
z_j,  \qquad  &\mbox{otherwise}.
\end{cases}
\]
Note that $\matT_{|c|^2}=\Psi(\uno_c,\TTT^{\DD_c},0)$ with
{
$\uno_c:=\{1\}_{j\in \DD_c}$,}
so we are left to prove that $\Psi$ is a diffeomorphism. 
First of all, we note that $\Psi$ is a homeomorphism: in the case $\DD_c=\ZZZ$, this is known \cite[Appendix B]{BMP1},
and
{
the case $\DD_c\subsetneq \ZZZ$}
can be dealt with in a completely analogous way. On the other hand,
using that $\TTT^{\DD_c}$ is a Banach manifold whose tangent space is isomorphic to $\ell^\io(\DD_c,\RRR)$ \cite{russi},
%setting
%
%\[
%\blue{
%\ell^{\io}_{\DD_c}(\RRR):= \Bigl\{\tx =\{\tx _j\}_{j\in  D_c}\in \ell^\io(\RRR) :  \|\tx \|_{\io}:= \sup_{j\in \DD_c} |\tx_j|  <\io\Bigr\} ,}
%\]
%
we find, by direct inspection, that the Jacobian
\[
D\Psi : \ell^{\io}(\DD_c,\RRR)\times \ell^{\io}(\DD_c,\RRR)\times  \mathtt{g}(s,\al;c)\to \mathtt{g}(s,\al)
\]
is invertible. Hence $\Psi$ is a diffeomorphism, so the assertion follows.
\qed
%%%%%%%%%%%%%%%%%%%%%%%%%%%%%%%%%%%%%%%%%%%%%%%%%%%%%%%%%%%%%%%%%%%%%%%%%%

%%%%%%%%%%%%%%%%%%%%%%%%%%%%%%%%%%%%%%%%%%%%%%%%%%%%%%%%%%%%%%%%%%%%%%%%%%
\begin{rmk} \label{crucial}
\emph{
In the argument used in the proof of Proposition \ref{giausato}, it is crucial to consider spaces
{modelled on $\ell^\io$,}
so that, once a uniform lower bound is given, one may work
componentwise. On the other hand, if $\delta(c)=0$ the map $\CCCCC_c$ is still continuous and injective,
but the induced topology is finer than the topology inherited from that of $\gsa$.
}
\end{rmk}
%%%%%%%%%%%%%%%%%%%%%%%%%%%%%%%%%%%%%%%%%%%%%%%%%%%%%%%%%%%%%%%%%%%%%%%%%%

%%%%%%%%%%%%%%%%%%%%%%%%%%%%%%%%%%%%%%%%%%%%%%%%%%%%%%%%%%%%%%%%%%%%%%%%%% 
\subsection{Embedded tori and Cantor foliations}\vspace{-.2cm}
%%%%%%%%%%%%%%%%%%%%%%%%%%%%%%%%%%%%%%%%%%%%%%%%%%%%%%%%%%%%%%%%%%%%%%%%%% 

As a consequence of Proposition \ref{giausato}, for any $c$ such that $\delta(c)>0$ the invariant set $\MM(c)$
is also an embedded torus. In particular, if $\DD_c=\ZZZ$ we obtain a maximal torus. Conversely, if $\DD_c$ is finite,
then the condition $\delta(c)>0$ is trivially satisfied and we obtain a finite-dimensional embedded torus.

We are now ready to provide a proof of Theorem \ref{biascoappeso}, which follows from the discussion above.

\vspace{.3cm}

%%%%%%%%%%%%%%%%%%%%%%%%%%%%%%%%%%%%%%%%%%%%%%%%%%%%%%%%%%%%%%%%%%%%%%%%%%
\noindent\emph{Proof of Theorem \ref{biascoappeso}}.
The existence of the a solution to \eqref{nls} of the form \eqref{unnumero},
such that the map $\mathfrak{i}_{\,\UU}\!:\TTT^\ZZZ_a\to\gsa$ is analytic,
follows from item 5 of Theorem \ref{moser}, combined with \eqref{glidounnome}.
The relation between the constant $a=a(\e,s,\al,N)$ appearing in Theorem \ref{biascoappeso}
and the constant $a=\aaa(\e,s,\al,\om)$ in item 5 of Theorem \ref{moser} is established by setting
\[
a(\e,s,\al,N) = \min_{(V,W) \in \Gamma_N(\g)} \bigl\{ \aaa(\e,s,\al,\om(\ze(V,W,\e))) \bigr\} ,
\]
where the minimum is easily estimated by taking into account that all bounds are uniform for $\ze(V,W,\e) \in \matK_N(\g)$.

Regarding the properties of the solution listed in the statement, we reason as follows.
Fix $\de>0$ and take $(V,W)\in\Gamma_N(\g)$ such that  {$W\in \matA_{V,\de}$}. Since
$\gotWW^{-1}(W;V,\e) - W =\Delta_2(v,W)$ is $O(\e)$, because of the definition \eqref{WWW} and the bound \eqref{dolore!exta},
then, if we take $|\e|<\e_\star$, for some $\e_\star \in(0,\e_*]$ depending on $\de$, and set $c:=\gotWW^{-1}(W;V,\e)$,
we find that $\DD_c=\ZZZ$, $\delta(c)\ge\de/2$ and, by definition, $(V,c)\in\Lambda_N(\g)$.
%$c\in \ol{\matU}_{1/2}(\gsa)$, 
But then the assertion follows combining Remark \ref{2/3b} and Proposition \ref{giausato}.
\qed
%%%%%%%%%%%%%%%%%%%%%%%%%%%%%%%%%%%%%%%%%%%%%%%%%%%%%%%%%%%%%%%%%%%%%%%%%%

%%%%%%%%%%%%%%%%%%%%%%%%%%%%%%%%%%%%%%%%%%%%%%%%%%%%%%%%%%%%%%%%%%%%%%%%%%
\begin{rmk} \label{crucial2}
\emph{
Using the notation in \eqref{c2}, 
for any set $\IIII \subset \RRR_+^\ZZZ \cap \ol{\matU}_{1/4} (\mathtt g(2s,\al ))$, let $\matT_{\IIII}$ be defined according to \eqref{TC}, i.e.
\[
\matT_\IIII = \Bigl\{ c \in \gsa \,:\, c = \sqrt{I} e^{\ii \theta} , \, I \in \IIII, \, \theta\in\TTT^\ZZZ \Bigr\} .
\]
Fix $V\in\mathcal{G}$, $\de>0$ and $\e\in(-\e_\star,\e_\star)$ as in Theorem \ref{biascoappeso}. Define
\[
\begin{aligned}
\IIII_0  & =\IIII_0 (V,\de):= \Bigl\{ I \in \RRR_+^\ZZZ \cap \ol{\matU}_{1/4} (\mathtt g(2s,\al )) : \sqrt{I} \in \calA_{\de}(\gsa) \Bigr\} , \\
\IIII  & =\IIII (V,\de,\e):= \Bigl\{ I \in \RRR_+^\ZZZ \cap \ol{\matU}_{1/4} (\mathtt g(2s,\al )) : \gotWW(\sqrt{I}) \in \matA_{V,\de} \Bigr\} ,
\end{aligned}
\]
and, similarly,
\[
\begin{aligned}
\IIIIO  & =\IIIIO(V,\de):= \Bigl\{ c \in \ol{\matU}_{1/2} (\mathtt g(s,\al )) :  |c|^2 \in \IIII_0 \Bigr\} , \\
\IIIIC & = \IIIIC(V,\de,\e) := \Bigl\{ c \in \ol{\matU}_{1/2} (\mathtt g(s,\al)) : |c|^2 \in \IIII \Bigr\} .
\end{aligned}
\]
Then, for all $I\in\IIII$, the set $\gotWW(\matT_I)$ is an invariant submanifold homeomorphic to the flat torus.
The set
\[
\gotWW(\matT_{\IIII}) = \bigcup_{I\in \IIII} \MM(\sqrt{I}) = \bigcup_{c \in \IIIIC} \MM(c) 
\]
provides a Cantor foliation of $\calA(\mathtt g(s,\al),\de/2)$ into invariant tori. If we define, instead,
\[
\IIII_1 := \IIII_1 (V,\e) = \Bigl\{ I \in \RRR_+^\ZZZ \cap \ol{\matU}_{1/4} (\mathtt g(2s,\al )) : \gotWW(\sqrt{I}) \in \TT_{V} \Bigr\} ,
\]
then, for any $I\in\IIII_1$, the set $\gotWW(\matT_I)$ is still an invariant torus. However, it may fail to be a submanifold.
In particular the set
\[
\gotWW(\matT_{\IIII_1}) = \bigcup_{I\in \IIII_1} \MM(\sqrt{I}) 
\]
is no more than a stratification in invariant tori, as well as \eqref{parafoli} when $\e=0$.
}
\end{rmk}
%%%%%%%%%%%%%%%%%%%%%%%%%%%%%%%%%%%%%%%%%%%%%%%%%%%%%%%%%%%%%%%%%%%%%%%%%%

%%%%%%%%%%%%%%%%%%%%%%%%%%%%%%%%%%%%%%%%%%%%%%%%%%%%%%%%%%%%%%%%%%%%%%%%%%
%%%%%%%%%%%%%%%%%%%%%%%%%%%%%%%%%%%%%%%%%%%%%%%%%%%%%%%%%%%%%%%%%%%%%%%%%%

\end{document}